\documentclass[11pt]{article}
\usepackage{amsmath,amsthm,amssymb,color,epic,eepic,
cite}
\renewcommand{\baselinestretch}{1.2}

\def\baselinestretch{1.4}
\newlength{\minitwocolumn}
\newcommand{\Z}{{\Bbb Z}} 
\newcommand{\R}{{\Bbb R}} 
\newcommand{\C}{{\Bbb C}} 
\newcommand{\FF}{{\Bbb F}} 

\newcommand{\F}{{\cal F}}

\newcommand{\cD}{{\cal D}}
\newcommand{\cA}{{\cal A}}

\newcommand{\B}{{\cal B}}
\newcommand{\cI}{{\cal I}}
\newcommand{\cH}{{\cal H}}
\newcommand{\cN}{{\cal N}}
\newcommand{\cR}{{\cal R}}
\newcommand{\cP}{{\cal P}}
\newcommand{\cE}{{\cal E}}
\newcommand{\cM}{{\cal M}}
\newcommand{\cQ}{{\cal Q}}

\newcommand{\hd}{\widehat{d}}
\newcommand{\hL}{\widehat{L}}
\newcommand{\hcL}{\widehat{\cL}}

\newcommand{\tR}{
\cR}
\newcommand{\tL}{
\cL}

\newcommand{\cL}{{\cal L}}
\newcommand{\cK}{{\cal K}}

\newcommand{\la}{\lambda}

\newcommand{\al}{\alpha}

\newcommand{\ep}{\epsilon}
\newcommand{\vep}{\varepsilon}

\newcommand{\bep}{\bar{\epsilon}}

\newcommand{\ha}{{\alpha}}
\newcommand{\hb}{{\beta}}

\newcommand{\hf}{\widehat{f}}

\newcommand{\hV}{\widehat{V}}

\newcommand{\bh}{{\bar{\h}}}

\newcommand{\hrho}{
{\rho}}
\newcommand{\hj}{\widehat{j}}
\newcommand{\hi}{\widehat{i}}
\newcommand{\hl}{\widehat{l}}
\newcommand{\hk}{\widehat{k}}

\newcommand{\nn}{{\nonumber}}
\newcommand{\bea}{\begin{eqnarray}}
\newcommand{\ena}{\end{eqnarray}}
\newcommand{\beit}{\begin{itemize}}
\newcommand{\enit}{\end{itemize}}

\newcommand{\be}{\begin{eqnarray*}}
\newcommand{\en}{\end{eqnarray*}}
\newcommand{\lb}[1]{\label{#1}}

\newcommand{\ds}[1]{{\displaystyle #1 }}

\newcommand{\End}{{\rm End}}

\newcommand{\id}{{\rm id}}

\newcommand{\td}{{\widetilde{d}}}

\def\infq4p#1{{(#1;q^4,p)_\infty}}

\newcommand{\hPsi}{\widehat{\Psi}}

\newcommand{\tot}{\widetilde{\otimes}}

\newcommand{\mmatrix}[1]{\begin{matrix} #1 \end{matrix}}
\newcommand{\mat}[1]{\left(\mmatrix{#1}\right)}

\font\teneufm=eufm10
\font\seveneufm=eufm7
\font\fiveeufm=eufm5
\newfam\eufmfam
\textfont\eufmfam=\teneufm
\scriptfont\eufmfam=\seveneufm
\scriptscriptfont\eufmfam=\fiveeufm

\let\goth\frak
\newcommand{\slth}{\widehat{\goth{sl}}_2}

\newcommand{\slnh}{\widehat{\goth{sl}}_N}

\newcommand{\g}{\goth{g}}

\newcommand{\Bqla}{{{\cal B}_{q,\lambda}}}

\newcommand{\UqpBN}{U_{q,p}(B_N^{(1)})}

\newcommand{\Bnh}{B_N^{(1)}}

\newcommand{\gl}{{\goth{gl}}}

\newcommand{\glth}{\widehat{\goth{gl}}_2}

\newcommand{\h}{\goth{h}}
\newcommand{\gh}{\widehat{\goth{g}}}
\newcommand{\hh}{\goth{h}}

\newcommand{\ghbig}{\widehat{\mbox{\fourteeneufm g}}}  

\font\fourteeneufm=eufm10 scaled\magstep2    
   
\makeatletter
\@addtoreset{equation}{section}
\makeatother

\newtheorem{thm}{Theorem}[section]
\newtheorem{prop}[thm]{Proposition}
\newtheorem{lem}[thm]{Lemma}
\newtheorem{cor}[thm]{Corollary}
\newtheorem{conj}[thm]{Conjecture}

\newtheorem{dfn}[thm]{Definition}

\begin{document}


\vspace{-1cm}
\begin{center}
{\bf\Large  Elliptic Quantum Group $U_{q,p}(B_N^{(1)})$ and Vertex Operators
\\[7mm] }
{\large  Hitoshi Konno${}^{\dagger}$ and Kazuyuki Oshima${}^{\star}$ }\\[6mm]
${}^\dagger${\it  Department of Mathematics, Tokyo University of Marine Science and 
Technology, \\Etchujima, Koto, Tokyo 135-8533, Japan\\
       hkonno0@kaiyodai.ac.jp}\\
${}^\star${\it Center for General Education, Aichi Institute of Technology, \\ Yakusa-cho, Toyota
470-0392, Japan }\\
oshima@aitech.ac.jp
\\[7mm]
\end{center}

\begin{center}
Dedicated to Professor Akihiro Tsuchiya on his 70th birthday. 
\end{center}

\begin{abstract}
\noindent 
Assuming the existence of the $L$-operators, 
we study the Hopf algebroid structure of $U_{q,p}(B_N^{(1)})$. 
As an application, we derive the type I and II vertex operators, which intertwine the 
$U_{q,p}(B_N^{(1)})$-modules of generic level,   
by assuming some analytic properties of the $L$-operators. 
For the level-1 case, we construct their free field realizations  
and show that the results satisfy the desired commutation relations with coefficients 
given by the elliptic dynamical $R$-matrices of the $\Bnh$ type.   
\end{abstract}
\nopagebreak

\section{Introduction}
The algebra $U_{q,p}(\gh)$\cite{K98,JKOS99,FKO} is an elliptic analogue of the quantum affine algebra $U_q(\gh)$ in the Drinfeld realization\cite{Drinfeld} associated with the affine Lie algebra $\gh$. The $U_{q,p}(\gh)$ is expected to give a realization of the face type elliptic quantum group \cite{Fronsdal, JKOStransfG} equipped with the Hopf algebroid structure.
In the previous works\cite{JKOS99, KojimaKonno, Konno09}, we have  constructed the $L$-operator of  $U_{q,p}(A_N^{(1)})$ in terms of  the elliptic currents, the generating functions of the Drinfeld generators of $U_{q,p}(A_N^{(1)})$. 
The $L$-operator  satisfies the $RLL$-relation with the elliptic dynamical $R$-matrix of the $A_N^{(1)}$ type\cite{JMO} and allows us to define the Hopf algebroid structure 
to $U_{q,p}(A_N^{(1)})$. 

The elliptic quantum group $U_{q,p}(A_N^{(1)})$ equipped with the Hopf algebroid structure  has proved to be quite useful in construction of both the finite and infinite dimensional representations as well as 
 their  intertwining operators, i.e.  the  vertex operators,  
  in terms of the free fields. Such construction becomes a central tool in the algebraic analysis 
of the face type solvable lattice models associated with the vector representation of $\gh$\cite{JMO} in the spirit of Jimbo and Miwa\cite{JM}. See for example \cite{LP,LashPug,KKW}.

The purpose of this paper is to continue the above study to the case $U_{q,p}(\Bnh)$.  
The $U_{q,p}(\Bnh)$ itself has an interesting connection to the deformation of Fateev-Lukyanov's\cite{FaLu} $W\! B_N$ algebra\cite{FKO}. Assuming the existence of the $L$-operators in the elliptic algebra $U_{q,p}(B_N^{(1)})$, 
we give an Hopf algebroid structure of $U_{q,p}(B_N^{(1)})$. 
We then define the type I and II vertex operators as the  intertwining operators of the   
$U_{q,p}(B_N^{(1)})$-modules of generic level.    
By assuming some analytic properties of the $L$-operators and the half currents, which are expected to be defined recursively through the
 Gauss decomposition of the $L$-operators,  we show that the components of the 
 both types of vertex operators are constructed by applying certain half currents to the top component.  
For the level-$1$ case, we construct their free field realizations  
and show that the results satisfy the desired commutation relations with coefficients 
given by the elliptic dynamical $R$-matrices of the $\Bnh$ type. These results give  elliptic  and dynamical analogues of those obtained for $U_q(B_N^{(1)})$ in \cite{JinMis,JinZho}. 

This paper is organized as follows. In Sec.2, we define the elliptic algebra $U_{q,p}(B_N^{(1)})$ as a certain topological algebra. In particular, we introduce the orthonormal basis type elliptic bosons and define the elliptic currents $k_{\pm j}(z)$. 
The Sec.3 is devoted to a conjecture of the construction of the $L$-operators in terms of the half currents of $U_{q,p}(B_N^{(1)})$. In Sec.4, assuming the existence of the $L$-operators we introduce the $H$-Hopf algebroid structure to  $U_{q,p}(B_N^{(1)})$. 
In Sec.5, we give the vector representation and the level-$1$ highest weight representation of  $U_{q,p}(B_N^{(1)})$. In Sec.6, after giving a construction of the type I and II vertex operators at generic level, we present a free field realization of the level 1 vertex operators and show that they satisfy the desired commutation relations with the coefficients given by the elliptic dynamical $R$-matrices. 
In Appendix A, we summarize a connection of $U_{q,p}(B_N^{(1)})$ to the quasi-Hopf formulation $\Bqla(B_N^{(1)})$ of the elliptic quantum group.  In Appendix B, we give a list of conjectural expressions for the half currents  of $U_{q,p}(B_N^{(1)})$. 


\section{Elliptic Algebra $U_{q,p}(B_N^{(1)})$}\lb{secuqpBn}
In this section,  we give a definition of the elliptic algebra $U_{q,p}(B_N^{(1)})$ 
associated with the affine Lie algebra $B^{(1)}_N$. 

\subsection{Definition}
Let  $A=(a_{ij})\ i,j\in \{0\}\cup I,\ I=\{1,\cdots,N\}$ be   
 the $B^{(1)}_N$ type generalized Cartan matrix\cite{Kac}. We denote by $B=(b_{ij})$,
 $b_{ij}=d_i a_{ij}$ the symmetrization of $A$ with $d_0=\cdots =d_{N-1}=1, d_N=1/2$. 
Let $q=e^{\hbar}\in \C[[\hbar]]$ and set $q_i=q^{d_i}$. 
Let $p$ be an indeterminate. 
We use the following notations.
\be 
&&[n]_q=\frac{q^n-q^{-n}}{q-q^{-1}}, \quad [n]_i=\frac{q_i^n-q_i^{-n}}{q_i-q_i^{-1}}, \qquad [n)_i=\frac{q^n-q^{-n}}{q_i-q_i^{-1}},\\
&&[n]_i!=[n]_i[n-1]_i\cdots [1]_i,\quad \left[\mmatrix{m\cr n\cr}\right]_i=\frac{[m]_i!}{[n]_i![m-n]_i!},\\
&&(x;q)_\infty=\prod_{n=0}^\infty(1-x q^n),\quad (x;q,t)_\infty=\prod_{n,m=0}^\infty(1-x q^n t^m),\quad
\Theta_p(z)=(z;p)_{\infty}(p/z;p)_\infty(p;p)_\infty.
\en

Let $\h=\widetilde{\h}\oplus \C d$, $\widetilde{\h}=\bar{\h}\oplus\C c$, $\bar{\h}=\oplus_{i\in I}\C h_i$ be the Cartan subalgebra of $B^{(1)}_N$.  
 Define  $\delta, \Lambda_0, \al_i\ (i\in I) \in \h^*$ 
by 
\bea
&&<\al_i,h_j>=a_{ji}, \ <\delta,d>=1=<\Lambda_0,c>,\lb{pairinghhs}
\ena
the other pairings are 0. We also define $\bar{\Lambda}_i\ (i \in I) \in \h^*$ by 
\bea
&&<\bar{\Lambda}_i,h_j>=\delta_{i,j} \lb{pairingLambdah}
\ena
We set $\bh^*=\oplus_{i\in I} \C \bar{\Lambda}_i, $ $\widetilde{\h}^*=\bh^*\oplus \C \Lambda_0$, 
 $\cQ=\oplus_{i\in I}\Z \al_i$ and 
$\cP=\oplus_{i\in I}\Z \bar{\Lambda}_i$. 
 Let $\{ \ep_j\ (1\leq j\leq N)\}$ be an orthonormal basis  in $\R^N$ 
with the inner product
$( \ep_j, \ep_k )=\delta_{j,k}$. 
We realize the simple roots by $\alpha_j=\ep_j-\ep_{j+1}\ (1\leq j\leq N-1),\ \al_N=\ep_N$ and the fundamental weights by $\bar{\Lambda}_j=\ep_1+\cdots+\ep_j \ (1\leq j\leq N-1), \ \bar{\Lambda}_N=(\ep_1+\cdots+\ep_N)/2 $. 
We define $h_{\ep_j}\in \bar{\h} \ (j\in I)$ by 
$<\ep_i,h_{\ep_j}>=(\ep_i,{\ep_j})$ and $h_\alpha \in \bar{\h}$ for $\alpha=\sum_j c_j \ep_j, \ c_j\in \C$ by
 $h_\alpha=\sum_j c_j h_{\ep_j}$.
We regard $\bar{\h}\oplus \bar{\h}^*$ as the Heisenberg 
algebra  by
\bea
&&~[h_{{\epsilon}_j},{\epsilon}_k]
=( {\epsilon}_j,{\epsilon}_k ),\qquad [h_{{\epsilon}_j},h_{{\epsilon}_k}]=0=
[{\epsilon}_j,{\epsilon}_k].\lb{HA1}
\ena
In particular, we have $[h_{j}, \alpha_k]=a_{j k}$.  We also set $h^j=h_{\bar{\Lambda}_j}$. 
 
Let $\{P_{{\alpha}}, Q_{{\beta}}\}\ 
({\alpha}, {\beta} \in \bar{\h}^*)$ be the Heisenberg algebra 
defined by the commutation relations
\begin{eqnarray}
&&[P_{{\epsilon}_j}, Q_{{\epsilon}_k}]=
( {\epsilon}_j, {\epsilon}_k ), \qquad 
[P_{{\epsilon}_j}, P_{{\epsilon}_k}]=0=
[Q_{{\epsilon}_j}, Q_{{\epsilon}_k}],\lb{HA2}
\end{eqnarray}
where
$P_\alpha=\sum_j c_j P_{\ep_j}$ for $\alpha=\sum_j c_j \ep_j$.
We set $P_{\bh}=\oplus_{j\in I}\C P_{\ep_j}, Q_{\bh}=\oplus_{j\in I}\C Q_{\ep_j}$  $P_{j}=P_{\alpha_j^\vee}, P^j=P_{\bar{\Lambda}_j}$ and 
$Q_{j}=Q_{\alpha_j}, Q^j=Q_{\bar{\Lambda}_j^\vee}$. 

For the abelian group 
$\cR_Q= \sum_{j=1}^N\Z Q_{\al_j}$, 
we denote by $\C[\cR_Q]$ the group algebra over $\C$ of $\cR_Q$. 
We denote by $e^{\al}$ the element of $\C[\cR_Q]$ corresponding to $\al\in \cR_Q$. 
These $e^\al$ satisfy $e^\al e^\beta=e^{\al+\beta}$
 and $(e^\al)^{-1}=e^{-\al}$. 
In particular, $e^0=1$ is the identity element. 

Now let us consider to double the Cartan subalgebra : $H=\widetilde{\h}\oplus P_{\bh}=\sum_{j}\C(P_{\bep_j}+h_{\bep_j})+\sum_j\C P_{\bep_j}+\C c$. We  denote its dual space by $H^*=\widetilde{\h}^*\oplus Q_{\bh}$. We define the paring by \eqref{pairingLambdah},  
$<Q_\al,P_\beta>=(\al,\beta)$ and  $<Q_\al,h_\beta>=<Q_\al,c>=<Q_\al,d>=0=<\al,P_\beta>=<\delta,P_\beta>=<\Lambda_0,P_\beta>$ . We define $\FF=\cM_{H^*}$ to be the field of meromorphic functions on $H^*$. 
\begin{dfn}\cite{FKO}\lb{defUqp}
The elliptic algebra $U_{q,p}(B_N^{(1)})$ is a topological algebra over $\FF[[p]]$ generated by $\cM_{H^*}$,  $e_{j,m}, f_{j,m}, \al^\vee_{j,n}, K^\pm_j$,  
$(j\in I, m\in \Z, n\in \Z_{\not=0})$, ${\hd}$ and the central element $c$. 
We assume $K^\pm_{j}$ are invertible and set 
\be
&&e_j(z)=\sum_{m\in \Z} e_{j,m}z^{-m}
,\quad f_j(z)=\sum_{m\in \Z} f_{j,m}z^{-m},\lb{deffn}\\
&&{\psi}_j^+(
q^{-\frac{c}{2}}
z)=K^+_{j}\exp\left(-(q_j-q_j^{-1})\sum_{n>0}\frac{\al^\vee_{j,-n}}{1-p^n}z^n\right)
\exp\left((q_j-q_j^{-1})\sum_{n>0}\frac{p^n\al^\vee_{j,n}}{1-p^n}z^{-n}\right),\\
&&{\psi}_j^-(q^{\frac{c}{2}}z)=K^-_{j} \exp\left(-(q_j-q_j^{-1})\sum_{n>0}\frac{p^n\al^\vee_{j,-n}}{1-p^n}z^n\right)
\exp\left((q_j-q_j^{-1})\sum_{n>0}\frac{\al^\vee_{j,n}}{1-p^n}z^{-n}\right). 
\en
We call $e_j(z), f_j(z), \psi^\pm_j(z)$ the elliptic currents. 
The  defining relations are as follows. For $g(P), g(P+h)\in \cM_{H^*}$, 
\bea
&&g({P+h})e_j(z)=e_j(z)g({P+h}),\quad g({P})e_j(z)=e_j(z)g(P-<Q_{\al_j},P>),\lb{ge}\\
&&g({P+h})f_j(z)=f_j(z)g(P+h-<{\al_j},P+h>),\quad g({P})f_j(z)=f_j(z)g(P),\lb{gf}
\ena
\bea
&&[g(P), \al^\vee_{i,m}]=[g(P+h),\al^\vee_{i,n}]=0,\lb{gboson}\\ 
&&g({P})K^\pm_j=K^\pm_jg(P-<Q_{\al_j},P>),\\
&&\ g({P+h})K^\pm_j=K^\pm_jg(P+h-<Q_{\al_j},P>),\lb{gKpm}
\\
&&[\hd, g(P+h)]=[\hd, g(P)]=0,
\quad\lb{dg}\\
&& [\hd, \al^\vee_{j,n}]=n\al^\vee_{j,n},\quad [\hd, e_j(z)]=-z\frac{\partial}{\partial z}e_j(z), \quad  [\hd, f_j(z)]=-z\frac{\partial}{\partial z}f_j(z),\quad \lb{dedf}\\
&&K_i^{\pm}e_j(z)=q_i^{\mp a_{ij}}e_j(z)K_i^{\pm},\quad 
K_i^{\pm}f_j(z)=q_i^{\pm a_{ij}}f_j(z)K_i^{\pm},
\\
&&[\al^\vee_{i,m},\al^\vee_{j,n}]=\delta_{m+n,0}\frac{[a_{ij}m]_i
[cm)_j
}{m}
\frac{1-p^m}{1-p^{*m}}
q^{-cm}
,\lb{ellboson}\\
&&
[\al^\vee_{i,m},e_j(z)]=\frac{[a_{ij}m]_i}{m}\frac{1-p^m}{1-p^{*m}}
q^{-cm}z^m e_j(z),
\lb{bosonve}\\
&&
[\al^\vee_{i,m},f_j(z)]=-\frac{[a_{ij}m]_i}{m}z^m f_j(z)
,\lb{bosonvf}\\
&&
z_1 \frac{(q^{b_{ij}}z_2/z_1;p^*)_\infty}{(p^*q^{-b_{ij}}z_2/z_1;p^*)_\infty}e_i(z_1)e_j(z_2)=
-z_2 \frac{(q^{b_{ij}}z_1/z_2;p^*)_\infty}{(p^*q^{-b_{ij}}z_1/z_2;p^*)_\infty}e_j(z_2)e_i(z_1),\lb{ee}\\
&&
z_1 \frac{(q^{-b_{ij}}z_2/z_1;p)_\infty}{(pq^{b_{ij}}z_2/z_1;p)_\infty}f_i(z_1)f_j(z_2)=
-z_2 \frac{(q^{-b_{ij}}z_1/z_2;p)_\infty}{(pq^{b_{ij}}z_1/z_2;p)_\infty}f_j(z_2)f_i(z_1),\lb{ff}\\
&&[e_i(z_1),f_j(z_2)]=\frac{\delta_{i,j}}{q_i-q_i^{-1}}
\left(\delta(
q^{-c}
z_1/z_2)
\psi_j^-(
q^{\frac{c}{2}}
z_2)-
\delta(
q^c
z_1/z_2)
\psi_j^+(
q^{-\frac{c}{2}}
z_2)
\right),\lb{eifj}
\ena
\bea
&&\sum_{\sigma\in S_{a}}\prod_{1\leq m< k \leq a}
\frac{(p^*q^{2}{z_{\sigma(k)}}/{z_{\sigma(m)}}; p^*)_{\infty}}
{(p^*q^{-2}{z_{\sigma(k)}}/{z_{\sigma(m)}}; p^*)_{\infty}}\nn\\
&&\quad\times\sum_{s=0}^{a}(-1)^{s}
\left[\begin{array}{c}
a\cr
s\cr
\end{array}\right]_{i}\prod_{1\leq m \leq s}\frac{(p^*q^{b_{ij}}{w}/{z_{\sigma(m)}}; p^*)_{\infty}}
{(p^*q^{-b_{ij}}{w}/{z_{\sigma(m)}}; p^*)_{\infty}}
\prod_{s+1\leq m \leq a}\frac{(p^*q^{b_{ij}}{z_{\sigma(m)}}/{w}; p^*)_{\infty}}
{(p^*q^{-b_{ij}}{z_{\sigma(m)}}/{w}; p^*)_{\infty}}
\nn\\
&&\quad\times
e_{i}(z_{\sigma(1)})\cdots e_{i}(z_{\sigma(s)})e_{j}(w)e_{i}(z_{\sigma(s+1)}) \cdots 
e_{i}(z_{\sigma(a)})=0,\label{serree}
\ena
\bea
&&\sum_{\sigma\in S_{a}}\prod_{1\leq m< k \leq a
}
\frac{(pq^{-2}{z_{\sigma(k)}}/{z_{\sigma(m)}}; p)_{\infty}}
{(pq^{2}{z_{\sigma(k)}}/{z_{\sigma(m)}}; p)_{\infty}}\nn\\
&&\quad\times\sum_{s=0}^{a}(-1)^{s}\left[\begin{array}{c}
a\cr
s\cr
\end{array}\right]_{i}\prod_{1\leq  m \leq s}
\frac{(pq^{-b_{ij}}{w}/{z_{\sigma(m)}}; p)_{\infty}}{(pq^{b_{ij}}{w}/{z_{\sigma(m)}}; p)_{\infty}}
\prod_{s+1\leq m \leq a}\frac{(pq^{-b_{ij}}{z_{\sigma(m)}}/{w}; p)_{\infty}}{(pq^{b_{ij}}{z_{\sigma(m)}}/{w}; p)_{\infty}}\nn\\
&&\times f_{i}(z_{\sigma(1)})\cdots f_{i}(z_{\sigma(s)})f_{j}(w)f_{i}(z_{\sigma(s+1)}) \cdots f_{i}(z_{\sigma(a)})=0
\quad(i\neq j, a=1-a_{ij}),\label{serref}
\end{eqnarray}
where $p^*=p
q^{-2c}
$ and $\delta(z)=\sum_{n\in \Z}z^n$.  
We also denote by ${U}'_{q,p}(\gh)$ the subalgebra obtained by removing ${d}$. 
\end{dfn}
We treat the relations \eqref{dedf}, \eqref{bosonve}-\eqref{serref} as formal Laurent series in $z, w$ and $z_j$'s. 
All the coefficients in $z_j$'s are well-defined in the $p$-adic topology.  

\subsection{The orthonormal basis type elliptic bosons}
It is convenient to introduce the simple root type generators $\al_{j,m}$ and $\al'_{j,m}$ defined by $\al_{j,m}=[d_j]_q\al^\vee_{j,m}$ 
and $\ds{\alpha'_{j,m}=\frac{1-p^{*m}}{1-p^m}q^{cm} \alpha_{j,m},\ (j\in I, n\neq 0)}$. 
From \eqref{ellboson}, \eqref{bosonve}, \eqref{bosonvf},  we have 
\bea
&&[\al_{i,m},\al_{j,n}]=\frac{[b_{ij}m]_q[cm]_q}{m}\frac{1-p^m}{1-p^{*m}}q^{-km}\delta_{m+n,0},\lb{alal}\\
&&[\al'_{i,m},\al'_{j,n}]=\frac{[b_{ij}m]_q[cm]_q}{m}\frac{1-p^{*m}}{1-p^{m}}q^{km}\delta_{m+n,0},\\
&&[\al_{i,m},\al'_{j,n}]=\frac{[b_{ij}m]_q[cm]_q}{m}\delta_{m+n,0},\\
&&[\al_{i,m},e_j(z)]=\frac{[b_{ij}m]_q}{m}\frac{1-p^m}{1-p^{*m}}q^{-cm}z^m e_j(z),
\lb{bosone}\\
&&[\al'_{i,m},f_j(z)]=-\frac{[b_{ij}m]_q}{m}\frac{1-p^{*m}}{1-p^{m}}q^{cm}z^m f_j(z).\lb{bosonf}
\ena

Let $\eta=-(2N-1)/2$. Let us further introduce the orthonormal basis type elliptic bosons $\cE^{\pm j}_m \ (j\in \{0\}\cup I, m\in \Z_{\not=0})$ \cite{FKO} by
\bea
&&\cE^{\pm j}_m=q^{\pm jm}C_m\left(q^{\pm \eta  m}\sum_{k=1}^{j-1}[km]_q\al_{k,m}
\pm\sum_{k=j}^N [ ( \eta+k ) m]_+\al_{k,m} \right),\\
&&\cE^0_m=\frac{[\frac{m}{2}]_q}{[m]_q}(\cE^{+N}_m+\cE^{-N}_m).
\ena
Here we set
\be
&&C_m=\frac{[\eta m]_q}{[m]_q^2[2\eta m]_q}, \qquad 
[m]_+=\frac{q^m+q^{-m}}{q-q^{-1}}.
\en

\begin{prop}\lb{cEcE}
\bea
&&[\cE^{\pm j}_m,\cE^{\pm j}_n]
=\delta_{m+n,0}\frac{[cm]_q[\eta m]_q [2(\eta+1)m]_q}{m(q-q^{-1})^2[m]_q^3[2\eta m]_q [(\eta+1)m]_q}\frac{1-p^m}{1-p^{*m}}q^{-cm}, 
\\
&&[\cE^{\pm j}_m,\cE^{\mp j}_n]
=\mp\delta_{m+n,0}\frac{q^{\pm jm}[cm]_q[\eta m]_q}{m[m]_q^3(q-q^{-1})^2[2\eta m]_q}
\frac{1-p^m}{1-p^{*m}}q^{-cm}
\left(q^{\pm(\eta +j)m}[m]_q\pm q^{\mp(j-1)m}[\eta m]_{+}\right),\nn\\
\\
&&[\cE^{\pm j}_m,\cE^{\pm k}_n]
=\mp{\rm sgn}(k-j)\delta_{m+n,0}q^{\mp({\rm sgn}(k-j) \eta +k-j)m}
\frac{[cm]_q[\eta m]_q}{m(q-q^{-1})[m]_q^2[2 \eta m]_q}\frac{1-p^m}{1-p^{*m}}q^{-cm},\nn\\
&&\\
&&
[\cE^{\pm j}_m,\cE^{\mp k}_n]=\mp \delta_{m+n,0}q^{\pm(\eta+j+k)m}\frac{[cm]_q[\eta m]_q}{m(q-q^{-1})[m]_q^2[2 \eta m]_q}\frac{1-p^m}{1-p^{*m}}q^{-cm},
\ena
Here
\be
{\rm sgn}(l-j)=\left\{\mmatrix{+&(l>j),\cr
                              -&(l<j). \cr }\right.
\en
\end{prop}

Then one can realize the simple root type $\al_{j,m}$ in terms of the orthonormal basis type $\cE^{\pm j}_m$ a follows.
\begin{prop}\lb{alphaAA}
\noindent 
\bea
&&\al_{j,m}=\pm[m]_q^2(q-q^{-1})(\cE^{\pm j}_m-q^{\mp m}\cE^{\pm(j+1)}_m), \\
&&\al_{N,m}=[m]_q(q^{m/2}-q^{-m/2})(q^{-m/2}\cE^{+N}_m-q^{m/2} \cE^{-N}_m).
\ena
\end{prop}

The following formulae are also useful.
\begin{prop}
For $1\leq i,j\leq N$, the following commutation relations hold. 
\bea
&&[\al_{i,m}, \cE^{\pm j}_n]=\pm \delta_{m+n,0}\frac{[cm]_q}{m(q^m-q^{-m})}\frac{1-p^m}{1-p^{*m}}q^{-cm}(q^{\mp m}\delta_{i,j}-\delta_{i,j-1}),\\
&&[\cE^{\pm i}_m, e_j(z)]=\pm\frac{q^{-cm}z^m}{m(q^m-q^{-m})}\frac{1-p^m}{1-p^{*m}}e_j(z)(q^{\pm m}\delta_{i,j}-\delta_{i-1,j}),\\
&&[\cE^{\pm i}_m, f_j(z)]=\mp\frac{z^m}{m(q^m-q^{-m})}f_j(z)(q^{\pm m}\delta_{i,j}-\delta_{i-1,j}).
\ena
\end{prop}

\subsection{The elliptic currents $k_{\pm j}(z)$}
Let us set 
\bea
&&\psi_j(z)=:\exp\left\{(q-q^{-1})\sum_{m\not=0}\frac{\al_{j,m}}{1-p^m}p^mz^{-m}\right\}:.
\ena
Then the elliptic currents $\psi_j^\pm(z)$ in Definition \ref{defUqp} can be written as 
\bea
&&\psi^+_j(q^{-\frac{c}{2}}z)=K^+_j\psi_j(z),\qquad \psi^-_j(q^{-\frac{c}{2}}z)=K^-_j\psi_j(pq^{-c}z).
\ena
Let us introduce the new elliptic currents $k_{\pm j}(z) \,  (j\in \{0\}\cup I)$ associated with $\cE^{\pm j}_m$ by
\bea
k_{\pm j}(z) &=& :\exp\left\{\sum_{m\not=0} \frac{[m]_q^2(q-q^{-1})^2}{1-p^m}p^m\cE^{\pm j}_mz^{-m} \right\}: ,\\
k_0(z) &=&:k_{-N}(q^{-1/2}z)\psi_N(q^{-1/2}z):=:k_{+N}(q^{1/2}z)\psi_N(q^{1/2}z)^{-1}:.
\ena
Then from Proposition \ref{alphaAA} we have the following decompositions.  
\begin{prop}\lb{psikk}
\bea
\psi_j(z)&=&:k_{+j}(z)k_{+(j+1)}(qz)^{-1}:=:k_{-j}(z)^{-1}k_{-(j+1)}(q^{-1}z):,\\
\psi_N(z)&=&:k_{+N}(z)k_0(q^{-1/2}z)^{-1}:=:k_{-N}(z)^{-1}k_0(q^{1/2}z):.
\ena
\end{prop}

In addition, from Proposition \ref{cEcE} we obtain the following commutation relations. 
\begin{thm}\lb{kk}
\be
&&k_{\pm j}(z_1)k_{\pm j}(z_2)=\frac{\tilde{\rho}^{+*}(z)}{\tilde{\rho}^+(z)}k_{\pm j}(z_2)k_{\pm j}(z_1),
\qquad  (1\leq j\leq N),\\
&&k_{+j}(q^jz_1)k_{+k}(q^kz_2)=\frac{\tilde{\rho}^{+*}(z)}{\tilde{\rho}^+(z)}\frac{\Theta_{p^*}(q^{-2}z)\Theta_{p}(z)}{\Theta_{p^*}(z)\Theta_{p}(q^{-2}z)} k_{+k}(q^kz_2)k_{+j}(q^jz_1) \qquad \qquad (1\leq j<k\leq N),\\
&&k_{-j}(q^{-j}z_1)k_{-k}(q^{-k}z_2)=\frac{\tilde{\rho}^{+*}(z)}{\tilde{\rho}^{+}(z)}\frac{\Theta_{p^*}(q^{-2}z)\Theta_{p}(z)}{\Theta_{p^*}(z)\Theta_{p}(q^{-2}z)} k_{-k}(q^{-k}z_2)k_{-j}(q^{-j}z_1) \qquad (1\leq k<j\leq N),\\
&&k_{+j}(q^{j}z_1)k_{-k}(q^{-k}\xi z_2)=\frac{\tilde{\rho}^{+*}(z)}{\tilde{\rho}^{+}(z)}\frac{\Theta_{p^*}(q^{-2}z)\Theta_{p}(z)}{\Theta_{p^*}(z)\Theta_{p}(q^{-2}z)} k_{-k}(q^{-k}\xi z_2)k_{+j}(q^{j}z_1) \qquad (j\not=k),\\
&&k_{+j}(q^jz_1)k_{-j}(q^{-j}\xi z_2)=\frac{\tilde{\rho}^{+*}(u)}{\tilde{\rho}^{+}(u)}
\frac{\Theta_{p^*}(q^{2j-2}\xi^{-1}z)\Theta_{p}(q^{2j}\xi^{-1}z)}{\Theta_{p^*}(q^{2j}\xi^{-1}z)\Theta_{p}(q^{2j-2}\xi^{-1}z)} \frac{\Theta_{p^*}(q^{-2}z)\Theta_{p}(z)}{\Theta_{p^*}(z)\Theta_{p}(q^{-2}z)} k_{-j}(q^{-j}\xi z_2)k_{+j}(q^jz_1),\\
&&k_0(z_1)k_0(z_2)=\frac{\tilde{\rho}^{+*}(u)}{\tilde{\rho}^{+}(u)}
\frac{\Theta_{p^*}(q^{-2}z)\Theta_{p}(q^{2}z)\Theta_{p^*}(qz)\Theta_{p}(q^{-1}z)}
{\Theta_{p^*}(q^{2}z)\Theta_{p}(q^{-2}z)\Theta_{p^*}(q^{-1}z)\Theta_{p}(qz)}
k_0(z_2)k_0(z_1),\\
&&k_{+j}(q^jz_1)k_0(q^{N-1/2}z_2)=\frac{\tilde{\rho}^{+*}(u)}{\tilde{\rho}^{+}(u)}
\frac{\Theta_{p^*}(q^{-2}z)\Theta_{p}(z)}
{\Theta_{p^*}(z)\Theta_{p}(q^{-2}z)}
k_0(q^{N-1/2}z_2)k_{+j}(q^jz_1)\quad (1\leq j\leq N),\\
&&k_{-j}(\xi q^{-j}z_1)k_0(q^{N-1/2}z_2)=\frac{\tilde{\rho}^{+*}(u)}{\tilde{\rho}^{+}(u)}
\frac{\Theta_{p^*}(z)\Theta_{p}(q^{2}z)}
{\Theta_{p^*}(q^{2}z)\Theta_{p}(z)}
k_0(q^{N-1/2}z_2)k_{-j}(\xi q^{-j}z_1)\quad (1\leq j\leq N),
\en
where $z=z_1/z_2$, and $\tilde\rho^+(z)$ is a function which appears associated with 
the elliptic dynamical $R$-matrices\cite{Konno06}. (See $\S \ref{R-mat}$ )
\bea
\tilde{\rho}^+(z)
&=&\frac{\{\xi z\}^{2}\{\xi^2q^{-2} z\}\{q^2 z\}}{\{\xi^2 z\}\{z\}\{\xi q^2 z\}\{\xi q^{-2} z\}}
\frac{\{p\xi^2/z\}\{p/z\}\{p\xi q^2/z\}\{p\xi q^{-2}/z\}}{\{p\xi/z\}^{2}\{p\xi^2q^{-2}/z\}\{pq^2/z\}}
\ena
where $\xi=q^{-2\eta}$, $\{z\}=(z;p,\xi^2)_\infty$. We also set $\tilde{\rho}^{+*}(z)=\tilde{\rho}^+(z)|_{p\mapsto p^*}$.  

\end{thm}

\begin{prop}\lb{kef}
\be
&&k_{\pm j}(z_1)e_j(z_2)=\frac{\Theta_{p^*}(q^{-c}z)}{\Theta_{p^*}(q^{-c\mp 2}z)}
e_j(z_2)k_{\pm j}(z_1)\qquad (1\leq j\leq N),\\
&&k_{\pm j}(z_1)e_{j-1}(z_2)=\frac{\Theta_{p^*}(q^{-c\mp 1}z)}{\Theta_{p^*}(q^{-c\pm 1}z)}
e_{j-1}(z_2)k_{\pm j}(z_1)\qquad (2\leq j\leq N),\\
&&k_{\pm j}(z_1)e_k(z_2)=e_k(z_2)k_{\pm j}(z_1)\qquad(k\not=j,j-1), \\
&&k_{\pm j}(z_1)f_j(z_2)=\frac{\Theta_{p}(q^{\mp 2}z)}{\Theta_{p}(z)}
f_j(z_2)k_{\pm j}(z_1)\qquad (1\leq j\leq N),\\
&&k_{\pm j}(z_1)f_{j-1}(z_2)=\frac{\Theta_{p}(q^{\pm 1}z)}{\Theta_{p}(q^{\mp 1}z)}
f_{j-1}(z_2)k_{\pm j}(z_1)\qquad (2\leq j\leq N),\\
&&k_{\pm j}(z_1)f_k(z_2)=f_k(z_2)k_{\pm j}(z_1)\qquad(k\not=j,j-1),\\
&&k_0(q^{N-1/2}z_1)e_N(z_2)=\frac{\Theta_{p^*}(q^{-c+N}z)\Theta_{p^*}(q^{-c+N-1}z)}
{\Theta_{p^*}(q^{-c+N-2}z)\Theta_{p}(q^{-c+N+1}z)}e_N(z_2)k_0(q^{N-1/2}z_1),\\
&&k_0(q^{N-1/2}z_1)e_j(z_2)=e_j(z_2)k_0(q^{N-1/2}z_1)\qquad (1\leq j\leq N-1),\\
&&k_0(q^{N-1/2}z_1)f_N(z_2)=\frac{\Theta_{p}(q^{N-2}z)\Theta_{p}(q^{N+1}z)}
{\Theta_{p}(q^{N}z)\Theta_{p}(q^{N-1}z)}f_N(z_2)k_0(q^{N-1/2}z_1),\\
&&k_0(q^{N-1/2}z_1)f_j(z_2)=f_j(z_2)k_0(q^{N-1/2}z_1)\quad (1\leq j\leq N-1).
\en
\end{prop}


\subsection{The $H$-algebra $U_{q,p}(B_N^{(1)})$}
Let $\cA$ be a complex associative algebra, $\cH$ be a finite dimensional commutative subalgebra of 
$\cA$, and $\cM_{\cH^*}$ be the 
field of meromorphic functions on $\cH^*$ the dual space of $\cH$. 

\begin{dfn}[$\cH$-algebra\cite{EV}] 
An $\cH$-algebra is an associative algebra $\cA$ with 1, which is bigraded over 
$\cH^*$, $\ds{\cA=\bigoplus_{\alpha,\beta\in \cH^*} \cA_{\al,\beta}}$, and equipped with two 
algebra embeddings $\mu_l, \mu_r : \cM_{\cH^*}\to \cA_{0,0}$ (the left and right moment maps), such that 
\be
\mu_l(\hf)a=a \mu_l(T_\al \hf), \quad \mu_r(\hf)a=a \mu_r(T_\beta \hf), \qquad 
a\in \cA_{\al,\beta},\ \hf\in \cM_{\cH^*},
\en
where $T_\al$ denotes the automorphism $(T_\al \hf)(\la)=\hf(\la+\al)$ of $\cM_{\cH^*}$.
\end{dfn}

\begin{prop}
$U=U_{q,p}(B_N^{(1)})$  is a $H$-algebra by 
\be
&&U=\bigoplus_{\al,\beta\in H^*}U_{\al,\beta}\\
&&U_{\al,\beta}=\left\{x\in U \left|\ q^{P+h}x q^{-(P+h)}=q^{<\al,P+h>}x,\quad q^{P}x q^{-P}=q^{<\beta,P>}x\ \forall P+h, P\in H\right.\right\}
\en
and $\mu_l, \mu_r : \FF \to U_{0,0}$ defined by 
\be
&&\mu_l(\hf)=\hf(P+h,p)\in \FF[[p]],\qquad \mu_r(\hf)=\hf(P,p^*)\in \FF[[p]].
\en
\end{prop}
We regard ${T}_{\al}=e^\al\in \C[\cR_Q]$ as 
the shift operator $\cM_{{H}^*}\to \cM_{{H}^*}$ 
\be
({T}_{\al}\widehat{f})=e^{\al}\hf(P,p^*)e^{-\al}={\hf}(P+<\al,P>,p^*).
\en
Hereafter we abbreviate 
$f(P+h,p)$ and $f(P,p^*)$ as $f(P+h)$ and
 $f^*(P)$, respectively.


We also consider the $H$-algebra of the shift operators\cite{EV} 
\be
&&\cD=\{\ \sum_\al \widehat{f}_\al{T}_{\al}\ |\ \widehat{f}_\al\in {M}_{H^*}, 
\al\in \cR_Q\ \},\\
&&\cD_{\ha,\ha}=\{\ \widehat{f}{T}_{-\ha}\ \},\quad \cD_{\ha,{\beta}}=0\ 
(\al\not=\beta),\\
&&\mu_l^{\cD}(\widehat{f})=
\mu_r^{\cD}(\widehat{f})=\widehat{f}{T}_0 \qquad \widehat{f}\in {M}_{H^*}.
\en
Then we have the $H$-algebra isomorphism 
\bea
 U\cong U\tot\cD\cong \cD\tot U. \lb{Diso}
\ena 

\section{The $L$-operators and  The Dynamical $RLL$-relations}
We introduce the elliptic dynamical $R$-matrix of the $B_N^{(1)}$ type as a certain gauge transformation of 
Jimbo-Miwa-Okado's face type Boltzmann weight given in \cite{JMO}. Then we propose a construction of the $L$-operator satisfying the 
$RLL$-relation by means of the elliptic currents of $\UqpBN$. 

Hereafter we regard $q, p$ as a generic complex number satisfying $|q|, |p|<1$ and 
set $p=q^{2r}$.  We also use the following theta functions.
\bea
&&[u]=q^{\frac{u^2}{r}-u}\Theta_p(z), \qquad  
[u]^*=q^{\frac{u^2}{r^*}-u}\Theta_{p^*}(z), \lb{thetafull}
\ena
where we set $z=q^{2u}$.

\subsection{The elliptic dynamical $R$-matrix of the $B_N^{(1)}$ type}\lb{R-mat}

Let $\cI=\{0,\pm1,\pm2,\cdots, \pm N\}$. 
We fix the order $1\prec 2\prec \cdots \prec N\prec 0\prec -N\prec \cdots \prec -2\prec -1$. 
Let us consider the elliptic dynamical $R$-matrix  of the $B_N^{(1)}$ type  given by
\begin{eqnarray}
R^+(u,s)&=&\hrho^+(u)\bar{R}^+(u,s),\\
\bar{R}^+(u,s)&=&\left\{
\sum_{j=1\atop j\not=0}^{-1}
E_{j,j}\otimes E_{j,j}+
\sum_{1 \preceq j_1\prec j_2 \preceq -1\atop j_2\not=-j_1}
\left(b_{}(u,s_{j_1j_2 })
E_{j_1,j_1}
\otimes E_{j_2,j_2}+
\bar{b}_{}(u)
E_{j_2,j_2}\otimes E_{j_1,j_1}
\right)\right. \nonumber\\
&&\qquad \left.+
\sum_{1 \preceq j_1\prec j_2 \preceq -1\atop j_2\not=-j_1}
\left(c_{}(u,s_{j_1j_2 })
E_{j_1,j_2}\otimes E_{j_2,j_1}+
\bar{c}_{}
(u,s_{j_1j_2 })E_{j_2,j_1}\otimes E_{j_1,j_2}
\right)\right.\nn\\
&&\qquad\left.+\sum_{1 \preceq j_1\prec j_2 \preceq -1}
\left(d_{}(u,s_{j_1},s_{j_2 })
E_{-j_2,j_1}\otimes E_{j_2,-j_1}+
\bar{d}_{}(u,s_{j_1},s_{j_2 })E_{-j_1,j_2}\otimes E_{j_1,-j_2}
\right)\right.\nn\\
&&\qquad\left.
+\sum_{j\in\{ 1,2,\cdots,-2,-1\}}e_j(u,s)E_{-j,j}\otimes E_{j,-j}
\right\},\lb{ellR}
\end{eqnarray}
where $s=P, P+h$, we set $s_{\pm j}\equiv \pm s_{\ep_j}$ for $1\preceq j\preceq N$,  $s_{ij}=s_i-s_j$, $s_0=-\frac{1}{2}$, and  
\be
&&{\hrho}^+(u)=\xi^{-\frac{1}{r}}C(u,\xi)^{\frac{1}{2}}{\rho}^+_0(u)
,\\
&&{\rho}_0^+(u)=q^{-1}z^{\frac{1}{r}}\tilde{\rho}^+(u),\qquad 
C(u,\xi)=\frac{\Theta_{\xi^2}(z)^2\Theta_{\xi^2}(\xi q^2z)\Theta_{\xi^2}(\xi q^{-2}z)}
{\Theta_{\xi^2}(\xi z)^2\Theta_{\xi^2}(q^2z)\Theta_{\xi^2}(q^{-2}z)},\\
&&b(u,s)=
\frac{[s+1][s-1][u]}{[s]^2[u+1]},\qquad 
\bar{b}(u)=
\frac{[u]}{[u+1]},\\
&&c(u,s)=\frac{[1][s+u]}{[s][u+1]},\qquad \bar{c}(u,s)=\frac{[1][s-u]}{[s][u+1]},
\en
\be
&&d(u,s_{j},s_{k})=G_{s_j}\frac{[u][1][s_{j}+s_{k}+1+\eta-u]}{[\eta-u][u+1][s_j+s_k+1]}\prod_{m=1}^{j-1}\frac{[s_j-s_m]}{[s_j-s_m+1]}\prod_{m=1}^{k-1}\frac{[s_k-s_m+1]}{[s_k-s_m]}\\
&&\quad (j\prec k \preceq 0)
\\
&&d(u,s_{-k},s_{-j})=G_{s_{-j}}\frac{[u][1][s_{-j}+s_{-k}+1+\eta-u]}{[\eta-u][u+1][s_{-j}+s_{-k}+1]}\prod_{m=1}^{j-1}\frac{[s_{-j}+s_m]}{[s_{-j}+s_m+1]}\prod_{m=1}^{k-1}\frac{[s_{-k}+s_m+1]}{[s_{-k}+s_m]}\\
&&\quad (0\preceq -k\prec -j)
\\
&&d(u,s_{j},s_{-k})=G_{s_j}G_{s_{-k}}\frac{[u][1][s_{j}+s_{-k}+1+\eta-u]}{[\eta-u][u+1][s_j+s_{-k}+1]}\prod_{m=1}^{j-1}\frac{[s_j-s_m]}{[s_j-s_m+1]}
\prod_{m=1}^{k-1}\frac{[s_{-k}+s_m]}{[s_{-k}+s_m+1]}\\
&&\quad (j\prec 0\prec -k),
\en
\be
&&\bar{d}(u,s_{j},s_{k})=G_{s_k}\frac{[u][1][s_{j}+s_{k}+1+\eta-u]}{[\eta-u][u+1][s_j+s_k+1]}\prod_{m=1}^{j-1}\frac{[s_j-s_m+1]}{[s_j-s_m]}\prod_{m=1}^{k-1}\frac{[s_k-s_m]}{[s_k-s_m+1]}\\
&&\quad (j\prec k \preceq 0),
\\
&&\bar{d}(u,s_{-k},s_{-j})=G_{s_{-k}}\frac{[u][1][s_{-j}+s_{-k}+1+\eta-u]}{[\eta-u][u+1][s_{-j}+s_{-k}+1]}\prod_{m=1}^{j-1}\frac{[s_{-j}+s_m+1]}{[s_{-j}+s_m]}\prod_{m=1}^{k-1}\frac{[s_{-k}+s_m]}{[s_{-k}+s_m+1]}\\
&&\quad (0\preceq -k\prec -j),
\\
&&\bar{d}(u,s_{j},s_{-k})=\frac{[u][1][s_{j}+s_{-k}+1+\eta-u]}{[\eta-u][u+1][s_j+s_{-k}+1]}\prod_{m=1}^{j-1}\frac{[s_j-s_m+1]}{[s_j-s_m]}
\prod_{m=1}^{k-1}\frac{[s_{-k}+s_m+1]}{[s_{-k}+s_m]}\\
&& (j\preceq 0\preceq -k),
\en
\be
&&e_j(u,s)=\frac{[1][2s_j+1-u]}{[u+1][2s_j+1]}+
\frac{[u][1][2s_j+1+\eta-u]}{[\eta-u][u+1][2s_j+1]}G_{s_j}\qquad (j\not=0),\\
&&{e}_0(u,s)=\frac{[\eta+u][1][2\eta-u]}{[\eta-u][u+1][2\eta]}-
\frac{[u][1]}{[u+1][2\eta]}H_{s}.
\en
where for $k, -k=0$ the product $\ds{\prod_{m=1}^{k-1}}$ should be understood as 
$\ds{\prod_{m=1}^N}$ etc.  We also set 
\be
&&G_{s_j}=\frac{[s_j+1]}{[s_j]}\prod_{k\in \cI\atop
\not=\pm j,0}\frac{[s_j-s_k+1]}{[s_j-s_k]}
=\frac{[s_j+1]}{[s_j]}\prod_{m=1\atop \not=|j|}^N\frac{[s_j-s_m+1]}{[s_j-s_m]}\prod_{m=1\atop \not=|j|}^N\frac{[s_j+s_m+1]}{[s_j+s_m]},\\
&&H_s=\sum_{k\in \cI \atop
\not=0}\frac{[s_k+\frac{1}{2}+2\eta]}{[s_k+\frac{1}{2}]}G_{s_k}.
\en
Note that 
\be
G_{s_j-1}=G_{s_{-j}}^{-1},\qquad G_{s_{-j}-1}=G_{s_{j}}^{-1}, 
\en
and 
\bea
&&\hrho^+(u)\hrho^+(-u)=1,\qquad \hrho^+(\eta-u)=\hrho^+(u)\frac{[u][\eta-u+1]}{[u+1][\eta-u]}.\label{rhoinversion}
\ena

The matrix $R^+(u,s)$ in \eqref{ellR} is related to  
Jimbo-Miwa-Okado's face type Botzmann weight\cite{JMO} by the gauge transformation.
For $j\in \cI$, we set $\hj=\ep_{j}$ for $1\preceq j\preceq N$,  
$\hj=-\ep_{|j|}$ for $-N\preceq j\preceq -1$ and $\widehat{0}=0$. 
Let us define $\displaystyle{F(s,s+\hj)=\left(\frac{G_{s_j}}{G_{s_j}(j)}\right)^{\frac{1}{2}}}$ with 
\be
G_{s_j}(j)&=&\left\{
\mmatrix{\mbox{$\displaystyle{\prod_{m=1}^{j-1}\frac{[s_j-s_m+1]}{[s_j-s_m]}}$}
&\mbox{if}\ j\prec 0\cr
\mbox{$\displaystyle{\frac{[s_j+1]}{[s_j]}\prod_{m=1\atop \not=|j|}^{N}\frac{[s_j-s_m+1]}{[s_j-s_m]}
\prod_{m=|j|+1}^{N}\frac{[s_j+s_m+1]}{[s_j+s_m]}}$}&
\mbox{if}\ 0\prec j 
\cr
}\right.
\en

For $a\in \widetilde{\hh}^*$, $\rho=\Lambda_0+\Lambda_1+\cdots+\Lambda_N$  and $a_j=<a+\rho,\hj>$, we 
identify $a_j$ with $s_j$. 
Then we have  
\bea
&&R^+(u,s)^{ij}_{kl}=W\left(\left.\mmatrix{a&a+\hi\cr a+\hl&a+\hi+\hj}\right|\ u\right)
\qquad (i+j=k+l),\nn\\
&&W\left(\left.\mmatrix{a&a+\hi\cr a+\hl&a+\hi+\hj}\right|\ u\right)
=\hrho^+(u)\frac{[\eta][1]}{[\eta-u][u+1]}\frac{F(a,a+\hi)F(a+\hi,a+\hi+\hj)}{
F(a,a+\hl)F(a+\hl,a+\hl+\hk)}\nn\\
&&\qquad\qquad\qquad\qquad \qquad\qquad\qquad \qquad \times W_{JMO}\left(\left.\mmatrix{a&a+\hi\cr a+\hl&a+\hi+\hj}\right|\ u\right)\lb{gauge}
\ena
One can derive the following relations from (2.10)-(2.13b) in \cite{JMO}.  
\begin{itemize}
\item[1)] Crossing symmetry
\bea
&&W\left(\left.\mmatrix{a&b\cr d&c}\right|\ u\right)
=\frac{F(b,c)F(c,b)}{
F(a,d)F(d,a)}\left(\frac{G_bG_d}{G_aG_c}\right)^{\frac{1}{2}} W\left(\left.\mmatrix{d&a\cr c&b}\right|\ \eta-u\right).
\lb{cross}\ena
\item[2)] Reflection symmetry
\be
&&W\left(\left.\mmatrix{a&b\cr d&c}\right|\ u\right)
=\left(\frac{F(a,b)F(b,c)}{
F(a,d)F(d,c)}\right)^2 W\left(\left.\mmatrix{a&d\cr b&c}\right|\ u\right).
\en
\item[3)] Unitarity
\be
&&\sum_{g}W\left(\left.\mmatrix{a&g\cr d&c}\right|\ u\right)W\left(\left.\mmatrix{a&b\cr g&c}\right|\ -u\right)=\delta_{bd}. 
\en
\item[4)] 2nd inversion relation
\be
&&\sum_{g}\left(\frac{G_aG_g}{G_bG_d}\right)^{\frac{1}{2}}W\left(\left.\mmatrix{a&b\cr d&g}\right|\ \eta- u\right)
\left(\frac{G_cG_g}{G_bG_d}\right)^{\frac{1}{2}}W\left(\left.\mmatrix{c&d\cr b&g}\right|\  \eta+u\right)
=\delta_{ac}. 
\en
\end{itemize}
Here 
\bea
&&G_a=\vep(a)\prod_{j=1}^N[a_j]\prod_{1\leq i<j\leq N}[a_i-a_j][a_i+a_j] \lb{gaugeG}
\ena
and $\vep(a)$ is a sign factor such that $\vep(a+\hj)/\vep(a)=1$. 

\noindent
{\it Remark.}\ The choice of the gauge \eqref{gauge} and the resultant $R$-matrix \eqref{ellR} is convenient to discuss the $RLL$-relations in the next sections, because it allows the $L$-operator $\widehat{L}^-(u)$ to be related to $\widehat{L}^+(u)$ simply by $\widehat{L}^-(u)=\widehat{L}^+(u+r-\frac{c}{2})$, i.e. one needs no extra modifications follow from Proposition 4.3 in \cite{JKOStransfG} (4.8). Note that $p$ in \cite{JKOStransfG} is $p^*$ in the present paper. See also \cite{JKOS99}. One drawback is that one needs to introduce a set of extra generators and a central extension to the group algebra $\C[\cR_Q]$ in order to remove constant gauge factors such as $q^{ 1/r}$ and $q^{ 1/r^*}$ in a realization of the proper modified elliptic currents, which  will be discussed in the next section. 
See Remark below Proposition 3.10 in \cite{KojimaKonno}. However in order to avoid such unessential complications, we hereafter treat the whole formulas up to those  
constant gauge factors.

\subsection{Modified elliptic currents}\lb{MEC}
Since our elliptic $R$-matrix is given by the theta functions \eqref{thetafull} accompanied by 
the fractional power of $z$, 
we need to introduce the following modifications of the elliptic currents.
\be
\!\!\!&&E_{ j}(u)=e_j(z) z^{-\frac{ P_{\al_j}-1}{r^*}}\qquad (1\leq j\leq N-1),\\
\!\!\!&&E_{ N}(u)=e_N(z) z^{-\frac{ P_{\al_N}-1/2}{r^*}},\\
\!\!\!&&F_{ j}(u)=f_j(z) z^{\frac{ P_{\al_j}+ h_{\al_j}-1}{r}}\qquad (1\leq j\leq N-1),\\
\!\!\!&&F_{ N}(u)=f_N(z) z^{\frac{ P_{\al_N}+ h_{\al_N}-1/2}{r}},\\
\!\!\!&&H^{\pm}_{j}(z)=\psi^{\pm}_{j}(z)(K_j^\pm)^{-1}e^{-Q_{\al_j}} 
(q^{\pm(r-\frac{c}{2})}z)^{-\frac{r-r^*}{rr^*}( P_{\al_j}-1)
+\frac{1}{r} h_{\al_j}},
\en
and 
\be
K^+_{+ j}(u)&=&k_{+j}(q^{j}z)e^{-Q_{\ep_j}}(q^{j}zq^{-r})^{-\frac{r-r^*}{rr^*}
P_{\ep_j}+\frac{1}{r} h_{\ep_j}},\\
K^+_{-j}(u)&=&k_{-j}(q^{-j}\xi z)e^{Q_{\ep_j}}(q^{-j}\xi zq^{-r})^{\frac{r-r^*}{rr^*}
P_{\ep_j}
-\frac{1}{r} h_{\ep_j}},\\
K^+_0(u)&=&k_0(\xi^{1/2}z),\\
K^-_{\pm j}(u)&=&K^+_{\pm j}\left(u+r-\frac{c}{2}\right).  
\en
for $1\leq j\leq N$. 
We also set
\be
\td=\hd+\frac{1}{2r^*}\sum_{j=1}^N(P_j+2)P^j-\frac{1}{2r}\sum_{j=1}^N((P+h)_j+2)(P+h)^j.
\en
We have 
\begin{prop}
\be
H^{\pm}_{j}(u)&=&:K^\pm_{+j}\left(u+\frac{c}{4}-\frac{j}{2}\right)K^\pm_{+(j+1)}\left(u+\frac{c}{4}-\frac{j}{2}\right)^{-1}:,
\\
&=&:K^\pm_{-j}\left(u+\frac{c}{4}+\frac{j}{2}+\eta\right)^{-1}K^\pm_{-(j+1)}\left(u+\frac{c}{4}+\frac{j}{2}+\eta\right):
\qquad (1\leq j\leq N-1),\\
H^{\pm}_{N}(u)&=&:K^\pm_{+N}\left(u+\frac{c}{4}-\frac{N}{2}\right)
K^\pm_{0}\left(u+\frac{c}{4}-\frac{N}{2}\right)^{-1}:,
\\
&=&:K^\pm_{-N}\left(u+\frac{c}{4}+\frac{N}{2}+\eta\right)^{-1}
K^\pm_{0}\left(u+\frac{c}{4}+\frac{N}{2}+\eta\right):.
\en
\end{prop}
Then one can rewrite the formulas in Theorem \ref{kk} and Proposition \ref{kef} as follows.
\begin{prop}\lb{relKK}
\be
&&K^+_{\pm j}(u_1)K^+_{\pm j}(u_2)=
\frac{{\rho}^{+*}(u_1-u_2)}{{\rho}^+(u_1-u_2)}K^+_{\pm j}(u_2)K^+_{\pm j}(u_1),\\
&&K^+_{+j}(u_1)K^+_{+l}(u_2)=
\frac{{\rho}^{+*}(u_1-u_2)}{{\rho}^+(u_1-u_2)}\frac{[u_1-u_2-1]^*[u_1-u_2]}{[u_1-u_2]^*[u_1-u_2-1]}K^+_{+l}(u_2)K^+_{+j}(u_1) \quad (1\preceq j\prec l\preceq 0),\\
&&K^+_{-j}(u_1)K^+_{-l}(u_2)=
\frac{{\rho}^+(u_1-u_2)}{{\rho}^{+*}(u_1-u_2)}\frac{[u_1-u_2]^*[u_1-u_2-1]}{[u_1-u_2-1]^*[u_1-u_2]}K^+_{-l}(u_2)K^+_{-j}(u_1) \quad (1\preceq j\prec l\preceq 0),\\
&&K^+_{+j}(u_1)K^+_{-j}(u_2)=
\frac{{\rho}^{+*}(u_1-u_2)}{{\rho}^+(u_1-u_2)}\frac{[u_1-u_2+\eta+j-1]^*[u_1-u_2-1]^*}{[u_1-u_2+\eta+j]^*[u_1-u_2]^*}\\
&&\qquad\qquad\qquad\qquad\qquad  \times \frac{[u_1-u_2+\eta+j][u_1-u_2]}{[u_1-u_2+\eta+j-1][u_1-u_2-1]}K^+_{-l}(u_2)K^+_{+j}(u_1),\\
&&K^+_{+j}(u_1)K^+_{-l}(u_2)=
\frac{{\rho}^{+*}(u_1-u_2)}{{\rho}^+(u_1-u_2)}\frac{[u_1-u_2-1]^*[u_1-u_2]}{[u_1-u_2]^*[u_1-u_2-1]}K^+_{-l}(u_2)K^+_{+j}(u_1) \quad (1\preceq j, l\preceq N,  j\not= l),\\
&&K_0^+(u_1)K_0^+(u_2)=
\frac{{\rho}^{+*}(u_1-u_2)}{{\rho}^+(u_1-u_2)}
\frac{[u_1-u_2-1]^*[u_1-u_2+\frac{1}{2}]^*}
{[u_1-u_2+1]^*[u_1-u_2-\frac{1}{2}]^*} \\
&&\qquad\qquad\qquad\qquad\qquad \times 
\frac{[u_1-u_2+1][u_1-u_2-\frac{1}{2}]}
{[u_1-u_2-1][u_1-u_2+\frac{1}{2}]}K_0^+(u_2)K_0^+(u_1).
\en
\end{prop}

\begin{prop}\lb{relEK}
\be
&&K^+_{+j}(u_1)E_j(u_2)=
\frac{\left[u_1-u_2+\frac{j-c}{2}\right]^*}{\left[u_1-u_2+\frac{j-c}{2}-1\right]^*}E_j(u_2)K^+_{+j}(u_1)\qquad (1\leq j\leq N),\\
&&K^+_{+j}(u_1)E_{j-1}(u_2)=
\frac{\left[u_1-u_2+\frac{j-1-c}{2}\right]^*}{\left[u_1-u_2+\frac{j-1-c}{2}+1\right]^*}E_{j-1}(u_2)K^+_{+j}(u_1)\qquad (2\leq j\leq N),\\
&&K^+_j(u_1)E_l(u_2)=E_l(u_2)K^+_j(u_1)\qquad (l\not=j,j-1),
\en
\be
&&K^+_{-j}(u_1)E_j(u_2)=
\frac{\left[u_1-u_2-\frac{j+c}{2}-\eta\right]^*}{\left[u_1-u_2-\frac{j+c}{2}-\eta+1\right]^*}E_j(u_2)K^+_{-j}(u_1)\qquad (1\leq j\leq N),\\
&&K^+_{-j}(u_1)E_{j-1}(u_2)=\frac{\left[u_1-u_2-\frac{j-1+c}{2}-\eta\right]^*}{\left[u_1-u_2-\frac{j-1+c}{2}-\eta-1\right]^*}E_{j-1}(u_2)K^+_{-j}(u_1)\qquad (2\leq j\leq N),\\
&&K^+_{-j}(u_1)E_l(u_2)=E_l(u_2)K^+_{-j}(u_1)\qquad (l\not=j,j-1),
\en
\be
&&K^+_{+j}(u_1)F_j(u_2)=
\frac{\left[u_1-u_2+\frac{j}{2}-1\right]}{\left[u_1-u_2+\frac{j}{2}\right]}F_j(u_2)K^+_{+j}(u_1)\qquad (1\leq j\leq N),\\
&&K^+_{+j}(u_1)F_{j-1}(u_2)=\frac{\left[u_1-u_2+\frac{j+1}{2}\right]}{\left[u_1-u_2+\frac{j+1}{2}-1\right]}F_{j-1}(u_2)K^+_{+j}(u_1)\qquad (2\leq j\leq N),\\
&&K^+_j(u_1)F_l(u_2)=F_l(u_2)K^+_j(u_1)\qquad (l\not=j,j-1),\\
&&K^+_{-j}(u_1)F_j(u_2)=
\frac{\left[u_1-u_2-\frac{j}{2}-\eta+1\right]}{\left[u_1-u_2-\frac{j}{2}-\eta\right]}F_j(u_2)K^+_{-j}(u_1)\qquad (1\leq j\leq N),\\
&&K^+_{-j}(u_1)F_{j-1}(u_2)=\frac{\left[u_1-u_2-\frac{j+1}{2}-\eta\right]}{\left[u_1-u_2-\frac{j+1}{2}-\eta+1\right]}F_{j-1}(u_2)K^+_{-j}(u_1)\qquad (2\leq j\leq N),\\
&&K^+_{-j}(u_1)F_l(u_2)=F_l(u_2)K^+_{-j}(u_1)\qquad (l\not=j,j-1),
\en
and
\be
&&K^+_{0}(u_1)E_{N}(u_2)=
\frac{\left[u_1-u_2+\frac{N-c}{2}\right]^*}{\left[u_1-u_2+\frac{N-c}{2}-1\right]^*}\frac{\left[u_1-u_2+\frac{N-c-1}{2}\right]^*}{\left[u_1-u_2+\frac{N-c+1}{2}\right]^*}E_{N}(u_2)K^+_{0}(u_1),\\
&&K^+_0(u_1)E_j(u_2)=E_j(u_2)K^+_0(u_1)\qquad (j\not=N,0),\\
&&K^+_{0}(u_1)F_{N}(u_2)=
\frac{\left[u_1-u_2+\frac{N}{2}-1\right]}{\left[u_1-u_2+\frac{N}{2}\right]}\frac{\left[u_1-u_2+\frac{N+1}{2}\right]}{\left[u_1-u_2+\frac{N-1}{2}\right]}F_{N}(u_2)K^+_{0}(u_1),\\
&&K^+_0(u_1)F_j(u_2)=F_j(u_2)K^+_0(u_1)\qquad (j\not=N,0).
\en
\end{prop}
In addition, the defining relations \eqref{ge}--\eqref{serref} of $U_{q,p}(B_N^{(1)})$ can be 
rewritten as follows in the sense of analytic continuation. 
\begin{prop}\lb{defrelfull}
\bea
&&[h_i,\al_{j,n}]=0,\quad [h_i,E_j(u)]=a_{ij} E_j(u),\quad 
[h_i,F_j(u)]=-a_{ij} F_j(u), 
\lb{u2}\\
&&[{\hd},h_i]=0,\quad [{\hd},\al_{i,n}]= n \al_{i,n},\quad\\ 
&&[{\hd},E_{i}(u)]=\left(-z\frac{\partial}{\partial z}+\frac{1}{r^*}\right)
E_i(u) , \quad
[{\hd},F_{i}(u)]=\left(-z\frac{\partial}{\partial z}+\frac{1}{r}\right)
F_i(u) ,
\lb{u3}
\ena
\bea
&& \left[u-v-\frac{b_{ij}}{2}\right]^* E_i(u)E_j(v)
=\left[u-v+\frac{b_{ij}}{2}\right]^* 
E_j(v)E_i(u),
\lb{u7}\\
&& \left[u-v+\frac{b_{ij}}{2} \right] F_i(u)F_j(v)
=\left[u-v-\frac{b_{ij}}{2}\right]
F_j(v)F_i(u),
\lb{u8}\\
&&[E_i(u),F_j(v)]=\frac{\delta_{i,j}}{q_i-q_i^{-1}}
\left(\delta\bigl(q^{-c}\frac{z}{w}\bigr)H^-_{i}(q^{c/2}w)
-\delta\bigl(q^{c}\frac{z}{w}\bigr)H^+_{i}(q^{-c/2}w)
\right),\lb{u9}
\ena
\bea
&&\sum_{\sigma\in S_a}z_{\sigma(1)}^{-\frac{1}{r^*}}\prod_{1\leq k<m\leq a}
\frac{(p^*q^2 z_{\sigma(m)}/
z_{\sigma(k)};p^*)_{\infty}}{(p^*q^{-2} z_{\sigma(m)}/
z_{\sigma(k)};p^*)_{\infty}}\nonumber\\
&&\qquad\times
\sum_{s=0}^a (-)^s\left[\mmatrix{a\cr
s\cr}\right]_i \prod_{s+1\leq m\leq a-1}\left(\frac{w}{z_{\sigma(m)}}\right)^{\frac{1}{r^*}}\prod_{a\leq m\leq s}\left(\frac{z_{\sigma(m)}}{w}\right)^{\frac{1}{r^*}}
\nn\\
&&\qquad\qquad\times\prod_{1\leq m\leq s} 
\frac{(p^*q^{b_{ij}} w/
z_{\sigma(m)};p^*)_{\infty}}{(p^*q^{-b_{ij}} w/
z_{\sigma(m)};p^*)_{\infty}}
\prod_{s+1\leq m\leq a}\frac{(p^*q^{b_{ij}} 
z_{\sigma(m)}/w;p^*)_{\infty}}{(p^*q^{-b_{ij}} 
z_{\sigma(m)}/w;p^*)_{\infty}}\\
&&\qquad\qquad\times E_{i}(u_{\sigma(1)})\cdots E_{i}(u_{\sigma(s)})
E_{j}(v) E_{i}(u_{\sigma(s+l)})\cdots E_{i}(u_{\sigma(a)})=0\qquad 
(i\not=j,\ a=1-a_{ij}),\nonumber\\
&&\sum_{\sigma\in S_a}z_{\sigma(1)}^{\frac{1}{r}}\prod_{1\leq k<m\leq a}
\frac{(pq^{-2} z_{\sigma(m)}/
z_{\sigma(k)};p)_{\infty}}{(pq^{2} z_{\sigma(m)}/
z_{\sigma(k)};p)_{\infty}}\nonumber\\
&&\qquad\times
\sum_{s=0}^a (-)^s\left[\mmatrix{a\cr
s\cr}\right]_i \prod_{s+1\leq m\leq a-1}\left(\frac{w}{z_{\sigma(m)}}\right)^{-\frac{1}{r}}\prod_{a\leq m\leq s}\left(\frac{z_{\sigma(m)}}{w}\right)^{-\frac{1}{r}}
\nn\\
&&\qquad\qquad\times\prod_{1\leq m\leq s} 
\frac{(pq^{-b_{ij}} w/
z_{\sigma(m)};p)_{\infty}}{(pq^{b_{ij}} w/
z_{\sigma(m)};p)_{\infty}}
\prod_{s+1\leq m\leq a}\frac{(pq^{-b_{ij}} 
z_{\sigma(m)}/w;p)_{\infty}}{(pq^{b_{ij}} 
z_{\sigma(m)}/w;p)_{\infty}}\\
&&\qquad\times F_{i}(u_{\sigma(1)})\cdots F_{i}(u_{\sigma(s)})
F_{j}(v) F_{i}(u_{\sigma(s+l)})\cdots F_{i}(u_{\sigma(a)})=0\qquad 
(i\not=j,\ a=1-a_{ij}).\nonumber
\ena
\end{prop}

\subsection{The half currents and the $L$-operators}
We next introduce the half currents $E^\pm_{i,j}(u), F^\pm_{j,i}(u)
\ (1\preceq i\prec j\preceq -1
)$ and propose a construction of the $L$-operators of $U_{q,p}(B_N^{(1)})$.

\begin{dfn}\lb{basicHC}
For $1\preceq j\preceq -2$, we define the basic half currents $E_{j+1,j}^+(u), 
F_{j,j+1}^+(u)
$ as follows. 
\be
&&F_{j,j+1}^+(u):=a_{j,j+1}\oint_C \frac{dz_j'}{2\pi i z_j'}
F_{j}(u_j')\frac{[u-u_j'+(P+h)_{j,j+1}+
\frac{j}{2}-1][1]}{[u-u_j'+
\frac{j}{2}][(P+h)_{j,j+1}-1]},\\
&&E_{j+1,j}^+(u):=
a_{j+1,j}^*\oint_{C^*} \frac{dz_j'}{2\pi i z_j'}
E_{j}(u_j')\frac{[u-u_j'+\frac{j-c}{2}+1-
P_{j,j+1}]^*[1]^*}
{[u-u_j'+\frac{j-c}{2}]^*
[P_{j,j+1}-1]^*}, 
\end{eqnarray*}
for $1\preceq j\preceq N$ and
\be
&&F_{-(j+1),-j}^+(u)
:=a_{-(j+1),-j}\oint_C \frac{dz_j'}{2\pi i z_j'}
F_{j}(u_j')\frac{[u-u_j'+(P+h)_{-(j+1),-j}-
\frac{j}{2}-\eta-1][1]}{[u-u_j'-
\frac{j}{2}-\eta][(P+h)_{-(j+1),-j}-1]},\\
&&E_{-j,-(j+1)}^+(u):=
a_{-j,-(j+1)}^*\oint_{C^*} \frac{dz_j'}{2\pi i z_j'}
E_{j}(u_j')\frac{[u-u_j'-\frac{j+c}{2}-\eta+1-
P_{-(j+1),-j}-\delta_{j,N}]^*[1]^*}
{[u-u_j'-\frac{j+c}{2}-\eta]^*
[P_{-(j+1),-j}-1+\delta_{j,N}]^*},~~
\end{eqnarray*}
for $-N\preceq -j\preceq -1$, where 
$N+1\equiv 0 \equiv -N-1$. We also define
\be
&&E^-_{j+1,j}(u):=E^+_{j+1,j}(u+r-\frac{c}{2}),\quad F^-_{j,j+1}(u):=F^+_{j,j+1}(u+r-\frac{c}{2})
. 
\en
\end{dfn}

By using Propositions \ref{relKK}-\ref{defrelfull}  one can derive the following relations. 
\begin{prop}\lb{relbasicHC}
\be
&&K_{j+1}^+(u_1)^{-1}E_{j+1,j}^+(u_2)K_{j+1}^+(u_1)
=E_{j+1,j}^+(u_2)\frac{1}{\bar{b}_{}^*(u)}
-E_{j+1,j}^+(u_1)\frac{c_{}^*(u,P_{j,j+1})}{\bar{b}_{}^*(u)},\lb{klelj}\\
%
%
%
%
&&K_{j+1}^+(u_1)F_{j,j+1}^+(u_2)
K_{j+1}^+(u_1)^{-1}=
\frac{1}{\bar{b}_{}(u)}F_{j,j+1}^+(u_2)-
\frac{\bar{c}_{}(u,P_{j,j+1}+h_{j,j+1})}{
\bar{b}_{}(u)}F_{j,j+1}^+(u_1),\lb{klfjl}\\
&&\frac{[1-u]^*}{[u]^*}E^+_{j+1,j}(u_1)E^+_{j+1,j}(u_2)+\frac{[1+u]^*}{[u]^*}
E^+_{j+1,j}(u_2)E^+_{j+1,j}(u_1)\nn\\
&&\qquad\qquad=E^+_{j+1,j}(u_1)^2\frac{[1]^*[P_{j,j+1}-2+u]^*}{[P_{j,j+1}-2]^*[u]^*}+
E^+_{j+1,j}(u_2)^2\frac{[1]^*[P_{j,j+1}-2-u]^*}{[P_{j,j+1}-2]^*[u]^*},\lb{eljelj}\\
%
%
%
%
&&\nn\\
&&\frac{[1+u]}{[u]}F^+_{j,j+1}(u_1)F^+_{j,j+1}(u_2)+\frac{[1-u]}{[u]}
F^+_{j,j+1}(u_2)F^+_{j,j+1}(u_1)\nn\\
&&\qquad\qquad=F^+_{j,j+1}(u_1)^2\frac{[1][P_{j,j+1}+h_{j,j+1}-2-u]}{[P_{j,j+1}+h_{j,j+1}-2][u]}+
F^+_{j,j+1}(u_2)^2\frac{[1][P_{j,j+1}+h_{j,j+1}-2+u]}{[P_{j,j+1}+h_{j,j+1}-2][u]},
\lb{fjlfjl}\\
&&\nn\\
&&[E_{j+1,j}^+(u_1),F_{j,j+1}^+(u_2)] \nn\\
&&\qquad\qquad=K_j^+(u_2)\frac{\bar{c}^*(u,P_{j,j+1})}{
\bar{b}^*(u)}K_{j+1}^+(u_2)^{-1}
-K_{j+1}^+(u_1)^{-1}
\frac{\bar{c}(u,(P+h)_{j,j+1})}{\bar{b}(u)}
K_j^+(u_1).
\en
\end{prop}
\begin{prop}
For $1\preceq j\preceq -2, j\not=N, 0$ the relations in Proposition \ref{relbasicHC} and those for $K^+_j(u), K^+_{j+1}(u)$ in Proposition \ref{relKK} coincide with 
the following $RLL$-relation of the $U_{q,p}(\glth)$ type. 
\begin{align}
R^{+(12)}_{j,j+1}(u_1-u_2,P+h)
&\widehat{L}^{+(1)}_{j,j+1}(u_1)\widehat{L}^{+(2)}_{j,j+1}(u_2)= \nn \\
& \widehat{L}^{+(2)}_{j,j+1}(u_2)\widehat{L}^{+(1)}_{j,j+1}(u_1)
R^{*+(12)}_{j,j+1}(u_1-u_2,P-h^{(1)}-h^{(2)}), \lb{RLLj}
\end{align}
where 
\be
R^{+}_{j,j+1}(u,s)&=&\hrho^+(u)\left(\mmatrix{1&0&0&0\cr
0&b_{}(u,s_{j,j+1 })&c_{}(u,s_{j,j+1})&0\cr
0&\bar{c}_{}
(u,s_{j,j+1})&\bar{b}_{}(u)&0\cr
0&0&0&1\cr}\right), \\
\widehat{L}_{j,j+1}(u)&=&\mat{1&F^+_{j,j+1}(u)\cr 0&1\cr}\mat{K_j^+(u)&0\cr 0&K_{j+1}^+(u)\cr}
\mat{1&0\cr
E^+_{j+1,j}(u)&1\cr}. 
\en
\end{prop}

\begin{dfn}
By means of  the basic half currents $E_{j+1,j}^\pm(u), 
F_{j,j+1}^\pm(u)
$ $(1\preceq j\preceq -2)$, we define the other half currents 
$E^\pm_{i,j}(u), F^\pm_{j,i}(u)\ (1\preceq i\prec j\preceq -1, j=i+1)$ 
by requiring the following conditions. 
\begin{itemize}
\item[1)] The half currents $E^\pm_{j,i}(u)$  and $F^\pm_{i,j}(u)\ (i\prec j)$ have the following 
series expansions.
\be
&&E^+_{j,i}(u)=\left(\sum_{n\in \Z_{\geq 0}}E^+_{j,i;-n}z^n+  
\sum_{n\in \Z_{> 0}}E^+_{j,i;n}p^nz^{-n}\right)z^{-\frac{P_{i,j}-1+\delta_{i,0}-\delta_{j,-i}}{r^*}},\\
&&F^+_{i,j}(u)=\left( \sum_{n\in \Z_{\geq 0}}F^+_{i,j;-n}z^n+  
\sum_{n\in \Z_{> 0}}F^+_{i,j;n}p^nz^{-n} \right)z^{\frac{(P+h)_{i,j}-1+\delta_{i,0}-\delta_{j,-i}}{r}}, 
\en
where $E^+_{j,i;\pm n}, F^+_{i,j;\pm n}\in U_{q,p}(B_N^{(1)})$ and
\bea
&&E^-_{j,i}(u):=E^+_{j,i}(u+r-\frac{c}{2}),\quad F^-_{i,j}(u):=F^+_{i,j}(u+r-\frac{c}{2}),\quad  
K^-_j(u)=K^+_j(u+r-\frac{c}{2}). \nn\\
&&\lb{EmFm}
\ena
\item[2)] 
Let us set  
\begin{eqnarray}
&&\hspace*{-2cm}\widehat{L}^+(u)=
\left(\begin{array}{ccccc}
1&F_{1,2}^+(u)&F_{1,3}^+(u)&\cdots&F_{1,-1}^+(u)\\
0&1&F_{2,3}^+(u)&\cdots&F_{2,-1}^+(u)\\
\vdots&\ddots&\ddots&\ddots&\vdots\\
\vdots&&\ddots&1&F_{-2,-1}^+(u)\\
0&\cdots&\cdots&0&1
\end{array}\right)
\left(
\begin{array}{cccc}
K^+_1(u)&0&\cdots&0\\
0&K^+_2(u)&&\vdots\\
\vdots&&\ddots&0\\
0&\cdots&0&K^+_{-1}(u)
\end{array}
\right)\nn\\
&\times&
\left(
\begin{array}{ccccc}
1&0&\cdots&\cdots&0\\
E^+_{2,1}(u)&1&\ddots&&\vdots\\
E^+_{3,1}(u)&E^+_{3,2}(u)&\ddots&\ddots&\vdots\\
\vdots&\vdots&\ddots&1&0\\
E^+_{-1,1}(u)&E^+_{-1,2}(u)
&\cdots&E^+_{-1,-2}(u)&1
\end{array}
\right).\lb{Lop}
\end{eqnarray}
Due to 1), the matrix products in \eqref{Lop}  are well defined in the $p$-adic 
topology. 
Then $\widehat{L}^+(u)$ and $\widehat{L}^-(u)=\widehat{L}^+(u+r-{c}/{2})$ satisfy the following $RLL$-relations.  
\begin{eqnarray}
&&R^{\pm(12)}(u,P+h)
\widehat{L}^{\pm(1)}(u_1)\widehat{L}^{\pm(2)}(u_2)=
\widehat{L}^{\pm(2)}(u_2)\widehat{L}^{\pm(1)}(u_1)
R^{*\pm(12)}(u,P-h^{(1)}-h^{(2)}),
\nn\\
\lb{RLLpm}\\
&&R^{\pm(12)}(u\pm\frac{c}{2},P+h)\hL^\pm(u_1) \hL^\mp(u_2)=\hL^\mp(u_2)\hL^\pm(u_1)R^{*\pm(12)}(u\mp\frac{c}{2},P-h^{(1)}-h^{(2)}). \nn\\
\lb{RLLmp}
\end{eqnarray}
Here 
\be
&&R^-(u,s)=R^+(u+r-\frac{c}{2},s). 
\en
Note also that $r^*+\frac{c}{2}=r-\frac{c}{2}$.
\end{itemize}
\end{dfn}

\noindent
{\it Remark.} Since $\rho^+(u)/\rho^{+*}(u)=\rho^+_0(u)/\rho^{+*}_0(u)$, 
the $RLL$-relations \eqref{RLLj}, \eqref{RLLpm}, \eqref{RLLmp} remain unchanged even 
if one uses $\rho^+_0(u)$ and $\rho_0^{+*}(u)$ instead of $\rho^+(u)$ and $\rho^{+*}(u)$, 
respectively. See Sec.\ref{VO}

\begin{conj}\lb{HCfromRLL}
The $RLL$-relation \eqref{RLLpm} and \eqref{RLLmp} determines the half currents 
$E^\pm_{j,i}(u), F^\pm_{i,j}(u)$ $(1\preceq i\prec j\preceq -1, j\not=i+1)$ recursively 
 and uniquely from the basic ones in Definition \ref{basicHC}.  
\end{conj}
In fact, the half currents $E^\pm_{j,i}(u), F^\pm_{i,j}(u)$ with $1\preceq i\prec j\prec 0$ or $0\prec i\prec j\preceq -1$ are determined recursively
 by the basic ones in the same way as for $U_{q,p}(A_N^{(1)})$ case\cite{DF,Konno12}. 
  As for the 
other half currents $E^\pm_{j,i}(u), F^\pm_{i,j}(u)$ with $1\preceq i\preceq 0\prec  j\preceq -1$ or 
$1\preceq i\prec 0\preceq  j\preceq -1$,  
we have observed that 
the combinations $E^+_{j,i}(u+c/4)-E^-_{j,i}(u-c/4)$ and $F^+_{i,j}(u-c/4)-F^-_{j,i}(u+c/4)$  satisfy a system of linear equations with the operator 
valued coefficients given by 
the {\it total} elliptic currents respectively.  In addition,  the half currents $E^\pm_{-1,1}(u)$ (reps. $F^\pm_{1,-1}(u)$) is determined by all the other half currents $E^\pm_{j,i}(u)$ (resp.  $F^\pm_{i,j}(u)$) $i\prec j$.  
An explicit expression for the half currents $E^\pm_{j,i}(u), F^\pm_{i,j}(u)$ ($1\preceq i\prec j\prec 0$ or $0\prec i\prec j\preceq -1$) and  conjectural expressions for the other half currents are given in Appendix B.

The existence of the operator $\hL^+(u)$ satisfying \eqref{RLLpm}-\eqref{RLLmp}, and hence the existence and the uniqueness of  the half currents,  can also be seen in the following argument. Consider the elliptic quantum group $\Bqla(\Bnh)$
realized by the Chevalley generators equipped with the quasi-Hopf  algebra structure\cite{JKOStransfG}.  See Appendix A. 
Note that $\Bqla(\gh)$ is isomorphic to the Drinfeld-Jimbo's quantum affine algebra  $U_q(\gh)$ \cite{Drinfeld85,Jimbo} as an associative algebra. 
In addition, we have shown the isomorphism\cite{FKO}  
\be
&&U_{q,p}(\gh)/pU_{q,p}(\gh)\cong (\cM_{H^*}\otimes_{\C}U_q(\gh)) \sharp \C[\cR_Q] 
\en
where $U_q(\gh)$ is the quantum affine algebra in the Drinfeld realization. Furthermore in \cite{JKOS99} Appendix A, we have obtained a realization of 
$U_{q,p}(\gh)$ in terms of the Drinfeld generators in $U_q(\gh)$ 
and a Heisenberg algebra $\C[P_{\al_j}, e^{Q_{\al_j}} (j\in I)]$\footnote{The Heisenberg generators $P_{\al_j}, Q_{\al_j}$ are related to $P_j, Q_j$ in \cite{JKOS99} by $P_{\al_j}=d_jP_j, Q_{\al_j}=-2Q_j$, respectively.}.
Note that such realization is well-defined in the $p$-adic topology. Hence applying the isomorphism between the Drinfeld-Jimbo  realization of $U_q(\gh)$ in terms of the Chevalley generators and the Drinfeld realization of the same algebra interms of the Drinfeld generators, one can expect the isomorphism 
\bea
&&U_{q,p}(\gh)\cong \B_{q,\la(r^*,P)}(\gh)\sharp \C[\cR_Q] \lb{UqpBqla}
\ena
as an associative algebra. Here $\la(r^*,P)$ is given in Appendix A. 
In fact one can derive the same 
$RLL$-relations as \eqref{RLLpm}-\eqref{RLLmp}  by using  
the universal  $R$-matrix of $\Bqla(\gh)$. There the $\cL^\pm(u)$ operators are the elements in $End_\C V\otimes \B_{q,\la(r^*,P)}(\gh)\sharp \C[\cR_Q]$  
and satisfy $\cL^-(u)=\cL^+(u+r^*+c/2)$.  Then by assuming the Gauss decomposition such as \eqref{Lop} in $\cL^\pm(u)$ and denoting their Gauss coordinates by 
$\cE^+_{j,i}(u), \F^+_{i,j}(u), \cK^+_j(u)$
one can show that for $1\leq j\leq N-1$, the 
\bea
&&\cE_j(u):=-\frac{1}{a^*_{j+1,j}[1]^*} \left(\cE^+_{j+1,j}(u-\frac{j}{2}+\frac{c}{2})-\cE^-_{j+1,j}(u-\frac{j}{2})\right), \lb{epem}\\
&&\F_j(u):=\frac{1}{a_{j,j+1}[1]}\left(\F^+_{j+1,j}(u-\frac{j}{2})-\F^-_{j+1,j}(u-\frac{j}{2}+\frac{c}{2})\right) \lb{fpfm}
\ena 
with 
\be
&&\frac{a^*_{j+1,j}a_{j,j+1}[1]}{q-q^{-1}}=1
\en
satisfy the same relations as the elliptic currents $E_{j}(u)$ and  $F_{j}(u)$ in Proposition \ref{defrelfull}\cite{DF,Konno12}. 
Note that the formulas in Definition \ref{basicHC} gives a solution to  \eqref{epem}-\eqref{fpfm}. However we have not yet succeeded to confirm similar formulas for the $j=N$ case 
due to a difficulty of  extracting the relations for these half currents from the $RLL$-relations.  

\section{Hopf Algebroid Structure}
In this section, we introduce an $H$-Hopf algebroid structure into  
the elliptic algebra $\UqpBN$ and formulate it as an elliptic 
quantum group.

\subsection{Definition of the $H$-Hopf algebroid}
Let us recall some basic facts on the $H$-Hopf algebroid following the works of  
Etingof and Varchenko\cite{EV} and of Koelink and Rosengren \cite{KR}.

\begin{dfn}[$H$-bialgebroid]
An $H$-bialgebroid is an $H$-algebra $A$ equipped with two $H$-algebra homomorphisms 
$\Delta:A\to A{{\tot}}A$ (the comultiplication) and $\vep : A\to \cD$ (the counit) such that 
\be
&&(\Delta \tot \id)\circ \Delta=(\id \tot \Delta)\circ \Delta,\\
&&(\vep \tot \id)\circ\Delta =\id =(\id \tot \vep)\circ \Delta,
\en
under the identification \eqref{Diso}.
\end{dfn}
 
\begin{dfn}[$H$-Hopf algebroid]\lb{defS}
An $H$-Hopf algebroid is an $H$-bialgebroid $A$ equipped with a $\C$-linear map $S : A\to A$ (the antipode), such that 
\be
&&S(\mu_r(\hf)a)=S(a)\mu_l(\hf),\quad S(a\mu_l(\hf))=\mu_r(\hf)S(a),\quad \forall a\in A, \hf\in \cM_{\h^*},\\
&&m\circ (\id \tot S)\circ\Delta(a)=\mu_l(\vep(a)1),\quad \forall a\in A,\\
&&m\circ (S\tot\id  )\circ\Delta(a)=\mu_r(T_{\al}(\vep(a)1)),\quad \forall a\in A_{\al\beta},
\en
where $m : A{{\tot}} A \to A$ denotes the multiplication and $\vep(a)1$ is the result of applying the difference operator $\vep(a)$ to the constant function $1\in \cM_{H^*}$.
\end{dfn}

The $H$-algebra $\cD$ is an $H$-Hopf algebroid with 
$\Delta_\cD : \cD\to \cD\tot \cD,\ \vep_\cD: \cD \to \cD,\ 
S_\cD : \cD \to \cD$ defined by 
\be
&&\Delta_\cD(\hf T_{-\al})=\hf T_{-\al} \tot T_{-\al},\\
&&\vep_\cD=\id,
\qquad  S_\cD(\hf T_{-\al})=T_{\al}\hf=(T_{\al}\hf)T_{\al}.
\en

\subsection{The $H$-Hopf algebroid $\UqpBN$}
Now let us consider the $H$-Hopf algebroid structure on $U=\UqpBN$. 
Let us consider the generating function of the $L$-operator matrix elements $L_{i,j}^+(u)$. 
We define two $H$-algebra homomorphisms, the co-unit $\vep : U\to \cD$ and the co-multiplication $\Delta : U\to U \widetilde{\otimes}U$ by
\bea
&&\vep(L^+_{i,j}(u))=\delta_{i,j}{T}_{Q_{\ep_i} }\quad (n\in \Z),
\qquad \vep(e^Q)=e^Q,\lb{counitUqp}\\
&&\vep(\mu_l({\hf}))= \vep(\mu_r(\hf))=\widehat{f}T_0, \lb{counitf}\\
&&\Delta(L^+_{i,j}(u))=\sum_{k}L^+_{i,k}(u)\widetilde{\otimes}
L^+_{k,j}(z),\lb{coproUqp}\\
&&\Delta(e^{Q})=e^{Q}\tot e^{Q},\\
&&\Delta(\mu_l(\hf))=\mu_l(\hf)\widetilde{\otimes} 1,\quad \Delta(\mu_r(\hf))=1\widetilde{\otimes} \mu_r(\hf).\lb{coprof}
\ena
In fact, one can check that $\Delta$ preserves the $RLL$-relations \eqref{RLLpm}-\eqref{RLLmp}. 

\begin{lem}\lb{counitcopro}
The maps $\vep$ and $\Delta$ satisfy
\bea
&&(\Delta\tot \id)\circ \Delta=(\id \tot \Delta)\circ \Delta,\lb{coaso}\\
&&(\vep \tot \id)\circ\Delta =\id =(\id \tot \vep)\circ \Delta.\lb{vepDelta}
\ena
\end{lem}

Combining this with the $H$-algebra structure, the set  $(U, \Delta, \cM_{H^*}, \mu_l, \mu_r, \vep)$ is an $H$-bialgebroid. 

From \eqref{coproUqp}, one can derive the following coproduct formulas for the basic half currents. 
\begin{prop}
\bea
\Delta(K_j^+(u))&=&K^+_{j+1}(u)\tot K^+_{j+1}(u)\left(1+1\tot l^{+'}_{j+1,j+1}(u)+ l^{+'}_{j+1,j+1}(u)\tot 1
+ l^{+'}_{j+1,j+1}(u)\tot l^{+'}_{j+1,j+1}(u)\right)\nn\\
&&\qquad\qquad\qquad\qquad\qquad \times\left(1+\Delta( l^{+'}_{j+1,j+1}(u))\right)^{-1}\nn\\
&=&\left(1+\Delta( l^{+''}_{j+1,j+1}(u))\right)^{-1}\left(1+1\tot l^{+''}_{j+1,j+1}(u)+ l^{+''}_{j+1,j+1}(u)\tot 1
+ l^{+''}_{j+1,j+1}(u)\tot l^{+''}_{j+1,j+1}(u)\right)\nn\\
&&\qquad\qquad\qquad\qquad\qquad \times K^+_{j+1}(u)\tot K^+_{j+1}(u),\nn
\ena
\bea
\Delta(E^+_{j+1,j}(u)) &=& \left(1+\Delta( l^{+'}_{j+1,j+1}(u))\right) \nn \\
&& \times \left(1+1\tot l^{+'}_{j+1,j+1}(u)+ l^{+'}_{j+1,j+1}(u)\tot 1 
 + l^{+'}_{j+1,j+1}(u)\tot l^{+'}_{j+1,j+1}(u)\right)^{-1} \nn\\
&& \times \Biggl(1\tot E^+_{j+1,j}(u)+E^+_{j+1,j}(u)\tot K^+_{j+1}(u)^{-1}K_j^+(u)+E^+_{j+1,j}(u)\tot l^{+'}_{j,j}(u)
\nn\\
&& \qquad + l^{+'}_{j+1,j+1}(u)\tot E^+_{j+1,j}(u)+1\tot l^{+'}_{j+1,j}(u)+ l^{+'}_{j+1,j}(u)\tot  K^+_{j+1}(u)^{-1}K_j^+(u)\nn\\
&& 
\qquad + l^{+'}_{j+1,j}(u)\tot l^{+'}_{j,j}(u)
+ l^{+'}_{j+1,j+1}(u)\tot l^{+'}_{j+1,j}(u)+\sum_{1\preceq k\preceq -1 \atop\not=j,j+1}\hL^{+'}_{j+1,k}(u)\tot
\hL^{+'}_{k,j}(u)\Biggr) \nn\\
&&\qquad  -\Delta(l^{+'}_{j+1,j}(u))\nn\\
\Delta(F^+_{j,j+1}(u))&=& \Biggl(K_j^+(u)K^+_{j+1}(u)^{-1}\tot F^+_{j,j+1}(u)+F^+_{j,j+1}(u)\tot 1+F^+_{j,j+1}(u)\tot l^{+''}_{j+1,j+1}(u)
\nn\\
&&\qquad + l^{+''}_{j,j}(u)\tot F^+_{j,j+1}(u)+ l^{+''}_{j,j+1}(u)\tot 1 + K_j^+(u)K^+_{j+1}(u)^{-1}\tot l^{+''}_{j,j+1}(u)\nn\\
&&
\qquad + l^{+''}_{j,j+1}(u)\tot l^{+''}_{j+1,j+1}(u)
+ l^{+''}_{j,j}(u)\tot l^{+''}_{j,j+1}(u)+\sum_{1\preceq k\preceq -1 \atop \not=j,j+1}\hL^{+''}_{j,k}(u)\tot
\hL^{+''}_{k,j+1}(u)\Biggr)
\nn\\
&& \times\left(1+1\tot l^{+''}_{j+1,j+1}(u)+ l^{+''}_{j+1,j+1}(u)\tot 1
+ l^{+''}_{j+1,j+1}(u)\tot l^{+''}_{j+1,j+1}(u)\right)^{-1} \nn \\
&& \times \left(1+\Delta( l^{+''}_{j+1,j+1}(u))\right) 
-\Delta(l^{+''}_{j,j+1}(u)).\nn
\ena
Here $l^{+'}_{k,l}(u)=K^+_{j+1}(u)^{-1}l^{+}_{k,l}(u),\ l^{+''}_{k,l}(u)=l^{+}_{k,l}(u)K^+_{j+1}(u)^{-1}$ and
\be
&&l^{+}_{j,j}(u)=\sum_{j\prec k}F^+_{j,k}(u)K_k^+(u)E^+_{k,j}(u),\\
&&l^{+}_{i,j}(u)=\sum_{j\prec k}F^+_{i,k}(u)K_k^+(u)E^+_{k,j}(u), \qquad (i\prec j)\\
&&l^{+}_{j,i}(u)=\sum_{j\prec k}F^+_{j,k}(u)K_k^+(u)E^+_{k,i}(u).
\en
\end{prop}

Now let us define formally an algebra antihomomorphism (the antipode) $S : U\to U$ by
\be
&&S(L^+(z))=L^+(z)^{-1},\\
&&S(e^Q)=e^{-Q},\quad S(\mu_r(\hf))=\mu_l(\hf),\quad S(\mu_l(\hf))=\mu_r(\hf).
\en
Then we have the following Lemma.
\begin{lem}\lb{antipode}
The map $S$ satisfies the antipode axioms 
\be
&&m\circ (\id \otimes S)\circ\Delta(x)=\mu_l(\vep(x)1),\quad \forall x\in U,\\
&&m\circ (S\otimes\id  )\circ\Delta(x)=\mu_r(T_{\ha}(\vep(x)1)),\quad \forall x\in (U)_{\ha \hb}.
\en
\end{lem}

From Lemmas \ref{counitcopro} and \ref{antipode}, we have 
\begin{thm}
The $H$-algebra $\UqpBN$ equipped with $(\Delta,\vep,S)$ is an $H$-Hopf algebroid. 
\end{thm}

\begin{dfn}
We call the $H$-Hopf algebroid $(\UqpBN,H,{\cM}_{H^*},\mu_l,\mu_r,\Delta,\vep,S)$ the \\ elliptic quantum group $\UqpBN$. 
\end{dfn}

\section{Representations}
\subsection{Dynamical representations}
Let us consider a vector space $\hV$ over $\FF=\cM_{H^*}$, which is  
${H}$-diagonalizable, i.e.  
\be
&&\hV=\bigoplus_{\la,\mu\in {H}^*}\hV_{\la,\mu},\ \hV_{\la,\mu}=\{ v\in \hV\ |\ q^{P+h}\cdot v=q^{<\la,P+h>} v,\ q^{P}\cdot v=q^{<\mu,P>} v\ \forall 
P+h, P\in 
{H}\}.
\en
Let us define the $H$-algebra $\cD_{H,\hV}$ of the $\C$-linear operators on $\hV$ by
\be
&&\cD_{H,\hV}=\bigoplus_{\al,\beta\in {H}^*}(\cD_{H,\hV})_{\al,\beta},\\
&&\hspace*{-10mm}(\cD_{H,\hV})_{\al,\beta}=
\left\{\ X\in \End_{\C}\hV\ \left|\ 
\mmatrix{ f(P+h)X=X f(P+h+<\alpha,P+h>),\cr 
f(P)X=X f(P+<\beta,P>)\cr
 f(P), f(P+h)\in \FF,\ X\cdot\hV_{\la,\mu}\subseteq 
 \hV_{\la+\al,\mu+\beta}\cr}  
 \right.\right\},\\
&&\mu_l^{\cD_{H,\hV}}(\widehat{f})v=f(<\la,P+h>,p)v,\quad 
\mu_r^{\cD_{H,\hV}}(\widehat{f})v=f(<\mu,P>,p^*)v,\quad \widehat{f}\in \FF,
\ v\in \hV_{\la,\mu}.
\en
\begin{dfn}\cite{EV,KR,Konno09} 
We define a dynamical representation of $U_{q,p}(B_N^{(1)})$ on $\hV$ to be  
 an $H$-algebra homomorphism ${\pi}: U_{q,p}(B_N^{(1)}) 
 \to \cD_{H,\hV}$. By the action of $U_{q,p}(B_N^{(1)})$ we regard $\hV$ as a 
$U_{q,p}(B_N^{(1)})$-module. 
\end{dfn}

\begin{dfn}
For $k\in \C$, we say that a $U_{q,p}(B_N^{(1)})$-module has  level $k$ if $c$ act 
as the scalar $k$ on it.  
\end{dfn}


\begin{dfn}
Let ${\cal H}$, ${\cal N}_+, {\cal N}_-$ be the subalgebras of 
$U_{q,p}(B_N^{(1)})$ 
generated by 
$c, d, 
K^\pm_{i}\ (i\in I)$, by $\al^\vee_{i,n}\ (i\in I, n\in \Z_{>0})$,  $e_{i, n}\ (i\in I, n\in \Z_{\geq 0})$  
$f_{i, n}\ (i\in I, n\in \Z_{>0})$ and by $\al^\vee_{i,-n}\ (i\in I, n\in \Z_{>0}),\ e_{i, -n}\ (i\in I, n\in \Z_{> 0}),\ f_{i, -n}\ (i\in I, n\in \Z_{\geq 0})$, respectively.   
\end{dfn}

\begin{dfn}
For $k\in\C$, $\la\in \h^*$ and $\mu\in H^*$, 
a (dynamical) $U_{q,p}(B_N^{(1)})$-module $\hV(\la,\mu)$ is called the 
level-$k$ highest weight module with the highest weight $(\la,\mu)$, if there exists a vector 
$v\in \hV(\la,\mu)$ such that
\be
&&\hV(\la,\mu)=U_{q,p}(B_N^{(1)})\cdot v,\qquad \cN_+\cdot v=0,\\
&&c\cdot v=kv, 
\quad  f({P})\cdot v =f({<\mu,P>})v,\quad f({P+h})\cdot v =f({<\la,P+h>})v.
\en
\end{dfn}

\subsection{Finite dimensional dynamical representation}

We here give an elliptic and dynamical analogue of the evaluation representation associated with the vector representation of $B_N^{(1)}$. Let us consider   $\displaystyle{\hV=\bigoplus_{1 \preceq m \preceq -1} \FF v_m\otimes 1}$ and $\hV_z=\hV[z,z^{-1}]$. Here $e^{Q_{\al}}\in \C[\cR_Q]$ acts on $f(P_\beta)v\otimes 1$  by $e^{Q_{\al}}(f(P_\beta)v\otimes 1)=f(P_\beta-(\al,\beta))v\otimes 1$. 
\begin{thm}\lb{vecRep}
Let $E_{j,k}\ (1\preceq j,k\preceq -1)$ denote the matrix units such that 
$E_{j,k}v_l=\delta_{k,l}v_j$. 
The following gives the $(2N+1)$-dimensional dynamical representation of $U_{q,p}(B_N^{(1)})$ 
on $\hV_z$.  
\be
 && \pi_z(e_j(w))=\left(\frac{(pq^2;p)_{\infty}}{(p;p)_{\infty}}E_{j,j+1}
     \delta\left(q^j\frac{z}{w}\right) 
    +\frac{(pq^2;p)_{\infty}}{(p;p)_{\infty}}E_{-j-1,-j}
    \delta\left(q^{-j}\xi\frac{z}{w}\right)\right)e^{-Q_{\al_j}}, \\
 && \pi_z(f_j(w))=\frac{(pq^{-2};p)_{\infty}}{(p;p)_{\infty}}E_{j+1,j}
     \delta\left(q^j\frac{z}{w}\right) 
    +\frac{(pq^{-2};p)_{\infty}}{(p;p)_{\infty}}E_{-j,-j-1}
    \delta\left(q^{-j}\xi\frac{z}{w}\right), \\
 && \pi_z(\psi_j^{-}(w,p))=q^{\pi(h_j)}e^{-Q_{\al_j}}
     \frac{\Theta_p(q^{j-2h_j^{+}} \frac{z}{w})}{\Theta_p(q^{j} \frac{z}{w})}
      \frac{\Theta_p(q^{-j-2h_j^{-}}\xi \frac{z}{w})}
                      {\Theta_p(q^{-j}\xi \frac{z}{w})}, \\
 && \pi_z(\psi_j^{+}(w,p))=q^{-\pi(h_j)}e^{-Q_{\al_j}}
     \frac{\Theta_p(q^{-j+2h_j^{+}} \frac{w}{z})}{\Theta_p(q^{-j} \frac{w}{z})}
     \frac{\Theta_p(q^{j+2h_j^{-}} \xi^{-1}\frac{w}{z})}
                      {\Theta_p(q^{j}\xi^{-1} \frac{w}{z})},\qquad (1\leq j\leq N-1) \\
 && \pi_z(e_N(w))=\frac{(pq^2;p)_{\infty}}{(p;p)_{\infty}} [2]_N E_{N,0}
     \delta\left(q^N \frac{z}{w}\right) 
    +\frac{(pq;p)_{\infty} (pq^{-2};p)_{\infty}}
          {(p;p)_{\infty}(pq^{-1};p)_{\infty}}E_{0,-N}
    \delta\left(q^{N-1}\frac{z}{w}\right), \\
 && \pi_z(f_N(w))=\frac{(pq^2;p)_{\infty}(pq^{-1};p)_{\infty}}
                       {(pq;p)_{\infty}(p;p)_{\infty}}[2]_N E_{-N,0}
     \delta\left(q^{N-1}\frac{z}{w}\right) 
    +\frac{(pq^{-2};p)_{\infty}}{(p;p)_{\infty}}E_{0,N}
    \delta\left(q^N \frac{z}{w}\right), \\
 && \pi_z(\psi_N^{-}(w))=q_N^{\pi(h_N)}e^{-Q_{\al_N}}
     \frac{\Theta_p(q^{N-\frac{1}{2}-\frac{3}{2}h_l^{+}} \frac{z}{w})}{\Theta_p(q^{N-\frac{1}{2}+\frac{1}{2}h_N^+} \frac{z}{w})}
      \frac{\Theta_p(q^{N-\frac{1}{2}-\frac{3}{2}h_N^{-}}\frac{z}{w})}
                      {\Theta_p(q^{N-\frac{1}{2}+\frac{1}{2}h_N^-} \frac{z}{w})},\\
 && \pi_z(\psi_N^{+}(w))=q_N^{-\pi(h_N)}e^{-Q_{\al_N}}
     \frac{\Theta_p(q^{-(N-\frac{1}{2})+\frac{3}{2}h_N^{+}} \frac{w}{z})}{\Theta_p(q^{-(N-\frac{1}{2})-\frac{1}{2}h^+_N} \frac{w}{z})}
      \frac{\Theta_p(q^{-(N-\frac{1}{2})+\frac{3}{2}h_N^{-}} \frac{w}{z})}
                      {\Theta_p(q^{-(N-\frac{1}{2})-\frac{1}{2}h^-_N} \frac{w}{z})}.
\en
Here $\pi(h_j)=h_j^++h_j^-$, $h_j^{+}=E_{j,j}-E_{j+1,j+1},\,h_j^{-}=E_{-j-1,-j-1}-E_{-j,-j}$, $\pi(h_N)=2(E_{N,N}-E_{-N,-N}),\ h_N^+=E_{N,N}-E_{0,0},\ h_N^-=E_{0,0}-E_{-N,-N}$. 
\end{thm}

\begin{thm}\lb{vecRepHC}
In terms of the half currents the dynamical representation $(\pi_z, \hV_z)$ is given as follows.
\be
\pi_{z}(K^+_{+j}(v))&=&\rho^+(v-u)\left\{\frac{[v-u]}{[v-u+1]}\sum_{1\preceq k\preceq j-1}E_{k,k}+E_{j,j}+\frac{[v-u-1]}{[v-u]}\sum_{j+1\preceq k\not=-j \preceq -1}E_{k,k}\right.\\
&&\left.\qquad+\frac{[v-u-1][v-u+j+\eta-1]}{[v-u][v-u+j+\eta]}E_{-j,-j}\right\}e^{-Q_{\vep_j}},\\
\pi_{z}(K^+_{-j}(v))&=&\rho^+(v-u)\left\{\frac{[v-u]}{[v-u+1]}\sum_{1\preceq k\not=j\preceq -(j+1)}E_{k,k}+\frac{[v-u][v-u-j-\eta]}{[v-u+1][v-u-j+\frac{2N-1}{2}+1]}E_{j,j}\right.\\
&&\left.\qquad+\frac{[v-u-1]}{[v-u]}\sum_{-(j-1)\preceq k \preceq -1}E_{k,k}+E_{-j,-j}
\right\}e^{Q_{\vep_j}},\\
\pi_{z}(K^+_{0}(v))&=&\rho^+(v-u)\left\{\frac{[v-u]}{[v-u+1]}\sum_{1\preceq k\preceq N}E_{k,k}+\frac{[v-u-1]}{[v-u]}\sum_{1\preceq k\preceq N}E_{-k,-k}\right.\\
&&\qquad\qquad \left.+\frac{[v-u+\frac{1}{2}][v-u-1]}{[v-u-\frac{1}{2}][v-u+1]}E_{0,0}\right\}.
\en
For $1\preceq j\prec l\preceq N+1\equiv 0$, 
\be
&&\hspace{-1.5cm}\pi_{z}(E^+_{l,j}(v))= e^{Q_{\vep_l}}\left\{-E_{j,l}\frac{[v-u-P_{j,l}][1]}{[v-u][P_{j,l}]}
+E_{-l,-j}\frac{[v-u+l-1+\eta-P_{j,l}][1]}{[v-u+l-1+\eta][P_{j,l}]}\prod_{m=j+1}^{l-1}\frac{[P_{j,m}+1]}{[P_{j,m}]}\right\}e^{-Q_{\vep_j}},\\
&&\hspace{-1cm}\pi_{z}(F^+_{j,l}(v))=
E_{1,j}\frac{[v-u+P_{j,l}][1]}{[v-u][P_{j,l}]}-E_{-j,-l}
\frac{[v-u+l-1+\eta+P_{j,l}][1]}{[v-u+l-1+\eta][P_{j,l}]}\prod_{m=j+1}^{l-1}\frac{[P_{j,m}-1]}{[P_{j,m}]},
\en
\be
&&\hspace{-1cm}\pi_{z}(E^+_{-j,-l}(v))\\
&&\hspace{-1cm}\qquad= e^{-Q_{\vep_j}}\left\{-E_{-l,-j}\frac{[v-u-P_{-l,-j}][1]}{[v-u][P_{-l,-j}]}
+E_{j,l}\frac{[v-u-j-\eta-P_{-l,-j}][1]}{[v-u-j-\eta][P_{-l,-j}]}\prod_{m=j+1}^{l-1}\frac{[P_{-l,-m}+1]}{[P_{-l,-m}]}\right\}e^{Q_{\vep_l}},\\
&&\hspace{-1cm}\pi_{z}(F^+_{-l,-j}(v))\\
&&\hspace{-1cm}\qquad =
E_{-j,-l}\frac{[v-u+P_{-l,-j}][1]}{[v-u][P_{-l,-j}]}
-E_{l,j}\frac{[v-u-j-\eta+P_{-l,-j}+\delta_{l,N+1}][1]}{[v-u-j-\eta][P_{-l,-j}+\delta_{l,N+1}]}\prod_{m=j+1}^{l-1}\frac{[P_{-l,-m}-1+\delta_{l,N+1}]}{[P_{-l,-m}+\delta_{l,N+1}]}
\en
\noindent
i) $j \prec k \prec N$,
\be
&&\pi_{z}(E_{-k,j}^+(v))\\
&&=e^{-Q_{\vep_k}}\left\{-E_{j,-k}\frac{[v-u-P_{j,-k}][1]}{[v-u][P_{j,-k}]}
\right.\\
&&\left.
+E_{k,-j}\frac{[v-u-k-\eta-P_{j,-k}][1]}{[v-u-k-\eta][P_{j,-k}]}
\frac{[P_j+1]}{[P_j]}\prod_{m=k+1}^{N}\frac{[P_{j,-m}+1]}{[P_{j,-m}]}
\prod_{m=j+1}^{N}\frac{[P_{j,m}+1]}{[P_{j,m}]}\right\}e^{-Q_{\vep_j}},\\
&&\pi_{z}(F_{j,-k}^+(v))\\
&&=E_{-k,j}\frac{[v-u+P_{j,-k}][1]}{[v-u][P_{j,-k}]}\\
&&
-E_{-j,k}\frac{[v-u-k-\eta+P_{j,-k}][1]}{[v-u-k-\eta][P_{j,-k}]}
\frac{[P_j-1]}{[P_j]}\prod_{m=k+1}^{N}\frac{[P_{j,-m}-1]}{[P_{j,-m}]}
\prod_{m=j+1}^{N}\frac{[P_{j,m}-1]}{[P_{j,m}]}.
\en

\noindent
ii) $j = k \prec N$,
\be
&&\pi_{z}(E_{-j,j}^+(v))\\
&&=e^{-Q_{\vep_j}}E_{j,-j}\left\{-\frac{[v-u-2P_{j}-1][1][v-u-j-\eta+1]}{[v-u][2P_{j}+1][v-u-j-\eta]}\right.\\
&&\qquad\qquad\qquad \left.+G_{P_j}\frac{[v-u-2P_{j}-j-\la][1]}{[v-u-j-\la][2P_{j}+1]}
\prod_{m=1}^{j-1}\left(\frac{[P_{j,-m}]}{[P_{j,-m}+1]}\frac{[P_{j,m}]}{[P_{j,m}+1]}
\right)\right\}e^{-Q_{\vep_j}},\\
&&\pi_{z}(F_{j,-j}^+(v))\\
&&=E_{-j,j}\left\{\frac{[v-u+2P_{j}-1][1][v-u-j-\eta+1]}{[v-u][2P_{j}-1][v-u-j-\eta]}\right.\\
&&\qquad\qquad\qquad \left.-G_{-P_j}\frac{[v-u+2P_{j}-j-\eta][1]}{[v-u-j-\eta][2P_{j}-1]}
\prod_{m=1}^{j-1}\left(\frac{[P_{j,-m}]}{[P_{j,-m}-1]}\frac{[P_{j,m}]}{[P_{j,m}-1]}
\right)\right\}.
\en

\noindent
iii) $k \prec j \preceq N+1\equiv 0$, 
\be
&&\pi_{z}(E_{-k,j}^+(v))\\
&&=e^{-Q_{\vep_k}}\left\{-E_{j,-k}\frac{[v-u-P_{j,-k}][1]}{[v-u][P_{j,-k}]}
\right.\\
&&\left.
+E_{k,-j}\frac{[v-u-k-\eta-P_{j,-k}][1]}{[v-u-k-\eta][P_{j,-k}]}
\frac{[P_j+1]}{[P_j]}\prod_{m=k+1\atop \not=j}^{N}\frac{[P_{j,-m}+1]}{[P_{j,-m}]}
\prod_{m=j+1}^{N}\frac{[P_{j,m}+1]}{[P_{j,m}]}\right\}e^{-Q_{\vep_j}},
\en
\be
&&\pi_{z}(F_{j,-k}^+(v))\\
&&=E_{-k,j}\frac{[v-u+P_{j,-k}][1]}{[v-u][P_{j,-k}]}\\
&&
-E_{-j,k}\frac{[v-u-k-\eta+P_{j,-k}][1]}{[v-u-k-\eta][P_{j,-k}]}
\frac{[P_j-1]}{[P_j]}\prod_{m=k+1\atop \not=j}^{N}\frac{[P_{j,-m}-1]}{[P_{j,-m}]}
\prod_{m=j+1}^{N}\frac{[P_{j,m}-1]}{[P_{j,m}]}.
\en

\end{thm}

In addition, it is also worth to remark the following formulas.
\be
\pi_{z}(H_j^{+}(v))&=& \Big\{ 
 \frac{[v-u-\frac{j}{2}-1]}{[v-u-\frac{j}{2}]}E_{j,j}
 +
  \frac{[v-u-\frac{j}{2}-1]}{[v-u-\frac{j}{2}]}E_{j+1,j+1}\\
&& \quad +\frac{[v-u+\eta+\frac{j}{2}+1]}
{[v-u+\eta+\frac{j}{2}]}E_{-j-1,-j-1}
+ \frac{[v-u+\eta+\frac{j}{2}-1]}{[v-u+\eta+\frac{j}{2}]}
E_{-j,-j}\Big\}e^{-Q_{\alpha_j}} \nn \\
&& \quad (1 \leqq j \leqq N-1), \\
 \pi_{z}(H_N^{+}(v))&=& \Big\{
   \frac{[v-u-\frac{N}{2}+1]}{[v-u-\frac{N}{2}]}E_{N,N}
+ 
\frac{[v-u-\frac{N}{2}-1][v-u-\frac{N}{2}-\frac{1}{2}]}{[v-u-\frac{N}{2}][v-u-\frac{N}{2}+\frac{1}{2}]}E_{0,0}
\\&& \quad\qquad\qquad 
+\frac{[v-u-\frac{N}{2}-\frac{1}{2}]}{[v-u-\frac{N}{2}+\frac{1}{2}]}E_{-N,-N} \Big\} e^{-Q_{\alpha_N}},\\
\pi_{z}(H_j^{-}(v))&=&\Big\{
\frac{[u-v+\frac{j}{2}-1]}{[u-v+\frac{j}{2}]}E_{j,j}
+ \frac{[u-v+\frac{j}{2}+1]}{[u-v+\frac{j}{2}]}E_{j+1,j+1}\\
&& \quad +\frac{[u-v-\eta-\frac{j}{2}-1]}
{[u-v-\eta-\frac{j}{2}]}E_{-j-1,-j-1}
+\frac{[u-v-\eta-\frac{j}{2}+1]}{[u-v-\eta-\frac{j}{2}]}
E_{-j,-j}\Big\}e^{-Q_{\alpha_j}} \nn \\
&& \quad (1 \leqq j \leqq N-1), \\
\pi_{z}(H_N^{-}(v))&=& \Big\{
    \frac{[u-v+\frac{N}{2}-1]}{[u-v+\frac{N}{2}]}E_{N,N}
+
\frac{[u-v+\frac{N}{2}+\frac{1}{2}][u-v+\frac{N}{2}-1]}{[u-v+\frac{N}{2}-\frac{1}{2}][u-v+\frac{N}{2}]}E_{0,0} \\
&&  \quad\qquad\qquad  +\frac{[u-v+\frac{N}{2}+\frac{1}{2}]}{[u-v+\frac{N}{2}-\frac{1}{2}]}E_{-N,-N} \Big\}e^{-Q_{\alpha_N}}.
\en
\noindent
{\it Remark.}  The statements in this theorem and the next one remain unchanged when one uses $\rho^+_0(u)$ and $\rho_0^{+*}(u)$ instead of  $\rho^+(u)$ and $\rho^{+*}(u)$, respectively. See Sec.\ref{VO}

Combining the formulas in Theorem \ref{vecRep}, it is not so hard to show the following. 
\begin{cor}\lb{repLR}
\be
&&\pi_z({\hL^+_{i,j}(v)})_{k,l}=R^+(v-u,P)_{ik}^{jl}.
\en
\end{cor}
\noindent
{\it Proof.} For example, for $-N\preceq -j\preceq -1$, we obtain 
\be 
\pi_{z}(L^+_{-j,-j}(v))&=&\pi_{z}\left(K^+_{-j}(v)+\sum_{-j\preceq -k\preceq -1}F^+_{-j,-k}(v)
K^+_{-k}(v)E^+_{-k,-j}(v)\right)\\
&=&\rho^+(v-u)\left\{\bar{b}(v-u)\sum_{1\preceq k\not=j\preceq -(j+1)}E_{k,k}+
\sum_{-(j-1)\preceq -k\preceq -1}b(v-u,P_{-j,-k})E_{-k,-k}\right.\\
&&\left.\qquad\qquad \qquad +\bar{d}(v-u,P_j,P_{-j})E_{j,j}
+E_{-j,-j}\mmatrix{\ \cr \ \cr }\right\}e^{Q_{\vep_j}}.
\en
Here we used the identity 
\be
&&\bar{d}(u,P_j,P_{-j})=\frac{[u]}{[u+1]}\left(\frac{[u-j-\eta]}{[u-j+1-\eta]}\right.\\
&&\quad \left. -\sum_{-(j-1)\preceq -k\preceq -1}\frac{[u-k-\eta+P_{-j,-k}][u-k-\eta-P_{-j,-k}][1]^2}{[u-k-\eta][u-k+1-\eta][P_{-j,-k}]^2} 
\prod_{m=k+1}^{j-1}\frac{[P_{-j,-m}-1][P_{-j,-m}+1]}{[P_{-j,-m}]^2}\right).
\en
In addition, we have
\be 
\pi_{z}(L^+_{0,0}(v))&=&\pi_{z}\left(K^+_{0}(v)+\sum_{-N\preceq -k\preceq -1}F^+_{0,-k}(v)
K^+_{-k}(v)E^+_{-k,0}(v)\right)\\
&=&\rho^+(v-u)\left\{\bar{b}(v-u)\sum_{1\preceq k\preceq N}E_{k,k}+e_0(v-u,P)E_{0,0}
\right.\\
&&\qquad\qquad\left. +\sum_{-N\preceq -k\preceq -1}b(v-u,P_{0,-k})E_{-k,-k} \right\}
e^{Q_{\vep_j}}
\en
This is due to the identity
\be
&&e_0(u,P)=\frac{[u-1][u+\frac{1}{2}]}{[u+1][u-\frac{1}{2}]}
-\frac{[u]}{[u+1]}\sum_{-N\preceq -k\preceq -1}
\frac{[u-k-\eta+\frac{1}{2}+P_{k}][u-k-\eta+\frac{1}{2}-P_{k}][1]^2}{[u-k-\eta][u-k+1-\eta][P_k+\frac{1}{2}][P_k-\frac{1}{2}]}.
\en
Furthermore, the coefficient of $E_{1, -1}$  in 
\be
\pi_{z}(\hL^+_{-1,1}(v))&=&\pi_{z}(K_{-1}^+(v)E^+_{-1,1}(v))\\
&=&\rho^+(v-u)E_{1, -1}\left\{
-\frac{[1][v-u-2P_1-1]}{[v-u+1][2P_1+1]}
+G_{P_1}\frac{[v-u][v-u-2P_1-1-\la][1]}{[v-u+1][v-u-\la][2P_1+1]}
\right\}e^{-Q_{\vep_1}}
\en
coincides with $e_1(v-u,P)$, and for $k \prec j \preceq N$ the coefficient of 
$E_{k,-j}$ in 
\be
\pi_{z}(L^+_{-k,j}(v))&=&\pi_{z}\left(K^+_{-k}(v)E^+_{-k,j}(v)+\sum_{-(k-1)\preceq -l\preceq -1}F^+_{-k,-l}(v)
K^+_{-l}(v)E^+_{-l,j}(v)\right)\\
&=&\rho^+(v-u)\left\{\bar{c}(v-u,P_{j,-k})E_{j,-k}\mmatrix{\quad\quad \cr \quad\cr}\right.\\
&&\left. +E_{k,-j}G_{P_j}\frac{[v-u][v-u-1-\eta-P_{j,-k}][1]}{[v-u+1][v-u-\eta][P_{j,-k}+1]}\prod_{m=1}^{j-1}\frac{[P_{j,m}]}{[P_{j,m}+1]}
\prod_{m=1}^{k-1}\frac{[P_{-k,-m}-1]}{[P_{-k,-m}]}
\right\}e^{-Q_{\vep_j}}
\en
coincides with  $\bar{d}(v-u,P_k,P_j)$.
\qed

\subsection{The level-1 representation}
Next we consider level-1 representation of $U_{q,p}(B_N^{(1)})$.
We follow the work \cite{FKO}.  Let  $e^{\alpha_i}\ (i\in I)$ be the generators of the group algebra $\C[\cQ]$ with the following central extension. 
 \be
&&e^{ \alpha_i}e^{ \alpha_j}=(-1)^{(\alpha_i, \alpha_j)+(\alpha_i, \alpha_i)(\alpha_j, \alpha_j)} e^{\alpha_j}e^{ \alpha_i}
\en
Let us consider the  Neveu-Schwartz ($NS$) fermion
$\{\Psi_n | n \in {\mathbb Z}+\frac{1}{2}\}$
and the Ramond ($R$) fermion
$\{\Psi_n | n \in {\mathbb Z}\}$ 
satisfying the following anti-commutation relations.
\begin{eqnarray*}
~\{\Psi_m,\Psi_n\}=\delta_{m+n,0}\cN({q^{m}+q^{-m}})
\end{eqnarray*}
with $\cN=1/(q^{\frac{1}{2}}+q^{-\frac{1}{2}})$. 
We  define 
\be
&&{\cal F}^{NS}={\mathbb C}[\Psi_{-\frac{1}{2}},\Psi_{-\frac{3}{2}} 
\cdots],\qquad 
\widetilde{\F}^{R}=\C[\Psi_{-1},\Psi_{-2},...]
\en
and their submodules $\F^{NS, R}_{even}$ (reps. $\F^{NS, R}_{odd}$)  
generated by the even (reps. odd) number of $\Psi_{-m}$'s.  
Due to the zero-mode $\Psi_0$ we have two degenerate vacuum states $1$ and $\Psi_01$. 
We hence consider the extended space 
\be
&&\widehat{\F}^R=\widetilde{\F}^R\otimes \C^2
\en
and realize the  $R$-fermions by   
\be
&&\widehat{\Psi}_m=\Psi_m\otimes \left(\mmatrix{1&0\cr 0&-1\cr}\right) 
\quad (m\in \Z_{\not =0} )
,\qquad
\widehat{\Psi}_0=\cN^{\frac{1}{2}}(1\otimes \left(\mmatrix{0&1\cr 1&0\cr}\right)).  
\en
Note that $\{\hPsi_m,\hPsi_n\}=\delta_{m+n,0}\cN({q^{m}+q^{-m}})$. We set 
\be
&&\F^R=\F^R_{even}
\otimes \C\left(\mmatrix{1\cr 1\cr}\right)
\oplus \F^R_{odd}\otimes\C\left(\mmatrix{1\cr -1\cr}\right).\lb{RFock}
\en
The action of $\Psi_{m}$ on $\F^{NS}$ is given by
\be
&&\Psi_{-m}\cdot u= \Psi_{-m}u,\qquad 
\Psi_{m}\cdot u=\{\Psi_{m}, u\}\qquad (m\in\frac{1}{2}+\Z_{\geq 0}),
\en
where $u\in\F^{NS}$, whereas $\widehat{\Psi}_{m}$ acts on $\F^{R}$ as
\be
&&\widehat{\Psi}_{-m}\cdot u\otimes v= \Psi_{-m}u\otimes \left(\mmatrix{1&0\cr 0&-1\cr}\right)v \quad (m\in\Z_{>0}),
\quad \widehat{\Psi}_{0}\cdot u\otimes v= u\otimes \left(\mmatrix{0&1\cr 1&0\cr}\right)v,\nn\\
&&\widehat{\Psi}_{m}\cdot u\otimes v=\{\Psi_{m}, u\}\otimes \left(\mmatrix{1&0\cr 0&-1\cr}\right)v\quad (m\in\Z_{>0}),
\en 
where $u\in\widetilde{\F}^{R},\ v\in \C^2$.  We define the fermion fields    
$\Psi^{NS}(z)$ and $\Psi^{R}(z)$ by
\begin{eqnarray*}
&&\Psi^{NS}(z)=\sum_{n \in {\mathbb{Z}}+\frac{1}{2}}
\Psi_n z^{-n},\qquad
\Psi^{R}(z)=\sum_{n \in {\mathbb{Z}}}
\widehat{\Psi}_n z^{-n}.
\end{eqnarray*}
Then we have the  operator product expansions 
\be
&&\Psi(z)\Psi(w)=:\Psi(z)\Psi(w):+<\Psi(z)\Psi(w)>,\\
&&<\Psi(z)\Psi(w)>=\left\{
\mmatrix{
    \frac{(zw)^{1/2}(z-w)}{(z-qw)(z-q^{-1}w)} & \mbox{for NS}\cr 
    \cN\frac{(z-w)(z+w)}{(z-qw)(z-q^{-1}w)} & \mbox{for R}.\cr
}
\right.
\en

Now we define 
\be
&&W(\Lambda_0)= \F^{NS}_{even}\otimes\C[\cQ_0]\oplus \F^{NS}_{odd}\otimes \C[\cQ_0]e^{\bar{\Lambda}_1},\\
&&W(\Lambda_1)=\F^{NS}_{even}\otimes \C[\cQ_0]e^{\bar{\Lambda}_1}\oplus \F^{NS}_{odd}\otimes \C[\cQ_0],\\
&&W(\Lambda_N)=\F^R\otimes \C[\cQ]e^{\bar{\Lambda}_N}
\cong \F^R\otimes \C[\cQ_0]e^{\bar{\Lambda}_N}\oplus \F^R\otimes \C[\cQ_0] e^{\bar{\Lambda}_1+\bar{\Lambda}_N},
\en
where $\cQ_0$ denotes the sublattice of $\cQ$ generated by the long roots. 
For generic $\mu\in \h^*$ and $a=0,1,N$, we set 
\be
\hV(\Lambda_{a}+\mu,\mu)&=&\F_{\al,1}\otimes_\C (\FF\otimes_\C W(\Lambda_a))\otimes e^{Q_{\bar{\mu}}}\C[\cR_Q].
\en 
Then we have the following decomposition.
\be
\hV(\Lambda_{a}+\mu,\mu)&=&\bigoplus_{\gamma\in \cQ_0,\kappa\in \cQ}\bigoplus_{\la\in {max}(\Lambda_a)\atop mod\ Q_0+\C\delta} 
\F_{\la,\gamma,\kappa}(\Lambda_a,\mu),
\en
where 
\be
&&\F_{\Lambda_0,\gamma,\kappa}(\Lambda_0,\mu)=\FF\otimes_\C(\F_{\al,1}\otimes \F^{NS}_{even}\otimes e^\gamma)\otimes e^{Q_{\bar{\mu}+\kappa}},\\
&&\F_{\Lambda_1,\gamma,\kappa}(\Lambda_0,\mu)=\FF\otimes_\C(\F_{\al,1}\otimes\F^{NS}_{odd}\otimes e^{\bar{\Lambda}_1+\gamma})\otimes e^{Q_{\bar{\mu}+\kappa}},\\
&&\F_{\Lambda_1,\gamma,\kappa}(\Lambda_1,\mu)=\FF\otimes_\C(\F_{\al,1}\otimes\F^{NS}_{even}\otimes e^{\bar{\Lambda}_1+\gamma})\otimes e^{Q_{\bar{\mu}+\kappa}},\\
&&\F_{\Lambda_0,\gamma,\kappa}(\Lambda_1,\mu)=\FF\otimes_\C(\F_{\al,1}\otimes\F^{NS}_{odd}\otimes e^{\gamma})\otimes e^{Q_{\bar{\mu}+\kappa}},\\
&&\F_{\Lambda_N,\gamma,\kappa}(\Lambda_N,\mu)=\FF\otimes_\C(\F_{\al,1}\otimes\F^{R}\otimes e^{\bar{\Lambda}_N+\gamma})\otimes e^{Q_{\bar{\mu}+\kappa}},\\
&&\F_{\Lambda_N-\al_N,\gamma,\kappa}(\Lambda_1,\mu)=\FF\otimes_\C(\F_{\al,1}\otimes\F^{R}\otimes e^{\bar{\Lambda}_N+\bar{\Lambda}_1+\gamma})\otimes e^{Q_{\bar{\mu}+\kappa}}.
\en

\begin{thm}\lb{levelone}\cite{FKO} 
The three spaces $\hV(\Lambda_{a}+\mu,\mu)$ ($a=0,1,N$) give the level-1 irreducible $U_{q,p}(\hat{B}_{l}^{(1)})$-modules with the higest weight $(\Lambda_a+\mu,\mu)$, where 
the highest weight vectors are given by $1\otimes  1\otimes 1\otimes e^{Q_{\bar{\mu}}}$ for $\hV(\Lambda_{0}+\mu,\mu)$, $1\otimes  1\otimes e^{\bar{\Lambda}_1}\otimes e^{Q_{\bar{\mu}}}$ for $\hV(\Lambda_{1}+\mu,\mu)$ and $1\otimes  1\otimes \left(\mmatrix{1\cr 1\cr}\right)\otimes e^{\bar{\Lambda}_N}\otimes e^{Q_{\bar{\mu}}}$ for $\hV(\Lambda_{N}+\mu,\mu)$,   
respectively. 
The  action of the elliptic currents on $\hV(\Lambda_{a}+\mu,\mu)$ ($a=0,1,N$) is given as follows. 
\begin{eqnarray}
 && E_j(v) =\; :\exp \left\{ -\sum_{n \ne 0} \frac{1}{[n]_q} \alpha_{j,n}w^{-n} \right\}:
 e^{\alpha_j}w^{h_{\alpha_j}}e^{-Q_{\alpha_j}}w^{-\frac{P_{\alpha_j}-1}{r^*}},  \\
 && E_N(v) =\; \frac{1}{\cN^{\frac{1}{2}}}:\exp \left\{ -\sum_{n \ne 0} \frac{1}{[n]_q} \alpha_{N,n}w^{-n} \right\}:\Psi(w)
 e^{\alpha_N}w^{h_{\alpha_N}+\frac{1}{2}}e^{-Q_{\alpha_N}}w^{-\frac{P_{\alpha_N}-1/2}{r^*}}, \\
 && F_j(v) =\; :\exp \left\{ \sum_{n \ne 0}\frac{1}{[n]_q} \frac{1-p^{*n}}{1-p^n}\alpha_{j,n}(q^{-1}w)^{-n}\right\}: e^{-\alpha_j}w^{-h_{\alpha_j}}
w^{\frac{P_{\alpha_j}+h_{\alpha_j}-1}{r}}, \\
 && F_N(v) =\; \frac{1}{\cN^{\frac{1}{2}}}:\exp \left\{ \sum_{n \ne 0} \frac{1}{[n]_q}\frac{1-p^{*n}}{1-p^n} \alpha_{N,n}(q^{-1}w)^{-n} \right\}:\Psi(w)
 e^{-\alpha_N}w^{-h_{\alpha_N}+\frac{1}{2}}w^{\frac{P_{\alpha_N}+h_{\alpha_N}-1/2}{r}}, \nn\\
\end{eqnarray}
$w=q^{2v}$, $(1 \le j \le N-1)$ together with $H_j^{\pm}(v), \,K_{j}^{+}(v)$  in Sec. \ref{MEC}. 
\end{thm}

\section{Vertex Operators of $U_{q,p}(B_N^{(1)})$}\lb{VO}
In this section we discuss the type I and II vertex operators of the $U_{q,p}(B_N^{(1)})$-modules. Through this section, we use $\rho^+_0(u)$ as the prefactor of the  $R$-matrix i.e. 
\be
&&R^+(u,s)=\rho^+_0(u)\bar{R}^+(u,s).
\en
In addition, we often use the following component form of the $RLL$-relation \eqref{RLLpm}.
\bea
&&\sum_{i',j'}R^+(u,P+h)_{i,j}^{i'j'}\hL^{+}_{i',i''}(u_1)\hL^+_{j',j''}(u_2) \nn \\
&& \qquad =\sum_{i',j'}\hL^+_{j,j'}(u_2)\hL^+_{i,i'}(u_1)R^{*+}(u,P-(\pi(h))_{i',i''}-(\pi(h))_{j',j''})_{i'j'}^{i''j''}.
\lb{RLLcomp}
\ena 
Here the components of the $R^{*+}$-matrix is evaluated in the same way as \eqref{ellR}. 
For example, the $(j_1,j_2), (j_1,j_2)\ (j_1\prec j_2)$ component is given by 
\be
&&b^*(u,P_{j_1,j_2}-(\pi(h_{j_1,j_2}))_{j_1,j_1}-(\pi(h_{j_1,j_2}))_{j_2,j_2}), 
\en
where $P_{j_1,j_2}=P_{\ep_{j_1}}-P_{\ep_{j_2}}$, $h_{j_1,j_2}=h_{\ep_{j_1}}-h_{\ep_{j_2}}$ and  $\pi(h_{\ep_j})=E_{j,j}-E_{-j,-j}\ (1\leq j\leq N)$.

\subsection{Definition}
The type $\mathrm{I}$ and $\mathrm{II}$ vertex operators are the intertwiners of 
the $U_{q,p}(B_N^{(1)})$-modules of the form
\bea
 \widehat{\Phi}(u) &:& \widehat{V}(\lambda,\mu) \to \widehat{V}(\lambda',\mu) \otimes \widehat{V}_z, \\
 \widehat{\Psi}^*(u) &:& \widehat{V}_z \otimes \widehat{V}(\lambda,\mu) \to \widehat{V}(\lambda,\mu'),
\ena
where $\la, \la'\in \h^*,  \mu, \mu' \in H^*$, $z=q^{2u}$. The $\widehat{V}_z$ denote the $(2N+1)$-dimensional  dynamical evaluation module of $U_{q,p}(B_N^{(1)})$ given in Theorem \ref{vecRep} and \ref{vecRepHC}, 
and $\widehat{V}(\lambda,\mu)$ denote the 
level-$k$ highest weight $U_{q,p}(B_N^{(1)})$-module with highest weight $(\lambda,\mu)$. The level-1 case is given in Theorem \ref{levelone}. 
The vertex operators satisfy the intertwining relations with respect to the comultiplication
$\Delta$ given in \eqref{coproUqp}
\bea
 && \Delta(x) \widehat{\Phi}(u) =\widehat{\Phi}(u)x,         \label{irI}\\
 && x \widehat{\Psi}^*(u) = \widehat{\Psi}^*(u) \Delta(x),   \qquad \forall x \in U_{q,p}(B_N^{(1)}).\label{irII}
\ena
These intertwining relations are equivalent to the following relations\cite{Konno08}:
\bea
&& \widehat{\Phi}^{(23)}(u_2)\widehat{L}^{+(12)}(u_1)=
   R^{+(13)}(u_1-u_2,P+h)\widehat{L}^{+(12)}(u_1)\widehat{\Phi}^{(23)}(u_2), 
 \label{typeI}\\
&& \widehat{L}^{+(13)}(u_1)\widehat{\Psi}^{*(23)}(u_2)=
   \widehat{\Psi}^{*(23)}(u_2)\widehat{L}^{+(13)}(u_1)R^{+*(12)}(u_1-u_2,P-h^{(1)}-h^{(2)}).
\label{typeII}
\ena
The relation (\ref{typeI}) (resp.  (\ref{typeII})) should be understood on $\widehat{V}_{z_1} \otimes \widehat{V}(\lambda,\mu)$ (resp. $\widehat{V}_{z_1} \otimes \widehat{V}_{z_2} \otimes \widehat{V}(\lambda,\mu)$). These relations are also expected\cite{JKOS99} from the quasi-Hopf algebra formulation of the face type elliptic quantum group $\Bqla(B_N^{(1)})$\cite{JKOStransfG} by using the connection given in Appendix A. 

We define the components of the vertex operators by
\begin{equation}
 \widehat{\Phi}(u+\frac{1}{2})=\sum_{1 \preceq m \preceq -1}\Phi_m \left(u \right) 
 \otimes v_m, \quad
 \widehat{\Psi}^*(u)(v_m \otimes u)=\Psi_m^* \left(u-\frac{c}{2} \right) \, u,
\end{equation}
where $v_m\in \widehat{V}_z$, $u \in \widehat{V}(\la,\mu)$, and the matrix elements of the $L$-operator $\widehat{L}^+(u)$ by
\begin{equation}
 \widehat{L}^+(u) v_m = \sum_{1 \preceq k \preceq -1} L_{k,m}^+(u)v_k.
\end{equation}
Using these and Corollary \ref{repLR}, the intertwining relations (\ref{typeI}), (\ref{typeII}) are
rewritten as follows:
\bea
&& \hspace{-0.7cm} \Phi_m(u_2) L_{k,j}^+(u_1) =\sum_{1 \preceq m',k' \preceq -1}
 R^+ \big(u_1-u_2+\frac{1}{2},P+h \big)_{km}^{k'm'} L_{k',j}^+(u_1)\Phi_{m'}(u_2),
\label{typeIc} \\
&& \hspace{-0.7cm} L_{k,j}^+(u_1)\Psi_m^*(u_2) =\sum_{1 \preceq j',m' \preceq -1}
\Psi_{m'}^*(u_2) L_{k,j'}^+(u_1)
  R^{+*} \big(u_1-u_2-\frac{c}{2},P-h^{(1)}-h^{(2)} \big)_{j' m'}^{jm}. 
\label{typeIIc}
\ena

\begin{prop}\lb{sufvo}
Let the half currents $E^+_{l,k}(u)$ and $F^+_{k,l}(u)$ $(1\preceq k\prec l\preceq -1)$ take their form as given in Definition \ref{basicHC}  and Appendix B. 
Assume that the top components $\Phi_{-1}(u)$ and $\Psi^*_{-1}(u)$ satisfy the following conditions:
\begin{itemize}
\item[i)] $K^+_{-1}(u_1)\Phi_{-1}(u_2)$ does not have a pole at $u_1-u_2=-\frac{3}{2}$
\item[ii)] $\Psi_{-1}(u_2)K^+_{-1}(u_1)$ does not have a pole at $u_1-u_2=\frac{c-2}{2}+r^*$.
\end{itemize} 
Then the  sufficient conditions for \eqref{typeIc} and \eqref{typeIIc} are given as follows. 
For the type I, 
\bea
&& \Phi_k(u_2)=F_{k,-1}^+(u_2-\frac{1}{2})\Phi_{-1}(u_2)\qquad (1 \preceq k \preceq -2),\lb{typeIk}
\ena
and
\bea
&& \Phi_{-1}(u_2)K_{-1}^+(u_1)=\rho_0^+(u_1-u_2+\frac{1}{2})K_{-1}^+(u_1)\Phi_{-1}(u_2),\\
&& [\Phi_{-1}(u_2),P_l]=0, \quad [\Phi_{-1}(u_2),E_l(u_1)]=0. \qquad (1\leq l\leq N),\\
&& \Phi_{-1}(u_2)(P+h)_{k,-1}=((P+h)_{k,-1}-1)\Phi_{-1}(u_2),
 \label{typIsufcnd3}\\
&& \Phi_{-1}(u_2)F_1(u_1)= \frac{[u_2-u_1-\eta]}
  {[u_2-u_1-\eta-1]}F_1(u_2)\Phi_{-1}(u_1),
\label{typIsufcnd4}\\
&&\Phi_{-1}(u_2)F_j(u_1)=F_j(u_1)\Phi_{-1}(u_2) \quad (2 \le j \le N).
\label{typIsufcnd5}
\ena
For the type II, 
\bea
&&\Psi_k^*(u_2)=\Psi^*_{-1}(u_2)E_{-1,k}^+(u_2+\frac{c}{2}+r^*)
\quad (1 \preceq k \preceq -2),
\ena
and
\bea
&& K_{-1}^+(u_1)\Psi_{-1}^*(u_2)=\Psi_{-1}^*(u_2)K_{-1}^+(u_1)\rho_0^{+*}(u_1-u_2-\frac{c}{2}),\label{typIIsufcnd1}\\
&& [ \Psi_{-1}^*(u_2),(P+h)_l]=0, \quad [\Psi_{-1}^*(u_2),F_l(u_1)]=0\qquad (1\leq l\leq N), 
 \label{typIIsufcnd2}\\
&& P_{j,-1}\Psi_{-1}^*(u_2)=\Psi_{-1}^*(u_2)(P_{j,-1}+1)\qquad  (j\prec -1),
 \label{typIIsufcnd3}\\
&& E_1(u_1)\Psi_{-1}^*(u_2)=\frac{[u_2-u_1-\eta+\frac{1}{2}]^*}{[u_2-u_1-\eta-\frac{1}{2}]^*}
 \Psi_{-1}^*(u_2)E_{1}(u_1),
 \label{typIIsufcnd4}\\
&&E_j(u_1)\Psi_{-1}^*(u_2)=\Psi_{-1}^*(u_2)E_j(u_1) \quad (2 \le j \le N).
\label{typIIsufcnd5}
\ena
\end{prop}
\noindent
{\it Proof.}\ We consider the type I case only. The type II case can be proved similarly. 
From the component $k=m \, (\ne 0)$ in (\ref{typeIc}), we have
\begin{equation}
\Phi_m(u_2)L_{m,j}^+(u_1)
 =\rho_0^+(u_1-u_2+\frac{1}{2})L_{m,j}^+(u_1)\Phi_m(u_2).
\label{typeI1}
\end{equation}
In particular, the component $m=j=-1$ of (\ref{typeI1}) is 
\begin{equation}
 \Phi_{-1}(u_2)K_{-1}^+(u_1)=\rho_0^+(u_1-u_2+\frac{1}{2})K_{-1}^+(u_1)\Phi_{-1}(u_2).
 \label{typIsufcnd1}
\end{equation}
Note the formula 
\begin{eqnarray}
&&{\rho}_0^+(u)=\frac{[u+1]}{\varphi(u)},\\
&&\varphi(u)=q^{-1} z^{\frac{1}{r}}[u-1]\frac{\{\xi^2z\}\{z\}\{\xi q^2z\}\{\xi q^{-2}z\}}{\{\xi z\}^2\{\xi^2 q^2z\}\{q^{-2}z\}}\frac{\{p\xi /z\}^2\{p\xi^2 q^2/z\}\{pq^{-2}/z\}}{\{p\xi^2/z\}\{p/z\}\{p\xi q^2/z\}\{p\xi q^{-2}/z\}}. 
\end{eqnarray}
Then from the assumption i), \eqref{typIsufcnd1} implies that $\Phi_{-1}(u_2)K_{-1}^+(u_1)$ has a zero at $u_1-u_2=-\frac{3}{2}$. 
We will check these points for the level-1 representation. 

In addition, from the component $m=-1 \succ j$ of (\ref{typeI1}), and putting the definition $L_{-1,j}^+(u)=K_{-1}^+(u)E_{-1,j}^+(u)$, we have
\begin{equation}
 \Phi_{-1}(u_2)E_{-1,j}^+(u_1)=E_{-1,j}^+(u_1)\Phi_{-1}(u_2).\lb{PhiE}
\end{equation}
From the conjectural expressions for $E_{-1,j}^+(u)$  in Appendix B, 
the sufficient conditions for \eqref{PhiE} are
\begin{equation}
 [\Phi_{-1}(u_2),P_l]=0, \quad [\Phi_{-1}(u_2),E_l(u_1)]=0. \qquad (1\leq l\leq N)
 \label{typIsufcnd2}
\end{equation}

Next, the component $k \ne \pm m \, (\ne 0), \, k \prec m$ in (\ref{typeIc}) is 
\begin{eqnarray}
 \Phi_m(u_2) L_{k,j}^+(u_1)&=&\rho_0^+(u_1-u_2+\frac{1}{2}) 
  \Big\{b(u_1-u_2+\frac{1}{2},(P+h)_{k,m}) L_{k,j}^+(u_1)\Phi_m(u_2) \nonumber \\
   && \quad +c(u_1-u_2+\frac{1}{2},(P+h)_{k,m}) L_{m,j}^+(u_1)\Phi_k(u_2) \Big\}.
\label{typeI2}
\end{eqnarray}
Then putting the definition $L_{-1,-1}^+(u)=K_{-1}^+(u)$ and 
$L_{k,-1}^+(u)=F_{k,-1}^+(u)K_{-1}^+(u)$ in the case $k \prec m=j=-1$ of (\ref{typeI2}),  we have
\begin{eqnarray}
 && \Phi_{-1}(u_2)F_{k,-1}^+(u_1)K_{-1}^+(u_1) \nonumber \\
&& \quad =\rho_0^+(u_1-u_2+\frac{1}{2}) \Big\{ \frac{[(P+h)_{k,-1}+1][(P+h)_{k,-1}-1][u_1-u_2+\frac{1}{2}]}
{[(P+h)_{k,-1}]^2[u_1-u_2+\frac{3}{2}]}F_{k,-1}^+(u_1)K_{-1}^+(u_1)\Phi_{-1}(u_2) \nonumber \\
&& \quad+\frac{[1][(P+h)_{k,-1}+u_1-u_2+\frac{1}{2}]}{[P+h)_{k,-1}][u_1-u_2+\frac{3}{2}]}K_{-1}^+(u_1)\Phi_k(u_2) \Big\}.
\label{intertwinI}
\end{eqnarray}
Putting $u_1-u_2=-\frac{3}{2}$, the left hand side of (\ref{intertwinI}) vanishes.
Then we obtain for $2 \preceq k \preceq -2$
\begin{eqnarray}
 \Phi_k(u_2) &=& K_{-1}^{+}(u_1)^{-1}\frac{[(P+h)_{k,-1}+1]}{[(P+h)_{k,-1}]}F_{k,-1}^{+}(u_1)K_{-1}^{+}(u_1)\Phi_{-1}(u_2) \nn \\
  &=& F_{k,-1}^{+}(u_1+1)\Phi_{-1}(u_2) \nn \\
  &=& F_{k,-1}^{+}(u_2-\frac{1}{2})\Phi_{-1}(u_2).
\label{typeIk2}
\end{eqnarray}

The remaining component $\Phi_1(u)$ is also obtained from $\Phi_{-1}(u)$ as follows. 
From the component $m=j=-1$ in (\ref{typeIc}), we have
\begin{eqnarray*}
\Phi_{-1}(u_2)F_{1,-1}^{+}(u_1)K_{-1}^{+}(u_1)
 &=& R^+(u_1-u_2+\frac{1}{2}, P+h)_{1-1}^{1-1}F_{1,-1}^{+}(u_1)K_{-1}^{+}(u_1)\Phi_{-1}(u_2) \nonumber \\
 &&+R^{+}(u_1-u_2+\frac{1}{2}, P+h)_{1-1}^{{-1}1}K_{-1}^{+}(u_1)\Phi_{1}(u_2) \nn \\
 &&+\sum_{2\preceq l\preceq -2} R^{+}(u_1-u_2+\frac{1}{2},P+h)_{1-1}^{{-l}l}F_{-l,-1}^+(u_1)K_{-1}^+(u_1)F_{l,-1}(u_2)\Phi_{-1}(u_2).
\end{eqnarray*}
Using \eqref{typeIk2} and the component $(i,j)=(1,-1)$, $(i'',j'')=(-1,-1)$ of
\eqref{RLLcomp}, we obtain
\be
&& \Phi_{-1}(u_2)F_{1,-1}^+(u_1)K_{-1}^+(u_1)\\
 &&=R^+(u_1-u_2+\frac{1}{2},P+h)_{1-1}^{{-1}1} K_{-1}^+(u_1)\Phi_{1}(u_2)\\
&&+\rho_0^+(u_1-u_2+\frac{1}{2})K_{-1}^+(u_2)F_{1,-1}^+(u_1)K_{-1}^+(u_2)^{-1}K_{-1}^+(u_1)\Phi_{-1}(u_2)\\
&& -R^+(u_1-u_2+\frac{1}{2},P+h)_{1-1}^{{-1}1}K_{-1}^+(u_1)F_{1,-1}^+(u_2)\Phi_{-1}(u_2).
\en
Then again setting $u_1-u_2= -\frac{3}{2}$,  the left hand side and 
the second term in right hand side vanish. Then we obtain
\bea 
&&\Phi_1(u_2)=F_{1,-1}^+(u_2-\frac{1}{2})\Phi_{-1}(u_2).\lb{Phi1}
\ena
Combining \eqref{typeIk2} and \eqref{Phi1}, we obtain  (\ref{typeIk}). 
 
Furthermore, substituting (\ref{typeIk}) into (\ref{intertwinI}), we obtain the sufficient conditions \eqref{typIsufcnd3}-\eqref{typIsufcnd5}. 

\begin{lem}\lb{PhiF}
For $ 1 \preceq k \preceq -2$,  we have
\begin{eqnarray}
 \Phi_{-1}(u_2)F_{k,-1}(u_1-\frac{1}{2}) &=& K_{-1}^+(u_2-\frac{1}{2})F_{k,-1}^+(u_1-\frac{1}{2})
 K_{-1}^+(u_2-\frac{1}{2})^{-1}\Phi_{-1}(u_2), \nn \\
 E_{-1,k}^+(u_1+\frac{c}{2})\Psi_{-1}^*(u_2) &=& \Psi_{-1}^*(u_2)K_{-1}^+(u_2+\frac{c}{2})^{-1}
 E_{-1,k}^+(u_1+\frac{c}{2})K_{-1}^+(u_2+\frac{c}{2}). \nn
\end{eqnarray}
\end{lem}
\noindent
{\it Proof.} \quad
From the component $m=-1, \, j=-1$ in the intertwining relation
(\ref{typeIc}), we have
\begin{eqnarray}
 &&\Phi_{-1}(u_2) F_{k,-1}^+(u_1)K_{-1}^+(u_1) \nn \\
 && =\rho_0^+(u_1-u_2+\frac{1}{2})b(u_1-u_2+\frac{1}{2},(P+h)_{k,-1})
     F_{k,-1}^+(u_1)K_{-1}^+(u_1)\Phi_{-1}(u_2) \nn \\
 && +\rho_0^+(u_1-u_2+\frac{1}{2})c(u_1-u_2+\frac{1}{2},(P+h)_{k,-1})
     K_{-1}^+(u_1)\Phi_{k}(u_2).
\end{eqnarray}
Using (\ref{typeIk}), we have
\begin{eqnarray}
 &&\Phi_{-1}(u_2) F_{k,-1}^+(u_1)K_{-1}^+(u_1) \nn \\
 && =\rho_0^+(u_1-u_2+\frac{1}{2})b(u_1-u_2+\frac{1}{2},(P+h)_{k,-1})
     F_{k,-1}^+(u_1)K_{-1}^+(u_1)\Phi_{-1}(u_2) \nn \\
 && +\rho_0^+(u_1-u_2+\frac{1}{2})c(u_1-u_2+\frac{1}{2},(P+h)_{k,-1})
     K_{-1}^+(u_1)F_{k,-1}^+(u_2-\frac{1}{2})\Phi_{-1}(u_2).
\end{eqnarray}
From (\ref{typIsufcnd1}), we have
\begin{eqnarray}
 && \Phi_{-1}(u_2) F_{k,-1}^+(u_1-\frac{1}{2}) \nn \\
 && =b(u_1-u_2,(P+h)_{k,-1})
     F_{k,-1}^+(u_1-\frac{1}{2})\Phi_{-1}(u_2) \nn \\
 && +c(u_1-u_2,(P+h)_{k,-1})
     K_{-1}^+(u_1-\frac{1}{2})F_{k,-1}^+(u_2-\frac{1}{2})K_{-1}(u_1-\frac{1}{2})^{-1}\Phi_{-1}(u_2).
\end{eqnarray}
Then it is sufficient to show
\begin{eqnarray}
&& b(u_1-u_2,(P+h)_{k,-1})F_{k,-1}^+(u_1)+
   c(u_1-u_2,(P+h)_{k,-1})K_{-1}^+(u_1)F_{k,-1}^+(u_2)K_{-1}^+(u_1)^{-1} \nn \\
&& =K_{-1}^+(u_2)F_{k,-1}^+(u_1)K_{-1}^+(u_2)^{-1}.
\end{eqnarray}
This is nothing but the component $(i,j)=(-1,-1), \, (i'',j'')=(k,-1)$ of 
\eqref{RLLcomp}.
\qed

\subsection{Level one vertex operators and commutation relations}
Next we consider a free field realization of the vertex operators fixing the 
representation level $c=1$.

From the sufficient conditions obtained in Propostion \ref{sufvo}, 
we can determine the free field realizations of 
vertex operators as follows:

\begin{prop}\lb{bareboson}
The highest components of the type I and type II vertex operators 
$\Phi_{-1}(u)$ and $\Psi_{-1}^*(u)$ are realized in terms of the free
field by
\begin{align}
 \Phi_{-1}(u) &= : \exp \left\{ \sum_{m \ne 0} (q^m-q^{-m})
  \frac{1-p^{*m}}{1-p^m} \cE_m^{-1}(q^{-3}\xi z)^{-m} \right\}:
  e^{\epsilon_1} (q^{-1} \xi z)^{h_{\epsilon_1}}(q^{-1} \xi z)^{-\frac{1}{r}(P+h)_{\epsilon_1}}
 , \label{typeIvo}\\
 \Psi_{-1}^*(u) &=: \exp \left\{ -\sum_{m \ne 0} (q^m-q^{-m})
   \cE_m^{-1}(q^{-1}\xi z)^{-m} \right\}: 
  e^{-\epsilon_1}e^{Q_{\epsilon_1}}(\xi z)^{-h_{\epsilon_1}} (\xi z)^{\frac{1}{r^*}P_{\epsilon_1}}
  . \label{typeIIvo}
\end{align}
These realizations satisfy the assumptions i) and ii) in Proposition \ref{sufvo}.
\end{prop}
\noindent
{\it Proof.} \quad
By straightforward calculations, we can show that (\ref{typeIvo}) (resp. (\ref{typeIIvo}))
satisfies the sufficient conditions (\ref{typIsufcnd1}), (\ref{typIsufcnd2}), and 
(\ref{typIsufcnd3})-(\ref{typIsufcnd5}) (resp. (\ref{typIIsufcnd1})-(\ref{typIIsufcnd5})).
\qed

\begin{thm}
The free field realizations of the type I $\Phi_j(u)$ and the type II  $\Psi_j^*(u)$ vertex operators satisfy the following commutation relations:
\begin{eqnarray}
 \Phi_{j_2}(u_2)\Phi_{j_1}(u_1) &=& \sum_{j_1',j_2'=1}^{-1}R(u_1-u_2,P+h)_{j_1 j_2}^{j_1' j_2'}
 \Phi_{j_1'}(u_1)\Phi_{j_2'}(u_2), \label{typeIcr}\\
 \Psi_{j_1}^*(u_1)\Psi_{j_2}^*(u_2) &=& \sum_{j_1',j_2'=1}^{-1}\Psi_{j_2'}^*(u_2)\Psi_{j_1'}^*(u_1)
 R^{*}(u_1-u_2,P-h^{(1)}-h^{(2)})_{j_1' j_2'}^{j_1 j_2}, \label{typeIIcr} \\
 \Phi_j(u_1)\Psi_k^*(u_2) &=& \chi(u_1-u_2)\Psi_k^*(u_2)\Phi_j(u_1). \label{typeI-II} 
\end{eqnarray}
Here we set
\begin{equation}
 R(u,P)=\mu(u)\bar{R}^+(u,P), \quad R^*(u,P)=\mu^*(u)\bar{R}^{+*}(u,P)
\end{equation}
with
\bea
&& \mu(u)=z^{-1+\frac{1}{r}} \frac{ \{p \xi^2 q^{-2} z\} \{p \xi z\}
 \{\xi z\} \{q^2 z \}}{\{p\xi q^{-2} z \} \{p z \}
 \{\xi^2 z \} \{\xi q^2 z \}} 
  \frac{ \{p\xi q^{-2} /z \} \{p /z \} 
  \{\xi^2 /z \} \{\xi q^2 /z \} }{\{p \xi^2 q^{-2}/z \}
  \{p \xi /z \} \{\xi /z \} \{q^2 /z \} },\\
&& \mu^*(u)=\mu(u) |_{r \to r^*}
\ena
and
\begin{equation}
 \chi(u)=\frac{\Theta_{\xi^2}(z)\Theta_{\xi^2}(q^{-2} \xi z)}
              {\Theta_{\xi^2}(\xi z)\Theta_{\xi^2}(q^{-2} \xi^2 z)}.
\end{equation}
\end{thm}
\noindent
{\it Proof.} \quad
Let us show the commutation relation of the type I vertex operators (\ref{typeIcr}). 
For $j_1=j_2=-1$, the equation
\begin{equation}
\Phi_{-1}(u_2)\Phi_{-1}(u_1) = \mu(u_1-u_2) \Phi_{-1}(u_1)\Phi_{-1}(u_2)
\label{typeIbarecr}
\end{equation}
can be shown by straightforward calculation with the use of the free field 
realization (\ref{typeIvo}).

For $1 \preceq j_1, \, j_2 \preceq -2$, 
using (\ref{typeIk}), (\ref{typeIbarecr}) and
Lemma \ref{PhiF},  the equation
(\ref{typeIcr}) is reduced to the following equation:
\begin{eqnarray}
&& F_{j_2,-1}^+(u_2)K_{-1}^+(u_2)F_{j_1,-1}^+(u_1)
K_{-1}^+(u_2)^{-1} \label{typeIred} \nn \\
&&=\sum_{j_1',j_2'=1}^{-1} \bar{R}^+(u,P+h)_{j_1 j_2}^{j_1' j_2'}
F_{j_1',-1}^+(u_1)K_{-1}^+(u_1)F_{j_2',-1}^+(u_2)
K_{-1}^+(u_1)^{-1}, 
\end{eqnarray}
where $u=u_1-u_2$.
From the component $(i,j)=(-1,-1),\,(i'',j'')=(j_1,j_2)$ of \eqref{RLLcomp}, we have
\begin{eqnarray}
&&\sum_{j_1',j_2'=1}^{-1}R^+(u,P+h)_{j_1 j_2}^{j_1'j_2'}
F_{j_1,-1}^+(u_1)K_{-1}^+(u_1)F_{j_2,-1}^+(u_2)K_{-1}^+(u_2) \nn  \\
&&=F_{j_2,-1}^+(u_2)K_{-1}^+(u_2)F_{j_1,-1}^+(u_1)K_{-1}^+(u_1)\tilde{\rho}^{*+}(u).
\end{eqnarray}
Multiplying the above equation by
\begin{equation}
 \tilde{\rho}^+(u)^{-1}K_{-1}^+(u_2)^{-1}K_{-1}^+(u_1)^{-1}=
 K_{-1}^+(u_1)^{-1}K_{-1}^+(u_2)^{-1}\tilde{\rho}^{*+}(u)^{-1}
 \label{commkk}
\end{equation}
from the right, we obtain the desired equation (\ref{typeIred}).

For $j_1=-1, \, 1 \preceq j_2 \preceq -2$, using (\ref{typeIk}), (\ref{typeIbarecr})
and Lemma \ref{PhiF}, the equation (\ref{typeIcr}) is reduced to the following equation
of the half currents:
\begin{eqnarray}
 &&F_{j_2,-1}^+(u_2) \label{typeI-1} \\
 && = \bar{R}^+(u,P+h)_{-1 j_2}^{-1 -1}
 +\sum_{j_2'=1}^{-2}\bar{R}^+(u,P+h)_{-1 j_2}^{-1 j_2'}K_{-1}^+(u_1-\frac{1}{2})
F_{j_2',-1}^+(u_2-\frac{1}{2})K_{-1}^+(u_1-\frac{1}{2})^{-1} \nn \\
 &&+\sum_{j_1'=1}^{-2}\bar{R}^+(u,P+h)_{-1 j_2}^{j_1' -1}F_{j_1',-1}^+(u_1-\frac{1}{2}) \nn \\
 &&+\sum_{j_1',j_2'=1}^{-2}\bar{R}^+(u,P+h)_{-1 j_2}^{j_1' j_2'}F_{j_1',-1}^+(u_1-\frac{1}{2})
K_{-1}^+(u_1-\frac{1}{2})F_{j_2',-1}^+(u_2-\frac{1}{2})K_{-1}^+(u_1-\frac{1}{2})^{-1}. \nn 
\end{eqnarray}
From the component $(i,j)=(-1,-1),\,(i'',j'')=(-1,j_2)$ of \eqref{RLLcomp}, we have
\begin{eqnarray}
 && R^+(u,P+h)_{-1 j_2}^{-1 -1}K_{-1}^+(u_1)K_{-1}^+(u_2) \nn \\
 && +\sum_{j'=1}^{-2}R^+(u,P+h)_{-1 j_2}^{-1 j'}K_{-1}^+(u_1)F_{j',-1}^+(u_2)K_{-1}^+(u_2) \nn \\
 && +\sum_{i'=1}^{-2}R^+(u,P+h)_{-1 j_2}^{i' -1}F_{i',-1}^+(u_1)K_{-1}^+(u_1)K_{-1}^+(u_2) \nn \\
 && +\sum_{i',j'=1}^{-2}R^+(u,P+h)_{-1 j_2}^{i' j'}F_{i',-1}^+(u_1)K_{-1}^+(u_1)
F_{j',-1}^+(u_2)K_{-1}^+(u_2) \nn \\
&&=F_{j_2,-1}^+(u_2)K_{-1}^+(u_2)K_{-1}^+(u_1)\rho^{+*}(u).
\end{eqnarray}
Multiplying the above equation by (\ref{commkk}) from the right, we obtain the desired equation
(\ref{typeI-1}).
The case $1 \preceq j_1 \preceq -2, \, j_2=-1$ can be proved in the same manner.

Similarly, one  can prove the commutation relation of the type II vertex operators (\ref{typeIIcr}).

Next, let us consider the relation (\ref{typeI-II}). 
The case $j=k=-1$ is a direct consequence from Proposition \ref{bareboson}. 
The cases $j=-1$ or $k=-1$ can be shown as follows:
consider the case $k=-1$ for instance. By (\ref{typeIk}),
\begin{eqnarray}
\Phi_j(u_1) \Psi_{-1}^*(u_2) &=& F_{j,-1}^+(u_1-\frac{1}{2})\Phi_{-1}(u_1)\Psi_{-1}^*(u_2) \nn \\
&=& F_{j,-1}^+(u_1-\frac{1}{2}) \chi(u_1-u_2)\Psi_{-1}^*(u_2)\Phi_{-1}(u_1) \nn \\
&=& \chi(u_1-u_2)\Psi_{-1}^*(u_2)F_{j,-1}^+(u_1-\frac{1}{2})\Phi_{-1}(u_1) \quad 
\nn \\
&=& \chi(u_1-u_2)\Psi_{-1}^*(u_2)\Phi_j(u_1).
\label{IjII-1}
\end{eqnarray}
Then the general case is proved as follows. Since both $\widehat{\Phi}(u_1)\widehat{\Psi}^*(u_2)$ and
$(\widehat{\Psi}^*(u_2) \otimes \id) (\id \otimes \widehat{\Phi}(u_1))$ commute with
$\Delta(x) \, (\forall x \in U_{q,p}(B_N^{(1)}))$, they act as scalars on the irreducible module
$V_{z_2} \otimes \widehat{V}(\lambda)$. In order to compare the scalars, we will see 
their actions on $v_{-1} \otimes |\lambda \rangle \in V_{z_{2}} \otimes \widehat{V}(\lambda)$.
\begin{eqnarray}
\widehat{\Phi}(u_1) \widehat{\Psi}^*(u_2)(v_{-1} \otimes | \lambda \rangle)
&=&\widehat{\Phi}(u_1) \Psi_{-1}^*(u_2-\frac{c}{2})|\lambda \rangle \nn \\
&=&\sum_{j} \Phi_j(u_1)\Psi_{-1}^*(u_2-\frac{c}{2})|\lambda \rangle \otimes v_j \nn \\
&=&\chi(u_1-u_2+\frac{1}{2}) \sum_{j} \Psi_{-1}^*(u_2-\frac{c}{2})\Phi_j(u_1+\frac{1}{2})
|\lambda \rangle \otimes v_j.
\label{I-IIl}
\end{eqnarray}
Here the last equality follows from (\ref{IjII-1}). On the other hand 
\begin{eqnarray}
(\widehat{\Psi}^*(u_2) \otimes \id)(\id \otimes \widehat{\Phi}(u_1))(v_{-1} \otimes | \lambda \rangle )
&=&(\widehat{\Psi}^*(u_2) \otimes \id)(v_{-1} \otimes \sum_{j}\Phi_j(u_1)|\lambda \rangle \otimes v_j )\nn \\
&=&\Psi_{-1}^*(u_2-\frac{c}{2})\sum_{j}\Phi_j(u_1+\frac{1}{2})|\lambda \rangle \otimes v_j.
\label{I-IIr}
\end{eqnarray}
Comparing (\ref{I-IIl}) and (\ref{I-IIr}), we get
\begin{equation}
\widehat{\Phi}(u_1) \widehat{\Psi}^*(u_2)=
\chi(u_1-u_2+\frac{1+c}{2})(\widehat{\Psi}^*(u_2) \otimes \id)(\id \otimes \widehat{\Phi}(u_1)).
\end{equation}
Hence comparing the components of the both sides, 
and changing variables $u_1+\frac{1}{2} \to u_1, \, u_2-\frac{c}{2} \to u_2$,
we obtain
\begin{equation}
\Phi_j(u_1) \Psi_k^*(u_2)=\chi(u_1-u_2) \Psi_k^*(u_2) \Phi_j(u_1).
\end{equation}
\qed

\section*{Acknowledgements}
The authors would like to thank Michio Jimbo, Masatoshi Noumi and Masato Okado for useful conversations. H.K is also grateful to the organizers of the workshop {\it Elliptic Integrable Systems and Hypergeometric Functions} for a kind invitation. K.O thanks the organizers of the conference {\it RIMS symposium : String
Theory, Integrable Systems and Representation Theory}, especially Koji Hasegawa, 
for an invitation. 
H.K was supported by the Grant-in -Aid for Scientific Research (C) 22540022 JSPS, Japan.

\begin{appendix}

\section{Relation to the Quasi-Hopf Formulation $\Bqla(\ghbig)$}

\subsection{Definition of $\Bqla(\Bnh)$}
Let $U_q=U_{q}(\Bnh)$ be the Drinfeld-Jimbo affine quantum group\cite{Drinfeld85,Jimbo}. Namely, $U_q(\Bnh)$ is a quasi-triangular
Hopf algebra realized by the Chevalley generators and equipped with the standard coproduct $\Delta_0$,  counit $\vep$, 
antipode $S$ and universal $R$ matrix $\cR$. 
Our conventions on the coalgebra structure follows \cite{JKOStransfG}. 
Let $\h$ and $\bar{\h}$ be the Cartan subalgebras as in Sec.2.1. 
We denote a basis and its dual basis of ${\h}$ by $\{\hat{h}_l\}$ and $\{\hat{h}^l\}$, respectively. More explicitly, they are given by $\{\hat{h}_l\}=\{d, c, h_j \}$ and $\{\hat{h}^l\}=\{c, d, h^j \}\ (1\leq j\leq N)$, where 
$\{{h}_j\}$ and $\{{h}^j\}$ are a basis  and a dual basis of 
$\bar{\h}$. 

The face type elliptic quantum group $\Bqla(\Bnh)$ is a quasi-Hopf
deformation of $U_{q}(\Bnh)$ by the face type twistor $F(\la)\ (\la \in \h)$.
The twistor $F(\la)$ is an invertible element in $U_q\otimes U_q$ 
satisfying 
\bea
&&({\rm id}\otimes \vep )F(\la)=1=F(\la)(\vep \otimes {\rm id}),\\
&&F^{(12)}(\la)(\Delta_0\otimes {\rm id})F(\la)
=F^{(23)}(\la+h^{(1)})({\rm id} \otimes \Delta_0)F(\la).
\lb{facecocy}
\lb{epF}
\ena
where $\la=\sum_l\la_l\hat{h}^l\  (\la_l\in \C)$, 
$\la+h^{(1)}=\sum_l(\la_l+\hat{h}_l^{(1)})\hat{h}^l$ and 
$\hat{h}_l^{(1)}=\hat{h}_l\otimes 1\otimes 1$.  
An explicit construction of the twistor $F(\la)$ is given in \cite{JKOStransfG}.
Then we define a new coproduct by
\bea
&&\Delta_\la(x)=F(\la)\Delta_0(x)F(\la)^{-1}\qquad \forall x \in U_q(\Bnh).
\ena
$\Delta_\la$ satisfies a weaker coassociativity
\bea
&&({\rm id}\otimes \Delta_\la )\Delta_\la(x)=\Phi(\la)
(\Delta_\la \otimes {\rm id})\Delta_\la(x)\Phi(\la)^{-1}\qquad \forall\ x \in U_q(\Bnh),\\
&&\Phi(\la)=F^{(23)}(\la)F^{(23)}(\la+h^{(1)})^{-1}.
\ena
The universal $R$-matrix is also deformed to
\bea
&&\cR(\la)=F^{(21)}(\la)\cR F^{(12)}(\la)^{-1}.
\ena

\begin{dfn}
\cite{JKOStransfG}~~
The face type elliptic quantum group $\Bqla(\Bnh)$ is a quasi-triangular
quasi-Hopf algebra 
$(U_{q}(\Bnh),$
$~\Delta_\la, \vep,~S,~\Phi(\la),~\alpha,~\beta,~\cR(\la))$,
where $\alpha,\ \beta$ are defined by 
\bea
&&\alpha=\sum_iS(k_i)l_i,\quad \beta=\sum_i m_iS(n_i).
\ena
Here we set $\sum_ik_i\otimes l_i=F(\la)^{-1},\ 
\sum_im_i\otimes n_i=F(\la)$.
\end{dfn}

The new universal $R$ matrix $\cR(\la)$ satisfies the dynamical Yang-Baxter equation.
\begin{equation}
\cR^{(12)}(\la+ h^{(3)} )\cR^{(13)}(\la)\cR^{(23)}(\la+h^{(1)})
=\cR^{(23)}(\la)\cR^{(13)}(\la+h^{(2)})\cR^{(12)}(\la).
\label{DYBE}
\end{equation}

In \cite{Konno06}, we derived vector representations of $\cR(\la)$ for $\gh=A_N^{(1)}, B_N^{(1)}, C_N^{(1)}, D_N^{(1)}$ and found  that if 
we parametrize $\la\in \hh^*$ as $\la=\la(r^*,P)=(r^*+h^\vee)d+\sum_{j=1}^{N}(P_{\alpha'_j}+1)\bar{\Lambda}_j$ with $\alpha'_j$ being the simple roots of the dual Lie algebra $\bar{\g}^\vee$ of $\bar{\g}$, the vector representation of  $R(\la)$ coincides with the corresponding face weight derived by Jimbo, Okado and Miwa\cite{JMO}. In particular, for the $B_N^{(1)}$ type, if we set    
\be
\tR^{+*}(z,P)&=&(\pi_{V,z}\otimes \pi_{V,1})\left({\rm Ad}\ z^{-\frac{\bar{\theta}(\la)}{r}}\otimes \id\right)\left(z^{\frac{\bar{T}}{r}}q^{c\otimes d+d\otimes c}\cR(\la)\right),
\\
\tR^+(z,P+h)&=&(\pi_{V,z}\otimes \pi_{V,1})\left({\rm Ad}\ z^{-\frac{\bar{\theta}(\la)}{r}}\otimes \id\right)\left(z^{\frac{\bar{T}}{r}}q^{c\otimes d+d\otimes c}\cR(\la+h)\right),
\en
with
\be
&&\bar{\theta}(\la)=-\bar{\la}+\bar{\rho}-\frac{1}{2}\sum_j\bar{h}_j\bar{h}^j,\\
&&\bar{T}=\sum_j\bar{h}_j\otimes \bar{h}^j
\en
for $\la=\la(r^*,P)$, then $\tR^{+*}(z,P)$ and $\tR^+(z,P+h)$ coincide with \eqref{ellR} up to a gauge transformation.  
Moreover we define the $L$ operators of $\Bqla(\gh)$ by  
\be
\tL^{+ }(z,P)&=&(\pi_{V,z}\otimes \id)\left({\rm Ad}\ z^{-\frac{\bar{\theta}(\la)}{r}}\otimes \id\right)\left(z^{\frac{\bar{T}}{r}}q^{c\otimes d+d\otimes c}\cR(\la)\right),\\
\tL^{-}(z,P)&=&(\pi_{V,z}\otimes \id)\left({\rm Ad}\ z^{-\frac{\bar{\theta}(\la)}{r}}\otimes \id\right)\left(z^{\frac{\bar{T}}{r}}\cR^{(21)}(\la)^{-1}q^{-c\otimes d-d\otimes c}\right). 
\en
Note that $\tL^{+ }(z,P)$ and $\tL^{-}(z,P)$ are not independent: 
we have 
\bea
&&\tL^-(z,P)=\tL^+(z p^* q^{c},P). \lb{LmfromLpBqla}
\ena
In addition, if we define 
\be
\tR^{-*}(z,P)&=&(\pi_{V,z_1}\otimes \pi_{V,z_2})\cR^{(21)}(\la)^{-1}q^{-c\otimes d-d\otimes c},\\
\tR^{-}(z,P+h)&=&(\pi_{V,z_1}\otimes \pi_{V,z_2})\cR^{(21)}(\la+h)^{-1}q^{-c\otimes d-d\otimes c}, 
\en
Then we have $\tR^{-*}(z,P)=\tR^{+*}(zp^*q^c,P), \tR^{-}(z,P+h)=\tR^{+}(zpq^{-c},P+h)$.
Combining these formulas we obtain from \eqref{DYBE} the following dynamical $RLL$ relations \cite{JKOS99}
\bea
&&\tR^\pm(z,P+h)\tL^\pm(z_1,P)\tL^-(z_2,P+h^{(1)})=\tL^\pm(z_2,P)\tL^\pm(z_1,P+h^{(2)})
\tR^{\pm*}(z,P),\lb{dynRLL}\\
&&\tR^\pm(zq^{\pm{c}},P+h)\tL^\pm(z_1,P)\tL^\mp(z_2,P+h^{(1)})=\tL^\mp(z_2,P)\tL^\pm(z_1,P+h^{(2)})\tR^{\pm*}(zq^{\mp{c}},P).\nn\\
&&\lb{dynRLLpm}
\ena
Furthermore define  
\be
&&\widehat{\cL}^\pm(z)=\tL^\pm(z,P)e^{-\sum_j\pi_V(h_{\bep_j})\otimes Q_{\bep_j}} \in \B_{q,\la(r^*,P)}(\gh)
[[z,z^{-1}]][z^{\pm \frac{1}{r^*}}, z^{\pm \frac{1}{r}}]\sharp \C[\cR_Q], 
\en
where 
 $\pi_V(h_{\vep_j})=E_{jj}-E_{-j-j}$ for the case $B_N^{(1)}
 $. 
Then one can verify that $\widehat{\cL}^\pm(z)$ satisfy the $RLL$ relations 
\bea
&&\tR^\pm(z,P+h)\hcL^\pm(z_1)\hcL^-(z_2)=\hcL^\pm(z_2)\hcL^\pm(z_1)
\tR^{\pm*}(z,P),\lb{RLLB}\\
&&\tR^\pm(zq^{\pm{c}},P+h)\hcL^\pm(z_1)\hL^\mp(z_2)=\hcL^\mp(z_2)\hL^\pm(z_1)\tR^{\pm*}(zq^{\mp{c}},P).\lb{RLLpmB}
\ena
These $RLL$-relations coincide with \eqref{RLLpm} and \eqref{RLLmp}. 

\section{Integral Expressions for the Half Currents}
For $X=E, F$, let us denote by $[X_{j_1}(v_1)\cdots X_{j_m}(v_m)]$ the product of the elliptic currents 
$X_{j_1}(v_1),\cdots ,X_{j_m}(v_m)$ where all the zero-modes, $w_j^{-\frac{P_{\al_j}-1}{r^*}}$ of $E_j(v_j)\ (j=1,\cdots ,N-1)$,  $w_N^{-\frac{P_{\al_N}-1/2}{r^*}}$ of $E_N(v_N)$, $w_j^{\frac{(P+h)_{\al_j}-1}{r}}$ of $F_j(v_j)\ (j=1,\cdots ,N-1)$ and 
 $w_N^{\frac{(P+h)_{\al_N}-1/2}{r}}$ of $F_N(v_N)$, are normally ordered, i.e. they are 
moved to the right of all $e_j(w_j)$ and $f_j(w_j)$ by using \eqref{ge} and \eqref{gf}.  

\noindent
i) The $j\prec k \prec N$ case: 
\bea
&&[E_j(v_j)E_{j+1}(v_{j+1})\cdots E_k(v_k)\cdots E_N(v_N)E_N(v_N')\cdots E_k(v_k')]\nn\\  
&&=e_j(v_j)e_{j+1}(v_{j+1})\cdots e_k(v_k)\cdots e_N(v_N)e_N(v_N')\cdots e_k(v_k')
\prod_{m=j}^{k-2}w_m^{-\frac{P_{\al_m}}{r^*}}
\cdot w_{k-1}^{-\frac{P_{\al_{k-1}}+1}{r^*}}\cdot w_k^{-\frac{P_{\al_k}-1}{r^*}} \nn\\
&&\qquad\times  \prod_{m=k+1}^{N-1}w_m^{-\frac{P_{\al_m}}{r^*}}\cdot 
w_N^{-\frac{P_{\al_N}-1/2}{r^*}} {w'}_N^{-\frac{P_{\al_N}+1/2}{r^*}}
 \prod_{m=k+1}^{N-1}{w'}_m^{-\frac{P_{\al_m}}{r^*}} \cdot {w'}_k^{-\frac{P_{\al_k}-1}{r^*}},\lb{jkNa}\\
&&[E_k(v_k')E_{k+1}(v_{k+1}')\cdots E_N(v_N')E_N(v_N)\cdots  E_k(v_k)\cdots E_k(v_j)]\nn\\  
&&= e_k(v_k')e_{k+1}(v_{k+1}')\cdots e_N(v_N')e_N(v_N)\cdots  e_k(v_k)\cdots e_k(v_j)\prod_{m=k}^{N-1}
{w'}_m^{-\frac{P_{\al_m}}{r^*}}
\cdot 
{w'}_N^{-\frac{P_{\al_N}-1/2}{r^*}} {w}_N^{-\frac{P_{\al_N}+1/2}{r^*}}\nn\\
&&\qquad \times 
 \prod_{m=j+1}^{N-1}{w}_m^{-\frac{P_{\al_m}}{r^*}} \cdot {w}_j^{-\frac{P_{\al_j}-1}{r^*}}.\lb{jkNb}
\ena

\noindent
ii) The $j = k \prec N$ case: 
\bea
&&[E_j(v_j)E_{j+1}(v_{j+1})\cdots E_N(v_N)E_N(v_N')\cdots E_j(v_j')]\nn\\
&&=e_j(v_j)e_{j+1}(v_{j+1})\cdots e_N(v_N)e_N(v_N')\cdots e_j(v_j')w_j^{-\frac{P_{\al_j}-1}{r^*}}\prod_{m=j+1}^{N-1}w_m^{-\frac{P_{\al_m}}{r^*}}\cdot 
w_N^{-\frac{P_{\al_N}-1/2}{r^*}} {w'}_N^{-\frac{P_{\al_N}+1/2}{r^*}}\nn\\
&&\qquad \times 
 \prod_{m=j+1}^{N-1}{w'}_m^{-\frac{P_{\al_m}}{r^*}} \cdot {w'}_j^{-\frac{P_{\al_j}-1}{r^*}},\lb{jjNa}\\
&&[E_j(v_j')E_{j+1}(v_{j+1}')\cdots E_N(v_N)E_N(v_N)\cdots E_j(v_j)]\nn\\
&&=e_j(v_j')e_{j+1}(v_{j+1}')\cdots e_N(v_N)e_N(v_N)\cdots e_j(v_j){w'}_j^{-\frac{P_{\al_j}-1}{r^*}}\prod_{m=j+1}^{N-1}{w'}_m^{-\frac{P_{\al_m}}{r^*}}\cdot 
{w'_N}^{-\frac{P_{\al_N}-1/2}{r^*}} {w}_N^{-\frac{P_{\al_N}+1/2}{r^*}}\nn\\
&&\qquad \times 
 \prod_{m=j+1}^{N-1}{w}_m^{-\frac{P_{\al_m}}{r^*}} \cdot {w}_j^{-\frac{P_{\al_j}-1}{r^*}}.\lb{jjNb}
\ena

\noindent
iii) The $k\prec j \prec N$ case: 
\bea
&&[E_j(v_j)E_{j+1}(v_{j+1})\cdots E_N(v_N)E_N(v_N')\cdots E_j(v_j')\cdots E_k(v_k')]\nn\\  
&&=e_j(v_j)e_{j+1}(v_{j+1})\cdots e_N(v_N)e_N(v_N')\cdots e_j(v_j')\cdots e_k(v_k')\prod_{m=j}^{N-1}w_m^{-\frac{P_{\al_m}}{r^*}}
\cdot 
w_N^{-\frac{P_{\al_N}-1/2}{r^*}} {w'}_N^{-\frac{P_{\al_N}+1/2}{r^*}}\nn\\
&&\qquad \times 
 \prod_{m=k+1}^{N-1}{w'}_m^{-\frac{P_{\al_m}}{r^*}} \cdot {w'}_k^{-\frac{P_{\al_k}-1}{r^*}},\lb{kjNa}\\
&&[E_j(v_k')E_{k+1}(v_{k+1}')\cdots E_j(v_j')\cdots E_N(v'_N)E_N(v_N)\cdots E_k(v_j)]\nn\\  
&&=e_j(v_k')e_{k+1}(v_{k+1}')\cdots e_j(v_j')\cdots e_N(v'_N)e_N(v_N)\cdots e_k(v_j)\prod_{m=k}^{j-2}{w'}_m^{-\frac{P_{\al_m}}{r^*}}
\cdot {w'}_{j-1}^{-\frac{P_{\al_{j-1}}+1}{r^*}}\cdot {w'_j}^{-\frac{P_{\al_j}-1}{r^*}} 
\nn\\
&&\qquad \times \cdot \prod_{m=j+1}^{N-1}{w'}_m^{-\frac{P_{\al_m}}{r^*}}\cdot 
{w'}_N^{-\frac{P_{\al_N}-1/2}{r^*}} {w}_N^{-\frac{P_{\al_N}+1/2}{r^*}}
 \prod_{m=j+1}^{N-1}{w}_m^{-\frac{P_{\al_m}}{r^*}} \cdot {w}_j^{-\frac{P_{\al_j}-1}{r^*}}.\lb{kjNb}
\ena

\noindent
iv) The $k\prec j=N$ case: 
\bea
&&[E_N(v_N)E_N(v_N')\cdots E_k(v_k')]\nn\\  
&&\hspace{-0.5cm}=e_N(v_N)e_N(v_N')\cdots e_j(v_j')\cdots e_k(v_k')
w_N^{-\frac{P_{\al_N}-1/2}{r^*}} {w'}_N^{-\frac{P_{\al_N}+1/2}{r^*}}
 \prod_{m=k+1}^{N-1}{w'}_m^{-\frac{P_{\al_m}}{r^*}} \cdot {w'}_k^{-\frac{P_{\al_k}-1}{r^*}},\lb{NNka}\\
&&[E_k(v_k')\cdots E_N(v_N')E_N(v_N)]\nn\\  
&&=e_k(v_k')\cdots e_N(v_N')e_N(v_N) \prod_{m=k}^{N-2}{w'}_m^{-\frac{P_{\al_m}}{r^*}}
\cdot {w'}_{N-1}^{-\frac{P_{\al_{N-1}}+1}{r^*}}\cdot {w'}_N^{-\frac{P_{\al_k}-3/2}{r^*}} 
{w}_N^{-\frac{P_{\al_N}-1/2}{r^*}}.\lb{kNNb}
\ena

\noindent
v) The $j\prec k=N$ case: 
\bea
&&[E_j(v_j)\cdots E_N(v_N)E_N(v_N')]\nn\\  
&&= e_k(v_j)\cdots e_N(v_N)e_N(v_N')\prod_{m=j}^{N-2}w_m^{-\frac{P_{\al_m}}{r^*}}
\cdot w_{N-1}^{-\frac{P_{\al_{N-1}}+1}{r^*}}\cdot w_N^{-\frac{P_{\al_k}-3/2}{r^*}} 
{w'}_N^{-\frac{P_{\al_N}-1/2}{r^*}},\lb{jNNa}\\
&&[E_N(v_N')E_N(v_N)\cdots E_k(v_j)] \nn\\
&&= e_N(v_N')e_N(v_N)\cdots e_k(v_j)
{w'}_N^{-\frac{P_{\al_N}-1/2}{r^*}} {w}_N^{-\frac{P_{\al_N}+1/2}{r^*}}
 \prod_{m=j+1}^{N-1}{w}_m^{-\frac{P_{\al_m}}{r^*}} \cdot {w}_j^{-\frac{P_{\al_j}-1}{r^*}}.\lb{jNNb}
\ena

\noindent
vi) The $j=k=N$ case: 
\bea
&&[E_N(v_N)E_N(v_N')]= e_N(v_N)e_N(v_N')w_N^{-\frac{P_{\al_k}-3/2}{r^*}} 
{w'}_N^{-\frac{P_{\al_N}-1/2}{r^*}},\lb{NNa}\\
&&[E_N(v_N')E_N(v_N)]=e_N(v_N')e_N(v_N){w'}_N^{-\frac{P_{\al_k}-3/2}{r^*}} 
{w}_N^{-\frac{P_{\al_N}-1/2}{r^*}}.\lb{NNb}
\ena

The  $F_j(v_j)$'s counterpart of the product $[E_{j_1}(v_1)\cdots E_{j_m}(v_m) ]$  is obtained by replacing $e_{j_k}(v_k)$ with $f_{j_k}(v_k)$, $w_j$ with $w_j^{-1}$ and $r^*$ with $r$.

\begin{conj}\lb{firstHC}
For $1\preceq j \prec l\preceq 0$, 
\be
&&F_{j,l}^+(v)
=a_{j,l}\oint_{C_{j,l}^+} \prod_{m=j}^{l-1}\frac{dw_m}{2\pi i w_m}
[F_{l-1}(v_{l-1})
F_{l-2}(v_{l-2})\cdots F_{j}(v_j)]f^+_{j,l}(v,v_j\cdots, v_{l-1},P+h)\nn\\
&&\qquad\qquad+a_{j,l}\oint_{C_{j,l}^{-}} \prod_{m=j}^{l-1}\frac{dw_m}{2\pi i w_m}
[F_{j}(v_j)
F_{j+1}(v_{j+1})\cdots F_{l-1}(v_{l-1})]
 f^{-}_{j,l}(v,v_j\cdots, v_{l-1},P+h), \nn\\
 &&F_{-l,-j}^+(v)
=a_{-l,-j}\oint_{C^+_{-l,-j}} \prod_{m=j}^{l-1}\frac{dw_m}{2\pi i w_m}
[F_{j}(v_j)
F_{j+1}(v_{j+1})\cdots F_{l-1}(v_{l-1})]
 f^{+}_{-l,-j}(v,v_j\cdots, v_{l-1},P+h)\nn\\
&&\qquad\qquad +a_{-l,-j}\oint_{C^-_{-l,-j}} \prod_{m=j}^{l-1}\frac{dw_m}{2\pi i w_m}
[F_{l-1}(v_{l-1})
F_{l-2}(v_{l-2})\cdots F_{j}(v_j)]f^{-}_{-l,-j}(v,v_j\cdots, v_{l-1},P+h),
\nn\\
&&E_{l,j}^+(v)=
a_{l,j}^*\oint_{C_{l,j}^{*+}} \prod_{m=j}^{l-1}\frac{dw_m}{2\pi i w_m}
[E_{j}(v_j)E_{j+1}(v_{j+1})\cdots E_{l-1}(v_{l-1})]g^+_{l,j}(v,v_j\cdots, v_{l-1},P)\nn
\\
&&\qquad\qquad+ 
a_{l,j}^*\oint_{C_{l,j}^{*-}} \prod_{m=j}^{l-1}\frac{dw_m}{2\pi i w_m}
[E_{l-1}(v_{l-1})E_{l-2}(v_{l-2})\cdots E_{j}(v_j)]g^-_{l,j}(v,v_j\cdots, v_{l-1},P),\nn
\\
&&E_{-j,-l}^+(v)=
a_{-j,-l}^*\oint_{C^{*+}_{-j,-l}} \prod_{m=j}^{l-1}\frac{dw_m}{2\pi i w_m}
 [E_{l-1}(v_{l-1})E_{l-2}(v_{l-2})\cdots E_{j}(v_j)]g^+_{-j,-l}(v,v_j,\cdots, v_{l-1},P)\\
 &&\qquad\qquad+a_{l,j}^*\oint_{C_{l,j}^{*-}} \prod_{m=j}^{l-1}\frac{dw_m}{2\pi i w_m}
[E_{j}(v_j)E_{j+1}(v_{j+1})\cdots E_{l-1}(v_{l-1})]g^-_{-j,-l}(v,v_j\cdots, v_{l-1},P),
\en
\be
&&\hspace{-1cm}f^\pm_{j,l}(v,v_j\cdots, v_{l-1},P+h)= \frac{[v-v_{l-1}+(P+h)_{j,l}+
\frac{l-1}{2}-1][1]}{[v-v_{l-1}+
\frac{l-1}{2}][(P+h)_{j,l}-1]}
\prod_{m=j}^{l-2}\frac{[v_{m+1}-v_{m}+(P+h)_{j,m+1}-\frac{1}{2}][1]}{
[v_{m+1}-v_{m}\pm\frac{1}{2}][(P+h)_{j,m+1}]},\\
&&\hspace{-1cm}f^{\pm}_{-l,-j}(v,v_j\cdots, v_{l-1},P+h)\nn\\
&&\hspace{-1cm}= \frac{[v-v_{j}+(P+h)_{-l,-j}
-\frac{j}{2}-\eta-1+\delta_{l,0}][1]}{[v-v_{j}
-\frac{j}{2}-\eta][(P+h)_{-l,-j}-1+\delta_{l,0}]}
\prod_{m=j}^{l-2}\frac{[v_{m}-v_{m+1}+(P+h)_{-l,-(m+1)}-\frac{1}{2}+\delta_{l,0}][1]}{
[v_{m}-v_{m+1}\pm\frac{1}{2}][(P+h)_{-l,-(m+1)}+\delta_{l,0}]},\\
&&\hspace{-1cm}g^\pm_{l,j}(v,v_j\cdots, v_{l-1},P)=\frac{[v-v_{l-1}+\frac{l-1-c}{2}+1-
P_{j,l}]^*[1]^*}
{[v-v_{l-1}+\frac{l-1-c}{2}]^*
[P_{j,l}-1]^*}
\prod_{m=j}^{l-2}\frac{[v_{m+1}-v_{m}-P_{j,m+1}+\frac{1}{2}]^*[1]^*}{
[v_{m+1}-v_{m}\pm\frac{1}{2}]^*[P_{j,m+1}-1]^*},~~\\
&&\hspace{-1cm}g^\pm_{-j,-l}(v,v_j,\cdots, v_{l-1},P)\nn\\
&&= \frac{[v-v_{j}-\frac{j+c}{2}-\eta+1-
P_{-l,-j}-\delta_{l,0}]^*[1]^*}
{[v-v_{j}-\frac{j+c}{2}-\eta]^*
[P_{-l,-j}-1+\delta_{l,0}]^*}
\prod_{m=j}^{l-2}
\frac{[v_{m}-v_{m+1}-P_{-l,-(m+1)}+\frac{1}{2}-\delta_{l,0}]^*[1]^*}{
[v_{m}-v_{m+1}\pm\frac{1}{2}]^*[P_{-l,-(m+1)}-1+\delta_{l,0}]^*}.
\en
\be
&&C^\pm_{j,l}\ :\ |pq^lw|<|w_{l-1}|<|q^lw|,\quad |pq^{\pm1}w_{m+1}|<|w_m|<|q^{\pm1}w_{m+1}|
,\\
&&C^\pm_{-l,-j}\ :\ |pq^{-j+1}\xi w|<|w_{j}|<|q^{-j+1}\xi w|,\quad |pq^{\pm1}w_{m}|<|w_{m+1}|<|q^{\pm1}w_{m}|
,\\
&&C^{*\pm}_{l,j}\ :\ |q^{l+c}w|<|w_{l-1}|<|p^{*-1}q^{l+c}w|,\quad |q^{\pm1}w_{m+1}|<|w_m|<|p^{*-1}q^{\pm1}w_{m+1}|
,\\
&&C^{*\pm}_{-j,-l}\ :\ |q^{-j+1+c}\xi w|<|w_{j}|<|p^{*-1}q^{-j+1+c}\xi w|,\quad |q^{\pm1}w_{m}|<|w_{m+1}|<|p^{*-1}q^{\pm1}w_{m}|
\en
for $j\leq m\leq l-2$. Here $N+1\equiv 0$. The case $l\not=0 $ can be proved in the same way as for $U_{q,p}(A_N^{(1)})$\cite{KojimaKonno}. 

\end{conj}

The following is a conjectural expression for the half currents of the second type, which we 
obtained by requiring that the integrand should be single-valued and the vector represention of the 
$L$-operator should reproduce the $R$-matrix.  
\begin{conj}\lb{conjEmkj}
i) For $j \prec k \preceq N$, 
\be
&&E_{-k,j}^+(v)\\
&&=
a_{-k,j}^*\oint_{C^{*+}_{-k,j}} \prod_{m=j}^{N}\frac{dw_m}{2\pi i w_m}
\prod_{m=k}^{N}\frac{dw'_m}{2\pi i w'_m}
 [E_{j}(v_{j})E_{j+1}(v_{j+1})\cdots E_{k}(v_{k})\cdots E_{N}(v_{N})E_{N}(v'_N)\cdots E_k(v'_k)]\nn\\
 &&\qquad\qquad\qquad\qquad \times g^+_{-k,j}(v,v_j,\cdots,v_k',P)\\
 &&+a_{-k,j}^*\frac{[1]^*}{[\frac{1}{2}]^*}\oint_{{C}^{*-}_{-k,j}} \prod_{m=j}^{N}\frac{dw_m}{2\pi i w_m}
\prod_{m=k}^{N}\frac{dw'_m}{2\pi i w'_m}
[ E_{k}(v'_{k})E_{k+1}(v'_{k+1})\cdots E_{N}(v_{N}')E_{N}(v_N)\cdots E_{k}(v_{k})\cdots E_k(v_j)]\\
&&\qquad\qquad\qquad\qquad \times g^-_{-k,j}(v,v_j,\cdots,v_k',P),
\en
\be
&&g^+_{-k,j}(v,v_j,\cdots,v_k',P)\nn\\
&&\quad= \frac{[v-v'_{k}-\frac{k+c}{2}-\eta-
P_{j,-k}+1]^*[1]^*}
{[v-v'_{k}-\frac{k+c}{2}-\eta]^*
[P_{j,-k}-1]^*}\prod_{m=k+1}^{N}
\frac{[v'_{m-1}-v'_{m}-P_{j,-m}+\frac{1}{2}]^*[1]^*}{
[v'_{m-1}-v'_{m}+\frac{1}{2}]^*[P_{j,-m}-1]^*}\\
&&\quad\times 
\frac{[v'_{N}-v_{N}-P_{j,0}+\frac{1}{2}]^*[1]^*}{
[v'_{N}-v_{N}-\frac{1}{2}]^*[P_{j,0}]^*}
\prod_{m=k+1}^{N}
\frac{[v_{m}-v_{m-1}-P_{j,m}+\frac{1}{2}]^*[1]^*}{
[v_{m}-v_{m-1}+\frac{1}{2}]^*[P_{j,m}-1]^*}
\frac{[v_{k}-v_{k-1}-P_{j,k}-\frac{1}{2}]^*[1]^*}{
[v_{k}-v_{k-1}+\frac{1}{2}]^*[P_{j,k}]^*}\\
&&\quad \times \prod_{m=j+1}^{k-1}
\frac{[v_{m}-v_{m-1}-P_{j,m}+\frac{1}{2}]^*[1]^*}{
[v_{m}-v_{m-1}+\frac{1}{2}]^*[P_{j,m}-1]^*},
\\
&&g^-_{-k,j}(v,v_j,\cdots,v_k',P)\\
&&\quad =\frac{[v-v'_{k}-\frac{k+c}{2}-\eta-
P_{j,-k}+1]^*[1]^*}
{[v-v'_{k}-\frac{k+c}{2}-\eta]^*
[P_{j,-k}-1]^*}\prod_{m=k+1}^{N}
\frac{[v'_{m-1}-v'_{m}-P_{j,-m}+\frac{1}{2}]^*[1]^*}{
[v'_{m-1}-v'_{m}-\frac{1}{2}]^*[P_{j,-m}-1]^*}\\
&&\quad\times 
\frac{[v'_{N}-v_{N}]^*}{
[v'_{N}-v_{N}+\frac{1}{2}]^*}
\frac{[v'_{N}-v_{N}-P_{j,0}]^*[1]^*}{
[v'_{N}-v_{N}-\frac{3}{2}]^*[P_{j,0}-\frac{3}{2}]^*}
\prod_{m=j+1}^{N}
\frac{[v_{m}-v_{m-1}-P_{j,m}+\frac{1}{2}]^*[1]^*}{
[v_{m}-v_{m-1}-\frac{1}{2}]^*[P_{j,m}-1]^*}.
\en
\be
&&C^{*\pm}_{-k,j}\ :\ |q^{-k+1+c}\xi w|<|w_{k}'|<|p^{*-1}q^{-k+1+c}\xi w|, \nn \\
&&\qquad\qquad |q^{\pm1}w'_{m-1}|<|w'_{m}|<|p^{*-1}q^{\pm1}w'_{m-1}|\quad (k+1\leq m\leq N), \\
&&\qquad\qquad  |q^{\mp1}w'_{N}|<|w_N|<|p^{*-1}q^{\mp1}w'_{N}|,\quad 
 |q^{\pm1}w_{n}|<|w_{n-1}|<|p^{*-1}q^{\pm1}w_{n}|\quad (j+1\leq n\leq N)
\en
\noindent
ii) For $j = k \prec N$, we obtain from \eqref{jjNa} and \eqref{NNa},  
\be
&&E_{-j,j}^+(v)\\
&&=
a_{-j,j}^*\oint_{C^{*+}_{-j,j}} \prod_{m=j}^{N}\frac{dw_m}{2\pi i w_m}
\prod_{m=j}^{N}\frac{dw'_m}{2\pi i w'_m}
 [E_{j}(v_{j})E_{j+1}(v_{j+1})\cdots E_{N}(v_{N})E_{N}(v'_N)\cdots E_k(v'_j)]\\
&&\qquad\qquad\qquad\qquad \times g^+_{-j,j}(v,v_j,\cdots,v_j',P)\\ 
&&+a_{-j,j}^*\frac{[1]^*}{[\frac{1}{2}]^*}\oint_{{C}^{*-}_{-j,j}} \prod_{m=j}^{N}\frac{dw_m}{2\pi i w_m}
\prod_{m=j}^{N}\frac{dw'_m}{2\pi i w'_m}
 [E_{j}(v'_{j})E_{j+1}(v'_{j+1})\cdots E_{N}(v'_{N})E_{N}(v_N)\cdots E_k(v_j)]\\
&&\qquad\qquad\qquad\qquad \times g^-_{-j,j}(v,v_j,\cdots,v_j',P),\\  
&&g^+_{-j,j}(v,v_j,\cdots,v_j',P)\\ 
&&\quad= \frac{[v-v'_{j}-\frac{j+c}{2}-\eta-2P_{j}+1]^*[1]^*}
{[v-v'_{j}-\frac{j+c}{2}-\eta]^*
[2P_{j}-1]^*}\prod_{m=j+1}^{N}
\frac{[v'_{m-1}-v'_{m}-P_{j,-m}+\frac{1}{2}]^*[1]^*}{
[v'_{m-1}-v'_{m}+\frac{1}{2}]^*[P_{j,-m}-1]^*}\\
&&\quad\times 
\frac{[v'_{N}-v_{N}-P_{j,0}+\frac{1}{2}]^*[1]^*}{
[v'_{N}-v_{N}-\frac{1}{2}]^*[P_{j,0}]^*}
\prod_{m=j+1}^{N}
\frac{[v_{m}-v_{m-1}-P_{j,m}+\frac{1}{2}]^*[1]^*}{
[v_{m}-v_{m-1}+\frac{1}{2}]^*[P_{j,m}-1]^*}
\frac{[v-v_{j}+\frac{-j+1+c}{2}-\eta+1]^*}{
[v-v_{j}+\frac{-j+1+c}{2}-\eta]^*},\\
&&g^-_{-j,j}(v,v_j,\cdots,v_j',P)\\  
&&\quad= \frac{[v-v'_{j}-\frac{j+c}{2}-\eta-
2P_{j}+2]^*[1]^*}
{[v-v'_{j}-\frac{j+c}{2}-\eta]^*
[2P_{j}-1]^*}\prod_{m=j+1}^{N}
\frac{[v'_{m-1}-v'_{m}-P_{j,-m}+\frac{1}{2}]^*[1]^*}{
[v'_{m-1}-v'_{m}-\frac{1}{2}]^*[P_{j,-m}-1]^*}\\
&&\quad\times 
\frac{[v'_{N}-v_{N}]^*}{
[v'_{N}-v_{N}+\frac{1}{2}]^*}
\frac{[v'_{N}-v_{N}-P_{j,0}]^*[1]^*}{
[v'_{N}-v_{N}-\frac{3}{2}]^*[P_{j,0}-\frac{3}{2}]^*}
\prod_{m=j+1}^{N}
\frac{[v_{m}-v_{m-1}-P_{j,m}+\frac{1}{2}]^*[1]^*}{
[v_{m}-v_{m-1}-\frac{1}{2}]^*[P_{j,m}-1]^*}.
\en
\be
&&C^{*\pm}_{-j,j}\ :\ |q^{-j+1+c}\xi w|<|w_{j}'|<|p^{*-1}q^{-j+1+c}\xi w|,\quad |q^{\pm1}w'_{m-1}|<|w'_{m}|<|p^{*-1}q^{\pm1}w'_{m-1}|, \\
&&\qquad\qquad |q^{\mp1}w'_{N}|<|w_N|<|p^{*-1}q^{\mp1}w'_{N}|,\quad 
 |q^{\pm1}w_{m}|<|w_{m-1}|<|p^{*-1}q^{\pm1}w_{m}|
\en
for $j+1\leq m\leq N$ with 
\be
 &&|q^{-j+1+c}\xi w|<|w_{j}|<|p^{*-1}q^{-j+1+c}\xi w| \qquad \mbox{for } \ C^{*+}_{-j,j},\\
&& |q^{-3}w'_{N}|<|w_N|<|p^{*-1}q^{-3}w'_{N}| \qquad \mbox{for } \ C^{*-}_{-j,j}.  
\en
\noindent
iii) For $k \prec j \preceq N$, we obtain from \eqref{kjNa}
\be
&&E_{-k,j}^+(v)\\
&&=
a_{-k,j}^*\oint_{C^{*+}_{-k,j}} \prod_{m=j}^{N}\frac{dw_m}{2\pi i w_m}
\prod_{m=k}^{N}\frac{dw'_m}{2\pi i w'_m}
[E_{j}(v_{j})E_{j+1}(v_{j+1})\cdots E_{N}(v_{N})E_{N}(v'_N)\cdots E_j(v'_j)\cdots 
E_k(v'_k)]\\
&&\qquad\qquad\qquad\qquad \times g^+_{-k,j}(v,v_j,\cdots,v_k',P)\\ 
&&+a_{-k,j}^*\frac{[1]^*}{[\frac{1}{2}]^*}\oint_{{C}^{*-}_{-k,j}} \prod_{m=j}^{N}\frac{dw_m}{2\pi i w_m}
\prod_{m=k}^{N}\frac{dw'_m}{2\pi i w'_m}
[ E_{k}(v'_{k})E_{k+1}(v_{k+1}')\cdots E_{j}(v'_{j})\cdots E_{N}(v'_{N})E_{N}(v_N)\cdots E_k(v_j)]\\
&&\qquad\qquad\qquad\qquad \times g^-_{-k,j}(v,v_j,\cdots,v_k',P),\\ 
&&g^+_{-k,j}(v,v_j,\cdots,v_k',P),\\ 
&&\quad= \frac{[v-v'_{k}-\frac{k+c}{2}-\eta-
P_{j,-k}+1]^*[1]^*}
{[v-v'_{k}-\frac{k+c}{2}-\eta]^*
[P_{j,-k}-1]^*}\prod_{m=k+1}^{N}
\frac{[v'_{m-1}-v'_{m}-P_{j,-m}+\frac{1}{2}]^*[1]^*}{
[v'_{m-1}-v'_{m}+\frac{1}{2}]^*[P_{j,-m}-1]^*}\\
&&\quad\times 
\frac{[v'_{N}-v_{N}-P_{j,0}+\frac{1}{2}]^*[1]^*}{
[v'_{N}-v_{N}-\frac{1}{2}]^*[P_{j,0}]^*}
\prod_{m=j+1}^{N}
\frac{[v_{m}-v_{m-1}-P_{j,m}+\frac{1}{2}]^*[1]^*}{
[v_{m}-v_{m-1}+\frac{1}{2}]^*[P_{j,m}-1]^*}, \\
&&g^-_{-k,j}(v,v_j,\cdots,v_k',P)\\ 
&&\quad= \frac{[v-v'_{k}-\frac{k+c}{2}-\eta-
P_{j,-k}+1]^*[1]^*}
{[v-v'_{k}-\frac{k+c}{2}-\eta]^*
[P_{j,-k}-1]^*}\prod_{m=k+1\atop m\not=j}^{N}
\frac{[v'_{m-1}-v'_{m}-P_{j,-m}+\frac{1}{2}]^*[1]^*}{
[v'_{m-1}-v'_{m}-\frac{1}{2}]^*[P_{j,-m}-1]^*}
\frac{[v'_{j-1}-v'_{j}-2P_{j}+\frac{3}{2}]^*[1]^*}{
[v'_{j-1}-v'_{j}-\frac{1}{2}]^*[2P_{j}-1]^*}\\
&&\quad\times 
\frac{[v'_{N}-v_{N}]^*}{
[v'_{N}-v_{N}+\frac{1}{2}]^*}
\frac{[v'_{N}-v_{N}-P_{j,0}]^*[1]^*}{
[v'_{N}-v_{N}-\frac{3}{2}]^*[P_{j,0}-\frac{3}{2}]^*}
\prod_{m=j+1}^{N}
\frac{[v_{m}-v_{m-1}-P_{j,m}+\frac{1}{2}]^*[1]^*}{
[v_{m}-v_{m-1}-\frac{1}{2}]^*[P_{j,m}-1]^*}.
\en
\be
&&C^{*\pm}_{-k,j}\ :\ |q^{-k+1+c}\xi w|<|w_{k}'|<|p^{*-1}q^{-k+1+c}\xi w|,\nn \\
&& \qquad\qquad |q^{\pm1}w'_{m-1}|<|w'_{m}|<|p^{*-1}q^{\pm1}w'_{m-1}|\quad (k+1\leq m\leq N), \\
&&\qquad\qquad  |q^{\mp1}w'_{N}|<|w_N|<|p^{*-1}q^{\mp1}w'_{N}|,\quad 
 |q^{\pm1}w_{n}|<|w_{n-1}|<|p^{*-1}q^{\pm1}w_{n}|\quad (j+1\leq n\leq N),
\en
in addition, for $C^{*-}_{-k,j}$
\be
&& |q^{-3}w'_{N}|<|w_N|<|p^{*-1}q^{-3}w'_{N}|. 
\en
vi) For $j = k = N$, we obtain from \eqref{NNa} and \eqref{NNb},  
\be
&&E_{-N,N}^+(v)=
a_{-N,N}^*\oint_{C^{*+}_{-NN}} \frac{dw_N}{2\pi i w_N}\frac{dw'_N}{2\pi i w'_N}
 [E_{N}(v_{N})E_{N}(v'_N)]g^+_{-N,N}(v,v_N,v_N',P)\\ 
&&\qquad\qquad\qquad +a_{-N,N}^*\oint_{C^{*-}_{-NN}} \frac{dw_N}{2\pi i w_N}\frac{dw'_N}{2\pi i w'_N}
[ E_{N}(v'_{N})E_{N}(v_N)]g^-_{-N,N}(v,v_N,v_N',P),\\ 
&&g^+_{-N,N}(v,v_N,v_N',P)
=\frac{[v-v'_{N}-\frac{N+c}{2}-\eta-2P_{N}+\frac{3}{2}]^*[1]^*}
{[v-v'_{N}-\frac{N+c}{2}-\eta]^*
[2P_{N}-1]^*}\\
&&\qquad\qquad\qquad\qquad\times
\frac{[v'_{N}-v_{N}-P_{N,0}+1]^*[1]^*}{
[v'_{N}-v_{N}-\frac{1}{2}]^*[P_{N,0}-\frac{3}{2}]^*}
\frac{[v-v_{N}-\frac{N+c}{2}-\eta+\frac{1}{2}]^*}{
[v-v_{N}-\frac{N+c}{2}-\eta]^*}, \\
&&g^-_{-N,N}(v,v_N,v_N',P)
=
\frac{[v-v'_{N}-\frac{N+c}{2}-\eta-2P_{N}+\frac{3}{2}]^*[1]^*}
{[v-v'_{N}-\frac{N+c}{2}-\eta]^*
[2P_{N}-1]^*}\\
&&\qquad\qquad\qquad\qquad\times 
\frac{[v'_{N}-v_{N}-P_{N,0}+1]^*[1]^*}{
[v'_{N}-v_{N}+\frac{1}{2}]^*[P_{N,0}-\frac{3}{2}]^*}
\frac{[v-v_{N}-\frac{N+c}{2}-\eta+\frac{1}{2}]^*}{
[v-v_{N}-\frac{N+c}{2}-\eta]^*}. 
\en
\be
&&C^{*\pm}_{-N,N}\ :\ |q^{-N+1+c}\xi w|<|w_N|, |w'_N|<|p^{*-1}q^{-N+1+c}\xi w|,\quad 
  |q^{\mp1}w'_{N}|<|w_N|<|p^{*-1}q^{\mp1}w'_{N}|. 
\en
\end{conj}

\begin{conj}\lb{conjEmkj}
i) For $j \prec k \preceq N$, 
\be
&&F_{k,-j}^+(v)\\
&&=
a_{k,-j}\frac{[1]}{[\frac{1}{2}]}\oint_{{C^+_{k-j}}} \prod_{m=k}^{N}\frac{dw_m}{2\pi i w_m}
\prod_{m=j}^{N}\frac{dw'_m}{2\pi i w'_m}
[ F_{k}(v_{k})F_{k+1}(v_{k+1})\cdots F_{N}(v_{N})F_{N}(v'_N)\cdots F_{k}(v'_{k})\cdots F_j(v'_j)]\nn\\
&&\qquad\qquad\qquad\qquad\times f^+_{k,-j}(v,v_k,\cdots,v_j',P+h)\\
&&\quad +a_{k,-j}\frac{[1]}{[\frac{1}{2}]}\oint_{{C^-_{k-j}}} \prod_{m=k}^{N}\frac{dw_m}{2\pi i w_m}
\prod_{m=j}^{N}\frac{dw'_m}{2\pi i w'_m}
 [F_{j}(v'_{j})F_{j+1}(v'_{j+1})\cdots F_{k}(v'_{k})\cdots F_{N}(v'_{N})F_{N}(v_N)\cdots F_k(v_k)]\nn\\
 &&\qquad\qquad\qquad\qquad\times f^-_{k,-j}(v,v_k,\cdots,v_j',P+h),\\
 &&f^+_{k,-j}(v,v_k,\cdots,v_j',P+h)\nn\\
&&\qquad= \frac{[v-v'_{j}-\frac{j}{2}-\eta-1+(P+h)_{k,-j}][1]}
{[v-v'_{j}-\frac{j}{2}-\eta]
[P_{k,-j}-1]}\prod_{m=j+1}^{N}
\frac{[v'_{m-1}-v'_{m}+(P+h)_{k,-m}-\frac{1}{2}][1]}{
[v'_{m-1}-v'_{m}-\frac{1}{2}][(P+h)_{k,-m}]}\\
&&\qquad\times 
\frac{[v'_{N}-v_{N}]}{
[v'_{N}-v_{N}+\frac{1}{2}]}
\frac{[v'_{N}-v_{N}+(P+h)_{k,0}-2][1]}{
[v'_{N}-v_{N}-\frac{3}{2}][(P+h)_{k,0}-\frac{1}{2}]}
\prod_{m=k+1}^{N}
\frac{[v_{m}-v_{m-1}+(P+h)_{k,m}-\frac{1}{2}][1]}{
[v_{m}-v_{m-1}-\frac{1}{2}][(P+h)_{k,m}]},\\
&&f^-_{k,-j}(v,v_k,\cdots,v_j',P+h)\nn\\
&&\qquad\times \frac{[v-v'_{j}-\frac{j}{2}-\eta-1+(P+h)_{k,-j}][1]}
{[v-v'_{j}-\frac{j}{2}-\eta]
[P_{k,-j}-1]}\prod_{m=j+1\atop \not=k}^{N}
\frac{[v'_{m-1}-v'_{m}+(P+h)_{k,-m}-\frac{1}{2}][1]}{
[v'_{m-1}-v'_{m}+\frac{1}{2}][(P+h)_{k,-m}]}\\
&&\qquad\times 
\frac{[v'_{k-1}-v'_{k}+2(P+h)_{k}-\frac{3}{2}][1]}{
[v'_{k-1}-v'_{k}+\frac{1}{2}][2(P+h)_{k}-1]}
\frac{[v'_{N}-v_{N}+(P+h)_{k,0}-\frac{3}{2}][1]}{
[v'_{N}-v_{N}-\frac{1}{2}][(P+h)_{k,0}-{2}]} \nn \\
&&\qquad \times 
\prod_{m=k+1}^{N}
\frac{[v_{m}-v_{m-1}+(P+h)_{k,m}-\frac{1}{2}][1]}{
[v_{m}-v_{m-1}+\frac{1}{2}][(P+h)_{k,m}]}.
\en
\be
&&C^{\pm}_{k,-j}\ :\ |pq^{-j+1}\xi w|<|w_{j}'|<|q^{-j+1}\xi w|,\quad |pq^{\mp1}w'_{m-1}|<|w'_{m}|<|q^{\mp1}w'_{m-1}|\quad (j+1\leq m\leq N), \\
&&\qquad\qquad  |pq^{\pm1}w'_{N}|<|w_N|<|q^{\pm1}w'_{N}|,\quad 
 |pq^{\mp1}w_{n}|<|w_{n-1}|<|q^{\mp1}w_{n}|\quad (k+1\leq n\leq N),
\en
in addition, for $C^{+}_{k,-j}$
\be
&& |pq^{-3}w'_{N}|<|w_N|<|q^{-3}w'_{N}|. 
\en

\noindent
ii) For $j = k \prec N$, we obtain from \eqref{jjNa} and \eqref{NNa},  
\be
&&F_{j,-j}^+(v)\\
&&=
a_{j,-j}\frac{[1]}{[\frac{1}{2}]}\oint_{C^+_{j-j}} \prod_{m=j}^{N}\frac{dw_m}{2\pi i w_m}
\prod_{m=j}^{N}\frac{dw'_m}{2\pi i w'_m}
[ F_{j}(v_{j})F_{j+1}(v_{j+1})\cdots F_{N}(v_{N})F_{N}(v'_N)\cdots F_j(v'_j)]
\nn\\
&&\qquad\qquad\qquad\qquad\times f^+_{j,-j}(v,v_j,\cdots,v_j',P+h)\\
&&\quad+a_{j,-j}\oint_{{C^-_{j-j}}} \prod_{m=j}^{N}\frac{dw_m}{2\pi i w_m}
\prod_{m=j}^{N}\frac{dw'_m}{2\pi i w'_m}
[ F_{j}(v'_{j})F_{j+1}(v'_{j+1})\cdots F_{N}(v'_{N})F_{N}(v_N)\cdots F_j(v_j)]
\nn\\
&&\qquad\qquad\qquad\qquad\times f^-_{j,-j}(v,v_j,\cdots,v_j',P+h),\\
&&f^+_{j,-j}(v,v_j,\cdots,v_j',P+h)\nn\\
&&\qquad= \frac{[v-v'_{j}-\frac{j}{2}-\eta-1+
2(P+h)_{j}-1][1]}
{[v-v'_{j}-\frac{j}{2}-\eta]
[2(P+h)_{j}-3]}\prod_{m=j+1}^{N}
\frac{[v'_{m-1}-v'_{m}+(P+h)_{j,-m}-\frac{3}{2}][1]}{
[v'_{m-1}-v'_{m}-\frac{1}{2}][(P+h)_{j,-m}-1]}\\
&&\qquad\times 
\frac{[v_N'-v_N]}{[v_N'-v_N+\frac{1}{2}]}\frac{[v'_{N}-v_{N}+(P+h)_{j,0}-3][1]}{
[v'_{N}-v_{N}-\frac{3}{2}][(P+h)_{j,0}-\frac{3}{2}]}
\prod_{m=j+1}^{N}
\frac{[v_{m}-v_{m-1}+(P+h)_{j,m}-\frac{3}{2}][1]}{
[v_{m}-v_{m-1}-\frac{1}{2}][(P+h)_{j,m}-1]},\\
&&f^-_{j,-j}(v,v_j,\cdots,v_j',P+h)\nn\\
&&\qquad= \frac{[v-v'_{j}-\frac{j}{2}-\eta-1+
2(P+h)_{j}-2][1]}
{[v-v'_{j}-\frac{j}{2}-\eta]
[2(P+h)_{j}-3]}\prod_{m=j+1}^{N}
\frac{[v'_{m-1}-v'_{m}+(P+h)_{j,-m}-\frac{1}{2}][1]}{
[v'_{m-1}-v'_{m}+\frac{1}{2}][(P+h)_{j,-m}]}\\
&&\qquad\times 
\frac{[v'_{N}-v_{N}+(P+h)_{j,0}-\frac{3}{2}][1]}{
[v'_{N}-v_{N}-\frac{1}{2}][(P+h)_{j,0}-{2}]}
\prod_{m=j+1}^{N}
\frac{[v_{m}-v_{m-1}+(P+h)_{j,m}-\frac{1}{2}][1]}{
[v_{m}-v_{m-1}+\frac{1}{2}][(P+h)_{j,m}]}\nn\\
&&\qquad\times\frac{[v-v'_{j}-\frac{3j-2}{2}-2\eta]}{
[v-v'_{j}-\frac{3j-2}{2}-2\eta-1]}.
\en
\be
&&C^{\pm}_{j,-j}\ :\ |pq^{-j+1}\xi w|<|w_{j}'|<|q^{-j+1}\xi w|,\quad |pq^{\mp1}w'_{m-1}|<|w'_{m}|<|q^{\mp1}w'_{m-1}|, \\
&&\qquad\qquad |pq^{\pm1}w'_{N}|<|w_N|<|q^{\pm1}w'_{N}|,\quad 
 |q^{\mp1}w_{m}|<|w_{m-1}|<|q^{\mp1}w_{m}|
\en
for $j+1\leq m\leq N$ with 
\be
 && |pq^{-3}w'_{N}|<|w_N|<|q^{-3}w'_{N}| \qquad \mbox{for } \ C^{+}_{j,-j},\\
&&|pq^{-3j-5}\xi^2 w|<|w'_{j}|<|q^{-3j-5}\xi^2 w|  \qquad \mbox{for } \ C^{-}_{j,-j}.  
\en

\noindent
iii) $k \prec j \preceq N$

\be
&&F_{k,-j}^+(v)\\
&&=
a_{k,-j}\frac{[1]}{[\frac{1}{2}]}\oint_{{C^+_{k,-j}}} \prod_{m=k}^{N}\frac{dw_m}{2\pi i w_m}
\prod_{m=j}^{N}\frac{dw'_m}{2\pi i w'_m}
[ F_{k}(v_{k})F_{k+1}(v_{k+1})\cdots F_j(v_j)\cdots  F_{N}(v_{N})F_{N}(v'_N)\cdots F_j(v'_j)]
\nn\\
&&\qquad\qquad\qquad\qquad\times f^+_{k,-j}(v,v_k,\cdots,v_j',P+h)\\
&&\quad+a_{k,-j}\frac{[1]}{[\frac{1}{2}]}\oint_{{C^-_{k,-j}}} \prod_{m=k}^{N}\frac{dw_m}{2\pi i w_m}
\prod_{m=j}^{N}\frac{dw'_m}{2\pi i w'_m}
 [F_{j}(v'_{j})F_{j+1}(v'_{j+1})\cdots  F_{N}(v'_{N})F_{N}(v_N)\cdots F_{j}(v_{j})\cdots F_k(v_k)]
\nn\\
&&\qquad\qquad\qquad\qquad\times f^-_{k,-j}(v,v_k,\cdots,v_j',P+h),\\
&&f^+_{k,-j}(v,v_k,\cdots,v_j',P+h)\\
&&\qquad=\frac{[v-v'_{j}-\frac{j}{2}-\eta-1+(P+h)_{k,-j}][1]}
{[v-v'_{j}-\frac{j}{2}-\eta]
[P_{k,-j}-1]}\prod_{m=j+1}^{N}
\frac{[v'_{m-1}-v'_{m}+(P+h)_{k,-m}-\frac{1}{2}][1]}{
[v'_{m-1}-v'_{m}-\frac{1}{2}][(P+h)_{k,-m}]}\\
&&\qquad\times 
\frac{[v'_{N}-v_{N}]}{
[v'_{N}-v_{N}+\frac{1}{2}]}
\frac{[v'_{N}-v_{N}+(P+h)_{k,0}-2][1]}{
[v'_{N}-v_{N}-\frac{3}{2}][(P+h)_{k,0}-\frac{1}{2}]}
\prod_{m=k+1\atop \not=j}^{N}
\frac{[v_{m}-v_{m-1}+(P+h)_{k,m}-\frac{1}{2}][1]}{
[v_{m}-v_{m-1}-\frac{1}{2}][(P+h)_{k,m}]}\nn\\
&&\qquad\times  \frac{[v_{j}-v_{j-1}+(P+h)_{k,j}+\frac{1}{2}][1]}{
[v_{j}-v_{j-1}-\frac{1}{2}][(P+h)_{k,j}]},\\ 
&&f^-_{k,-j}(v,v_k,\cdots,v_j',P+h)\nn\\
&&\qquad= \frac{[v-v'_{j}-\frac{j}{2}-\eta-1+(P+h)_{k,-j}][1]}
{[v-v'_{j}-\frac{j}{2}-\eta]
[P_{k,-j}-1]}\prod_{m=j+1}^{N}
\frac{[v'_{m-1}-v'_{m}+(P+h)_{k,-m}-\frac{1}{2}][1]}{
[v'_{m-1}-v'_{m}+\frac{1}{2}][(P+h)_{k,-m}]}\\
&&\qquad\times 
\frac{[v'_{N}-v_{N}+(P+h)_{k,0}-\frac{3}{2}][1]}{
[v'_{N}-v_{N}-\frac{1}{2}][(P+h)_{k,0}-{2}]}
\prod_{m=k+1}^{N}
\frac{[v_{m}-v_{m-1}+(P+h)_{k,m}-\frac{1}{2}][1]}{
[v_{m}-v_{m-1}+\frac{1}{2}][(P+h)_{k,m}]}.
\en
\be
&&C^{\pm}_{k,-j}\ :\ |pq^{-j+1}\xi w|<|w_{j}'|<|q^{-j+1}\xi w|,\quad |pq^{\mp1}w'_{m-1}|<|w'_{m}|<|q^{\mp1}w'_{m-1}|\quad (j+1\leq m\leq N), \\
&&\qquad\qquad  |pq^{\pm1}w'_{N}|<|w_N|<|q^{\pm1}w'_{N}|,\quad 
 |pq^{\mp1}w_{n}|<|w_{n-1}|<|q^{\mp1}w_{n}|\quad (k+1\leq n\leq N),
\en
in addition, for $C^{+}_{k,-j}$
\be
&& |pq^{-3}w'_{N}|<|w_N|<|q^{-3}w'_{N}|. 
\en

vi) $j = k = N$  
\be
&&F_{N,-N}^+(v)
=
a_{N,-N}\oint_{C^+_{N-N}} \frac{dw_N}{2\pi i w_N}\frac{dw'_N}{2\pi i w'_N}
[ F_{N}(v'_{N})F_{N}(v_N)] f^+_{N,-N}(v,v_N,v_N',P+h)\\
&&\qquad\qquad\quad+a_{N,-N}\oint_{C^-_{N-N}} \frac{dw_N}{2\pi i w_N}\frac{dw'_N}{2\pi i w'_N}
[ F_{N}(v_{N})F_{N}(v'_N)]f^-_{N,-N}(v,v_N,v_N',P+h),\\
&&f^+_{N,-N}(v,v_N,v_N',P+h)=
\frac{[v-v'_{N}-\frac{N}{2}-\eta+2(P+h)_{N}-\frac{5}{2}][1]}
{[v-v'_{N}-\frac{N}{2}-\eta]
[2(P+h)_{N}-3]}\\
&&\qquad\qquad\qquad\qquad\quad\times 
\frac{[v'_{N}-v_{N}+(P+h)_{N,0}-2][1]}{
[v'_{N}-v_{N}-\frac{1}{2}][(P+h)_{N,0}-\frac{3}{2}]}
\frac{[v-v_{N}-\frac{N-1}{2}-\eta+\frac{1}{2}]}{
[v-v_{N}-\frac{N-1}{2}-\eta]}, \\
&&f^-_{N,-N}(v,v_N,v_N',P+h)=
\frac{[v-v'_{N}-\frac{N}{2}-\eta+2(P+h)_{N}-\frac{5}{2}][1]}
{[v-v'_{N}-\frac{N}{2}-\eta]
[2(P+h)_{N}-3]}\\
&&\qquad\qquad\qquad\qquad\quad\times 
\frac{[v'_{N}-v_{N}+(P+h)_{N,0}-2][1]}{
[v'_{N}-v_{N}+\frac{1}{2}][(P+h)_{N,0}-\frac{3}{2}]}
\frac{[v-v_{N}-\frac{N-1}{2}-\eta+\frac{1}{2}]}{
[v-v_{N}-\frac{N-1}{2}-\eta]}. 
\en
\be
&&C^{\pm}_{N,-N}\ :\ |pq^{-N+1}\xi w|<|w_N|, |w'_N|<|q^{-N+1}\xi w|,\quad 
  |pq^{\mp1}w'_{N}|<|w_N|<|q^{\mp1}w'_{N}|. 
\en
\end{conj}

\end{appendix}

\renewcommand{\baselinestretch}{0.7}

\end{document}